\documentclass[final,abstracton]{scrartcl}
\usepackage{graphicx}
\usepackage{mathptmx}
\usepackage{subfigure}
\usepackage{amsfonts,amssymb,amsbsy,amsmath}
\usepackage{bm}
\graphicspath{{./graphics/}{./}{./graphics/VOF}{./graphics/Intro}}
\usepackage{psfrag}
\usepackage{microtype}
\usepackage{aicescover}
\usepackage[numbers]{natbib}
\aicescovertitle{Deforming fluid domains within the finite element method \\ \Large Five mesh-based tracking methods in comparison}
\aicescoverauthor{S. Elgeti \and H. Sauerland}

\begin{document}

\aicescoverpage

\title{Deforming fluid domains within the finite element method}
\subtitle{Five mesh-based tracking methods in comparison}

\author{S. Elgeti \footnote{elgeti@cats.rwth-aachen.de} \and H. Sauerland} 

\publishers{\small Chair for Computational Analysis of Technical Systems, CCES, RWTH Aachen University, Germany}

\date{}

\maketitle

 \hrule

\begin{abstract}
Fluid flow applications can involve a number of coupled problems. One is the simulation of free-surface flows, which require the solution of a free-boundary problem. Within this problem, the governing equations of fluid flow are coupled with a domain deformation approach. This work reviews five of those approaches:  interface tracking using a boundary-conforming mesh and, in the interface capturing context, the level-set method, the volume-of-fluid method, particle methods, as well as the phase-field method. The history of each method is presented in combination with the most recent developments in the field. Particularly, the topics of extended finite elements (XFEM) and NURBS-based methods, such as Isogeometric Analysis (IGA), are addressed. For illustration purposes, two applications have been chosen: two-phase flow involving drops or bubbles and sloshing tanks. The challenges of these applications, such as the geometrically correct representation of the free surface or the incorporation of surface tension forces, are discussed.\\ 

\end{abstract}

\setlength{\footnotesep}{15pt}
\renewcommand{\thefootnote}{}
\footnote{Preprint submitted to Archives of Computational Methods in Engineering, December 19, 2014 } 

\hrule

 \section{Introduction} \label{s-intro}

In the numerical analysis of fluid flow, we often encounter free-boundary value problems: apart from the flow solution in the bulk domain, the position of (a portion of) the boundary is also unknown. This boundary can either be an external boundary or an interface between subdomains. At the boundary/interface, certain boundary conditions need to be fulfilled, which specify the position of the boundary. These conditions relate the variables of the flow (velocity, pressure, stress) across the domains under consideration of external influences, such as for example surface tension.

\begin{figure*}
\label{fig:ngram}
\centering
\includegraphics[scale=0.23]{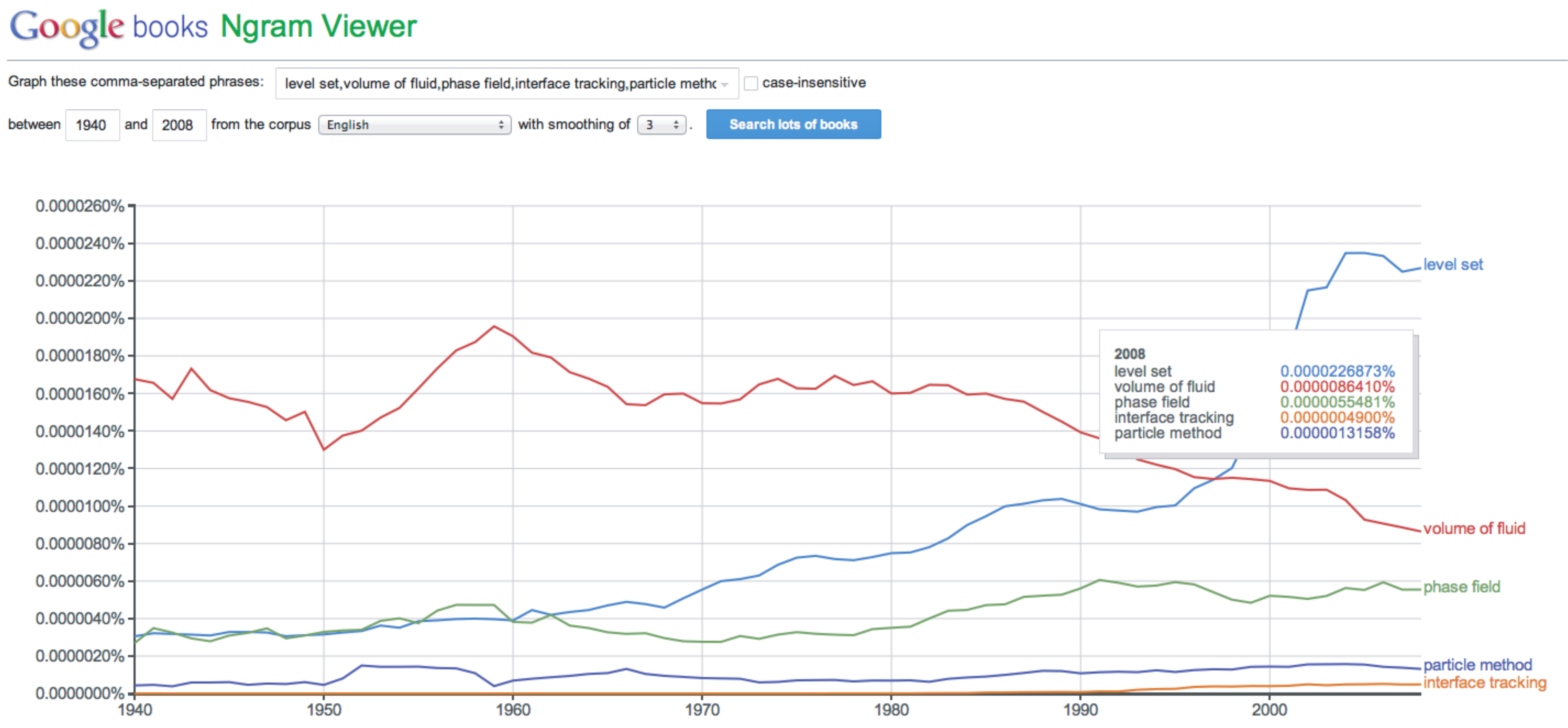}
\caption{Google Ngram Viewer: quantitative analysis of popularity of five main interface description methods.}
\end{figure*}

Numerically, in order to be able to compute the flow solution as well as the boundary/interface geometry, a measure to track the boundary starting from an initial position needs to be incorporated. Recently, Coutinho gave an overview of the most common mesh-based methods in this area \--- level-set, volume-of-fluid, and phase-field \--- from a very different point of view: a pure quantitative analysis of popularity using Google books Ngram Viewer \cite{Coutinho2014}. Ngram Viewer is a tool, which analyzes all 30 million books digitalized by Google with respect to certain keywords \cite{Michel2011}. The return value is the fraction of these books, in which the specific keyword appears. For this paper, we have added particle methods and interface tracking to the search, resulting in Figure~\ref{fig:ngram}. From the diagram, we get the notion that some methods already have a very long history of unremitting usage, such as the phase-field method, which \--- as we will later see \--- dates back already to 1871. The use of the level-set method rocketed within the last two decades, no doubt due to the large number of very effective recent developments. From pure numbers, it seems as if the volume-of-fluid methods may have already passed their zenith. Particle methods and interface tracking have so far not been able to compete with the more established methods on the grounds of popularity, but seem to be reserved for niche applications. Notwithstanding these crude quantifications, this paper aims to highlight the advantages and disadvantages of the above methods, illuminating the history of each method, but focusing on the most recent advancements in the individual fields. 

Numerous fluid flow applications involve deformable domains: drops and bubbles, die swell, dam break, liquid storage tanks, dendritic growth, spinodal decomposition, up to and including the topic of topology optimization, which can rely on the same boundary tracking methods. 

\begin{figure}
\psfrag{a}[rb][lt]{${\Omega}_1$}
\psfrag{b}[rb][lt]{${\Omega}_2$}
\psfrag{c}[cc][r]{${\Omega}_2$}
\psfrag{n}{${\bf n}$}
\psfrag{t }{${\bf t}$}
\centering
\includegraphics[scale=1.0]{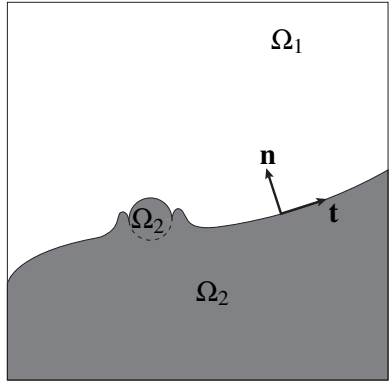}
\caption{Illustration of the general two-phase flow scenario.}
\label{two-phase-flow}
\end{figure}

The scenario we consider here are flow solutions with either one or with two immiscible fluids. The general case is illustrated in Figure~\ref{two-phase-flow}. Consider a $d$-dimensional computational domain $\Omega\subset \mathbb{R}^d$ with boundary $\Gamma = \partial \Omega$. This domain contains one, two or more immiscible fluids, which are enclosed in subdomains $\Omega_i(t)$. Position, number and shape of the individual subdomains may vary over time, i.e., both the domains themselves and the flow field are part of the solution. The main task in this context is to account for the interface $\Gamma^{int}$, which separates the distinct fluid domains and is generally in motion. Apart from its exact position, the computation of geometrical quantities of $\Gamma^{int}$, as for example its normal and tangential vectors ${\bf n}$ and ${\bf t}$ or its curvature $\kappa$ are of interest. The corresponding computational tasks can be summarized as: (1) define the shape and location of the interface, (2) track the time advancement of the interface, and (3) set boundary conditions along the interface.

The paper is organized as follows. In Section~\ref{s-governing}, the governing equations for the fluid flow with the appropriate boundary conditions for the interface are introduced. Particular challenges that come with the numerical treatment of these equations are highlighted in Section~\ref{s-challenges}. The subsequent sections concentrate on the five numerical methods considered within the scope of this paper: particle methods in Section~\ref{s-particlemethod}, the volume-of-fluid method in Section~\ref{s-volumeoffluid}, the level-set method in Section~\ref{s-levelset}, the phase-field method in Section~\ref{s-phasefield}, and mesh-conforming interface tracking in Section~\ref{s-interfacetracking}. To illustrate the methods, two common applications, drops and sloshing tanks, are the topics of Sections~\ref{s-drops}  and \ref{s-tank}.
 \section{Governing equations} \label{s-governing} 
 
In general, the governing equations for incompressible fluid flow are the instationary, incompressible Navier-Stokes equations. Consider an instationary fluid flow problem with any number (in this paper, one or two) of immiscible Newtonian phases. The computational domain at each instant in time, denoted by \(\Omega_t \), is a subset of \(\mathbb{R}^{nsd}\), where \(nsd\) is the number of space dimensions. Then at each point in time \(t \in [0,T] \), the velocity, \( {\bf u}({\bf x},t) \), and the pressure, \( p({\bf x},t) \), in each phase are governed by the following equations:

\begin{align}
\label{NavierStokes1}
\rho_i \Big( \frac{\partial {\bf u}}{\partial t} + {\bf u \cdot \nabla u} - {\bf f} \Big) -  \nabla \cdot \bm{\sigma}_i &= {\bf 0} \ \mbox{on} \ (\Omega_t)_i \ \forall \  t \in [0,T] \,, \\
\label{NavierStokes2}
{\bf \nabla \cdot u } &= 0 \ \mbox{on} \  (\Omega_t)_i  \ \forall t \in [0,T] \,,
\end{align}
for $ i=1, \dots,  np$ (number of phases) and \( \rho_i\) as the density of the respective fluid. The stress tensor \( \bm{\sigma}_i\) is defined as 

\begin{align}
\bm{\sigma}_i ( {\bf u}, p) &= -p {\bf I} + 2 \mu_i \bm{\varepsilon}( {\bf u}) \quad \mbox{on} \quad (\Omega_t)_i \,, 
\label{eqn:Newtonian}
\end{align}
with
\begin{align} \label{eqn:strain}
\bm{\varepsilon}({\bf u}) &= \frac{1}{2} (\nabla {\bf u} + (\nabla {\bf u})^T) \,,
\end{align}
where \( \mu_i  \) denotes the dynamic viscosity.  $ {\bf f} $ includes all external body forces referred to the unit mass of fluid. The computational domain $ \Omega_t $ is divided into parts $(\Omega_t)_i$, each occupied by fluid $i$. Note that the spatial domain is time-dependent, which is indicated by subscript \(t\). A subdomain may never be occupied by more than one fluid. The interface of two subdomains, \(\partial \Omega_1 \cap \partial \Omega_2\), is denoted by \( \Gamma_t^{int} \). \\

To allow a better grasp of the Navier-Stokes equations \eqref{NavierStokes1} -- \eqref{NavierStokes2}, we will illustrate the purpose of the individual terms, following in large parts the description in \cite{Cline}. 

\begin{itemize}
\item \textbf{Temporal derivative} $\frac{ \partial {\bf u}}{\partial t}$ {\bf :} The change in velocity with respect to time. This change over time is governed by the following influence factors:
\item \textbf{Inertia term} ${\bf u \cdot \nabla u}$ {\bf :} This term is a convection term arising from the conservation of momentum. The momentum of each portion of fluid needs to be conserved. Therefore, it needs to move with the fluid \--- it is convected with the fluid. 
\item \textbf{Pressure term} $-\nabla p$  {\bf :} This term appears when the Newtonian constitutive equation \eqref{eqn:Newtonian} is inserted into Equation~\eqref{NavierStokes1}. It includes forces resulting from pressure differences within the fluid into the formulation. Especially in the case of incompressible fluids, it is very important to consolidate this term with Equation~\eqref{NavierStokes2}.
\item \textbf{Friction term} $\mu_i \nabla^2 {\bf u}$ {\bf :} Again, this term arises when the Newtonian constitutive equation \eqref{eqn:Newtonian} is inserted into Equation~\eqref{NavierStokes1}. We obtain a diffusion operator equalizing the velocity of neighbouring elements. The more viscous the fluid, the stronger is the friction between neighbouring particles and thus the equalizing effect. 
\item \textbf{External forces} {\bf f} {\bf :}  This term includes external body forces such as gravity or vibrational excitation. 
\end{itemize}

For creeping flows (i.e., Reynolds number $\ll1$), the advective term in the Navier-Stokes equations is often neglected, giving rise to the Stokes equations. If furthermore the solution remains unaltered over time, the stationary Stokes equations can be employed:

\begin{align}
\label{StatStokes_1}
	-\nabla \cdot  \bm{\sigma}_i &= {\bf f} \quad \mbox{on} \quad \Omega_i \,, \\
	\nabla \cdot {\bf u} &= 0 \quad \mbox{on} \quad \Omega_i \,.
	\label{StatStokes_2}
\end{align} 
The constitutive equations for Newtonian fluids remain the same as above.

\subsection{Boundary, initial, and interface conditions} \label{sec:bc}

In the transient case, a divergence-free velocity field for the whole computational domain is needed as an initial condition:

\begin{align}
{\bf u} ( {\bf x}, 0 ) = \hat{{\bf u}}^0( {\bf x}) \quad \mbox{in}\ \Omega \ \mbox{at} \  t=0 \,.
\end{align} 

In order to obtain a well-posed system, boundary  conditions have to be imposed on the external boundary of $\Omega$, denoted as $\Gamma$. Here, we distinguish between Dirichlet and Neumann boundary conditions given by:

\begin{align}
{\bf u} &= \hat{{\bf u}} \quad \mbox{on} \ \Gamma_u,\ t \in [0,T] \,, \\
{\bf n} \cdot \bm{\sigma} &= \hat{{\bf h}} \quad \mbox{on} \ \Gamma_h,\ t \in [0,T] \label{eqn:Neumanncondition} \,,
\end{align} 

\noindent where $\hat{{\bf u}}$ and $\hat{{\bf h}}$ are prescribed velocity and stress values. $ \Gamma_u $ and $ \Gamma_h $ denote the Dirichlet and Neumann part of the boundary, forming a complementary subset of $\Gamma$, i.e., $ \Gamma_u \cup \Gamma_h = \Gamma$ and $ \Gamma_u \cap \Gamma_h = \emptyset$. {\bf n} refers to the outer normal vector on $ \Gamma_h $. 

In the case of two-phase flow, at the interface between the two phases, we impose interface conditions, which couple the velocity and stress between the two domains. The first interface condition implies that the velocities are continuous across the interface: 

\begin{align}\label{eq:contvel}
{\bf u}_1 = {\bf u}_2. 
\end{align}

The second interface condition is based on the Laplace-Young equation \cite{Batchelor1967}: 

\begin{align}\label{eq:SurfTens}
(\bm{\sigma}_2  - \bm{\sigma}_1  ) {\bf n}_1 &= \gamma \kappa {\bf n}_1\, .
\end{align}
Here, \( {\bf n}_1 \) is the, with respect to $\Omega_1$, outward unit normal vector on \( \Gamma_t^{int} \), \( \gamma \) the surface tension coefficient and \( \kappa \) the sum of the principal curvatures of \( \Gamma_t^{int} \). 

In addition to the interface conditions, the jump in density and viscosity between the fields leads to non-smooth behavior of the field quantities. If we examine condition \eqref{eq:contvel} closer, we note that it states that the normal component of the velocity $u_n = {\bf u} \cdot {\bf n}_1$ as well as the tangential velocity component $u_{\hat{t}} = {\bf u} \cdot {\bf t}_1$ must be continuous across the interface. Unaffected by this, it can however be shown that the normal gradient $\frac{\partial}{\partial n}$ of the tangential velocity is discontinuous \cite{Li2001}:

\begin{align}
\frac{\partial u_{{\hat{t}},1}}{\partial n} - \frac{\partial u_{{\hat{t}},2}}{\partial n} = - (\mu_1 - \mu_2) \frac{\partial u_{n}}{\partial {\hat{t}}} \,,
\end{align}

\noindent with $ \frac{\partial}{\partial {\hat{t}}} $ denoting the tangential derivative. In consequence, the velocity has a kink across the interface if $\mu_1 \neq \mu_2$.

Inserting the definitions for stress and strain (Equations~\eqref{eqn:Newtonian} and \eqref{eqn:strain}) into interface condition \eqref{eq:SurfTens} results in:

\begin{align}
[ -p_1 {\bf I} + \mu_1  (\nabla {\bf u}_1 + (\nabla {\bf u}_1)^T) +  p_2 {\bf I} -  \mu_2  (\nabla {\bf u}_2 + (\nabla {\bf u}_2)^T) ] {\bf n}_1 = \gamma \kappa {\bf n}_1 \,.
\label{eqn:pressurejump1}
\end{align}

If we restrict ourselves to only the normal component of relation \eqref{eqn:pressurejump1}, an expression for the pressure jump across the interface can be derived:

\begin{align}
{\bf n}_1^T [ -p_1 {\bf I} + \mu_1  (\nabla {\bf u}_1 + (\nabla {\bf u}_1)^T) +  p_2 {\bf I} -  \mu_2  (\nabla {\bf u}_2 + (\nabla {\bf u}_2)^T) ] \cdot {\bf n}_1 = \gamma \kappa \,, 
\end{align}
\begin{align}
\Leftrightarrow -p_1 + p_2 + [2\mu_1 (\nabla {\bf u}_1 {\bf n}_1) \cdot {\bf n}_1 - 2\mu_2 (\nabla {\bf u}_2 {\bf n}_1) \cdot {\bf n}_1] = \gamma \kappa \,,
\end{align}
\begin{align}
\Leftrightarrow -p_1 + p_2 + [2\mu_1 \frac{\partial u_{n,1}}{\partial n}  - 2\mu_2 \frac{\partial u_{n,2}}{\partial n}  ] = \gamma \kappa \,.
\label{eqn:pressurejump2}
\end{align}

From Equation~\eqref{eqn:pressurejump2} we can deduce that the pressure jump across the interface depends on the surface tension coefficient, the curvature of the interface, and the jump in the viscosity weighted by the normal derivative of the normal velocity component \cite{Axelsson1994,Idelsohn2009,Kang2000}. In particular, this has the consequence that pressure jumps across the interface are possible even when no surface tension effects are present, a possibility often disregarded in the modelling process. Usually however, this effect is indeed negligible as either the difference in viscosity or the change in the normal velocity are not high enough to actually make an impact \cite{Idelsohn2009}. 

A further influence on the pressure distribution is the volume force per unit mass of fluid ${\bf f}$. To see this, we examine the hydrostatic case, i.e., ${\bf u} = {\bf 0}$, where the momentum equation~\eqref{NavierStokes1} reduces to 

\begin{align}
{\boldsymbol{\nabla}} p_i = \rho_i {\bf f}_i \quad \mbox{in} \ \Omega.
\end{align}

Focusing on the interface $\Gamma^{int}$, a jump condition for the pressure \emph{gradient} is obtained:

\begin{align}
{\boldsymbol{\nabla}} p_1 - {\boldsymbol{\nabla}} p_2 = \rho_1 {\bf f}_1 - \rho_2 {\bf f}_2 = (\rho_1 - \rho_2) {\bf f}.
\end{align}

The last simplification can be made since the volume force is usually constant across the entire domain. We note that we obtain a jump in the pressure gradient proportional to the jump in density, implying a kink in the pressure distribution. This behavior is extendible to ${\bf u} \neq {\bf 0}$.

\begin{figure}[htbp]
\center
\includegraphics[width=5.0cm]{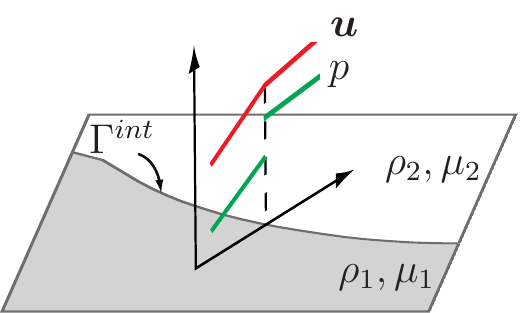}
\caption{General situation in two-phase flow: A kink in the velocity (due to $\mu_1 \neq \mu_2$) and a jump (due to surface tension)  as well as a  kink (due to volume forces) in the pressure across the interface.}
\label{scenario}
\end{figure}

All effects combined (depicted in Figure~\ref{scenario}), in the general two-phase flow scenario we have to account for a kink in the velocity as well as a jump and kink in the pressure at the interface.

\subsection{Variational form}

In the finite element method, we work with the variational form of the governing equations. One important measure when deriving the variational form is the integration by parts with the goal of reducing the differentiability requirements on the trial functions (i.e., the solution). Integration by parts furthermore naturally includes Neumann type boundary conditions as in our case Equations~\eqref{eqn:Neumanncondition} and \eqref{eq:SurfTens}. \cite{Hughes2000} 

The variational form of Equations~\eqref{NavierStokes1}--(\ref{NavierStokes2}) can be expressed as follows:  Find \( {\bf u} \) and \(p  \) such that \( \forall {\bf w}  \),  \( \forall  q  \):
\begin{align}
\begin{split}
\sum_{i=1}^{np} \bigg[ \int_{(\Omega_t)_i} {\bf w} \cdot \rho_i \left( \frac{\partial {\bf u}}{\partial t} + {{\bf u} \cdot {\bf \nabla u}} - {\bf f} \right)\;d{\bf x} + \\ \int_{(\Omega_t)_i} \mbox{\boldmath$\varepsilon$}({\bf w}):\mbox{\boldmath$\sigma$}_i (p,{\bf u})\;d{\bf x} 
+ \int_{(\Omega_t)_i} q\nabla \cdot {\bf u}^h\;d{\bf x} \bigg] \\
=  \int_{\Gamma_{h}} {\bf w} \cdot \hat{{\bf h}}\;d{\bf x} -  \gamma \int_{\Gamma^{int}_t} \kappa{\bf w} \cdot {\bf n}_1\;d{\bf x} .
\end{split}
\label{eqn:weakform}
\end{align}

Sections \ref{s-particlemethod} -- \ref{s-interfacetracking} will introduce a variety of methods that can be employed to solve the above equation numerically.
 \section{The challenges in free-surface flow} \label{s-challenges} 
 In this section, we will analyze particular challenges of Equation~\eqref{eqn:weakform}: the free boundary problem it represents, the treatment of the surface tension term, as well as the solution of additional equations on the interface. 
 
\subsection{The free boundary problem} \label{sec:ITIC}

Equations~\eqref{NavierStokes1} -- \eqref{NavierStokes2} describe a scenario where a partial differential equation is solved for an unknown function \--- in this case the fluid velocity ${\bf u}({\bf x}, t)$ and the pressure $p({\bf x}, t)$ \--- but at the same time the exact extent of the computational domain (or portions of it) is also unknown. We are dealing with a free boundary problem. Consequently, one of the computational tasks will be to monitor the domain shape \--- and especially the interface/boundary \--- throughout the simulation. Broadly, two significant approaches can be distinguished: interface capturing and interface tracking. The main difference can be located in the viewpoint, which can be either Eulerian or Lagrangian, as we will see in the following two sections. 
 
\subsubsection{Interface Capturing}

Interface capturing approaches are based on an Eulerian description and widely used in deformable domain problems. They define the interface implicitly on a fixed mesh. A characteristic scalar field $\phi$ is used to identify the two phases as well as the interface along the boundaries of the individual fluid domains. Depending on the different methods, this scalar field may be described for example by a discontinuous Heaviside function or a signed-distance function. In order to account for the interface motion, a standard advection equation

\begin{align}
\label{eqn:advectionequation}
\frac{\partial \phi}{\partial t} + {\bf u} \cdot {\bf \nabla} \phi = 0 \quad \mathrm{in} \ \Omega,\, t \in [0,T] \,,
\end{align}

\noindent is solved with $ {\bf u} $ as the fluid velocity. The most common representatives of this category are particle methods (Section~\ref{s-particlemethod}), the volume-of-fluid method (Section~\ref{s-volumeoffluid}), and the level-set method (Section~\ref{s-levelset}). In addition, the phase-field method (Section~\ref{s-phasefield}) shares some features with the interface capturing methods. 

A great advantage of the interface capturing approaches is that they are inherently able to account for topological changes of the interface. This allows for a much more flexible interface description than in interface tracking approaches. What needs to be considered however are the aspects of discontinuity treatment across the interface, mass conservation, and the application of boundary conditions along the interface, which remain a challenge in interface capturing. 

\subsubsection{Interface Tracking} \label{sec:interfacetracking1}
 
In interface tracking approaches, the interface is described explicitly on a boundary/interface conforming mesh. The idea is to track the position of the mesh nodes ${\bf x}$ in a Lagrangian fashion by integrating the evolution equation

\begin{align}
\frac{\partial {\bf x}}{\partial t} = {\bf u}({\bf x}, t) \,,
\label{eqn:Lagrangian}
\end{align}  
 
\noindent with $  {\bf u} $ as the fluid velocity in the domain. We can distinguish between fully Lagrangian approaches, where all mesh nodes are treated in a Lagrangian fashion, and Arbitrary Lagrangian-Eulerian (ALE) approaches, where the Lagrangian treatment is only applied to a portion of the mesh nodes, usually those along the boundary/interface (cf. Section~\ref{s-interfacetracking}).
 
Interface tracking approaches offer an accurate and computationally efficient approximation of the boundary/interface. Furthermore, the imposition of boundary conditions at the interface is simple with mesh nodes lying on the interface itself. However, firstly the mesh quality will usually degrade significantly in the course of large deformations and secondly topological changes require special treatment. In both cases, the last resort is remeshing with respect to the new interface. One consequence of remeshing is the need to project the field values from the old to the new mesh, a ceaseless cause for errors, and at the same time computationally costly especially in 3D. Consequently, it is desirable to keep the frequency of remeshing low (even to the level of no remeshing) \cite{Tezduyar92b}. 
 
\subsection{Surface tension}
 
Much of the computational accuracy and stability of two-phase flow depends on the discretization of the term responsible for the surface tension force; in the variational form Equation~\eqref{eqn:weakform} this is:
 
 \begin{align}
 \gamma \int_{\Gamma^{int}_t} \kappa{\bf w} \cdot {\bf n}\;d{\bf x},
 \label{st-term}
\end{align}   
 
\noindent with surface tension coefficient $\gamma$,  $\kappa$ as the sum over all principle curvatures, $ {\bf w} $ denoting the test function  and ${\bf n}$ the normal vector. For most numerical approaches \--- although we will see some exceptions later \---  the evaluation of the curvature $\kappa$ is nonviable as it contains second derivatives. As an alternative, two main strategies have evolved to avoid the direct evaluation of Equation \eqref{st-term}: the Laplace-Beltrami technique and the Continuum Surface Force approach. These approaches are applicable irrespective of the interface discretization method. They will be detailed in the following.
 
 \subsubsection{The Laplace-Beltrami technique} 
 \label{sec:Laplace-Beltrami}

The idea of the Laplace-Beltrami technique goes back to Dziuk \cite{Dziuk1991} and has been closely described in \cite{Baensch2001,Ganesan_2007,Gross2010}. The method replaces the above surface integral \eqref{st-term} requiring the curvature (and therefore the existence of second derivatives) by an integral requiring only gradient computations. In the following description, we will show two derivations of the method. The first is relevant when after discretization, an interface-conforming mesh (explicit interface representation) will be used and the second is suitable for implicit interface descriptions.

When the interface will eventually be resolved  explicitly  by the mesh, our starting point is the following relation regarding the curvature vector ${\bm{\kappa}}$ of the interface \cite{Enschenburg2007}: 

\begin{align}
    {\bm{\kappa}} &= \Delta^g {\bf X} \,, \\
    {\bf X}&: D \subset \mathbb{R}^{nsd-1}  \mapsto \mathbb{R}^{nsd} \,,
\end{align}

\noindent where $ {\bf X} $ is an immersion, i.e., a mapping from an arbitrary parameter space to a boundary or interface. An immersion can be understood as a differentiable map between differentiable manifolds whose derivative has full rank. In the case of free-surface flows, we deal with smooth embeddings, which are injective immersions. The curvature vector $\boldsymbol{\kappa}$ always points in normal direction: It can therefore be written as ($\kappa {\bf n}$). We obtain:

\begin{align}
    \kappa  {\bf n} = \Delta^g {\bf X}.
    \label{eqn:kappa-explicit}
\end{align}

This means that we have arrived at the expression needed in term \eqref{st-term}. Again in the case of smooth embeddings, $ \Delta^g $ takes on the form of the so-called Laplace-Beltrami operator, a generalization of the Laplace operator to manifolds, whose exact definition will be given later in Equation \eqref{eqn:Laplace-Beltrami-expl}. Before doing so, we define a quantity essential for any geometrical computation on our immersion ${\bf X}$: the metric $ {\bf g} $ given as 

\begin{align}
      {\bf g} = (g_{ij}) = J{\bf X}^TJ{\bf X} \quad \in \mathbb{R}^{(nsd-1)\times(nsd-1)}.
\end{align}

The metric is intimately related to the measurement of distances and angles, and as such, also an important requisite when integrals are to be defined in such a way that they are invariant to coordinate transformations. For an arbitrary scalar function $f$, the coordinate transformation yields: 

\begin{align}
      \int_{\Gamma^{int}} f({\bf x}) d{\bf x} \stackrel{def.}{=} \int_{\Gamma^{int}_\xi} f({\bf X}({\boldsymbol \xi})) \sqrt{\det g({\boldsymbol \xi})} d{\boldsymbol \xi}.
      \label{eqn:integration-on-reference-element}
\end{align}

Note that in the cases, where the coordinate transformation maps $\mathbb{R}^{nsd} \mapsto \mathbb{R}^{nsd}$,  $ J {\bf X} $ is square. As a consequence, $ \det ({\bf g}) = (\det(J {\bf X}))^2 $ and the definition includes the well-known transformation (change of variables) rule used routinely within the finite element method. 

Within the Laplace-Beltrami technique, the general procedure in modifying the integral in \eqref{st-term} will now be to replace $ \kappa {\bf n}  $ by the Laplace-Beltrami operator applied to the immersion ${\bf X}$, integrate by parts and then transform the integral to reference coordinates. The  definition used here for the Laplace-Beltrami operator can be found in \cite{Enschenburg2007}. The Laplace-Beltrami operator applied to a scalar function $ f $ is:

\begin{align}
    \Delta^g f = \frac{1}{\sqrt{\det g}} \partial_i(g^{ij} \sqrt{\det g}\ \partial_j f) \,,
    \label{eqn:Laplace-Beltrami-expl}
\end{align}

\noindent with $(g^{ij})={\bf g}^{-1}$. Inserted into the surface integral and integrated by parts we obtain: 

\begin{align}
\begin{split}
 \gamma \int_{\Gamma^{int}} \kappa{\bf w} \cdot {\bf n}\;d{\bf x} &\stackrel{(1)}{=}  \gamma \int_{\Gamma^{int}}  {\bf w} \cdot  \Delta^g {\bf {\tilde{X}}}\;d{\bf x} \\
&\stackrel{(2)}{=} \gamma \int_{\Gamma^{int}}  {\bf w} \cdot \partial_i(\tilde{g}^{ij} \sqrt{\det \tilde{g}}\ \partial_j {\bf {\tilde{X}}}) \ d{\boldsymbol{\tilde{\xi}}} \\
&\stackrel{(3)}{=} - \gamma \int_{\Gamma^{int}_{\tilde{\xi}}}  \left( \frac{\partial}{\partial \tilde{\xi}^i} {\bf w}^k \right) \tilde{g}^{ij} \left( \frac{\partial}{\partial \tilde{\xi}^j} {\bf {\tilde{X}}}_k \right) \sqrt{\det \tilde{g}}\ d{\boldsymbol{\tilde{\xi}}} \\
&\stackrel{(4)}{=} \sum_{\mathrm{interface\ elements}} - \gamma \int_{I}  \left( \frac{\partial}{\partial \xi^i} {\bf w}^k \right) g^{ij} \left( \frac{\partial}{\partial {\xi}^j} ({\bf {{X}}}^e)_k \right) \sqrt{\det g}\ d{\boldsymbol{\xi}} \,.
\end{split}
\label{st-explicit}
\end{align}

\begin{figure}
\centering
\includegraphics[scale=0.28]{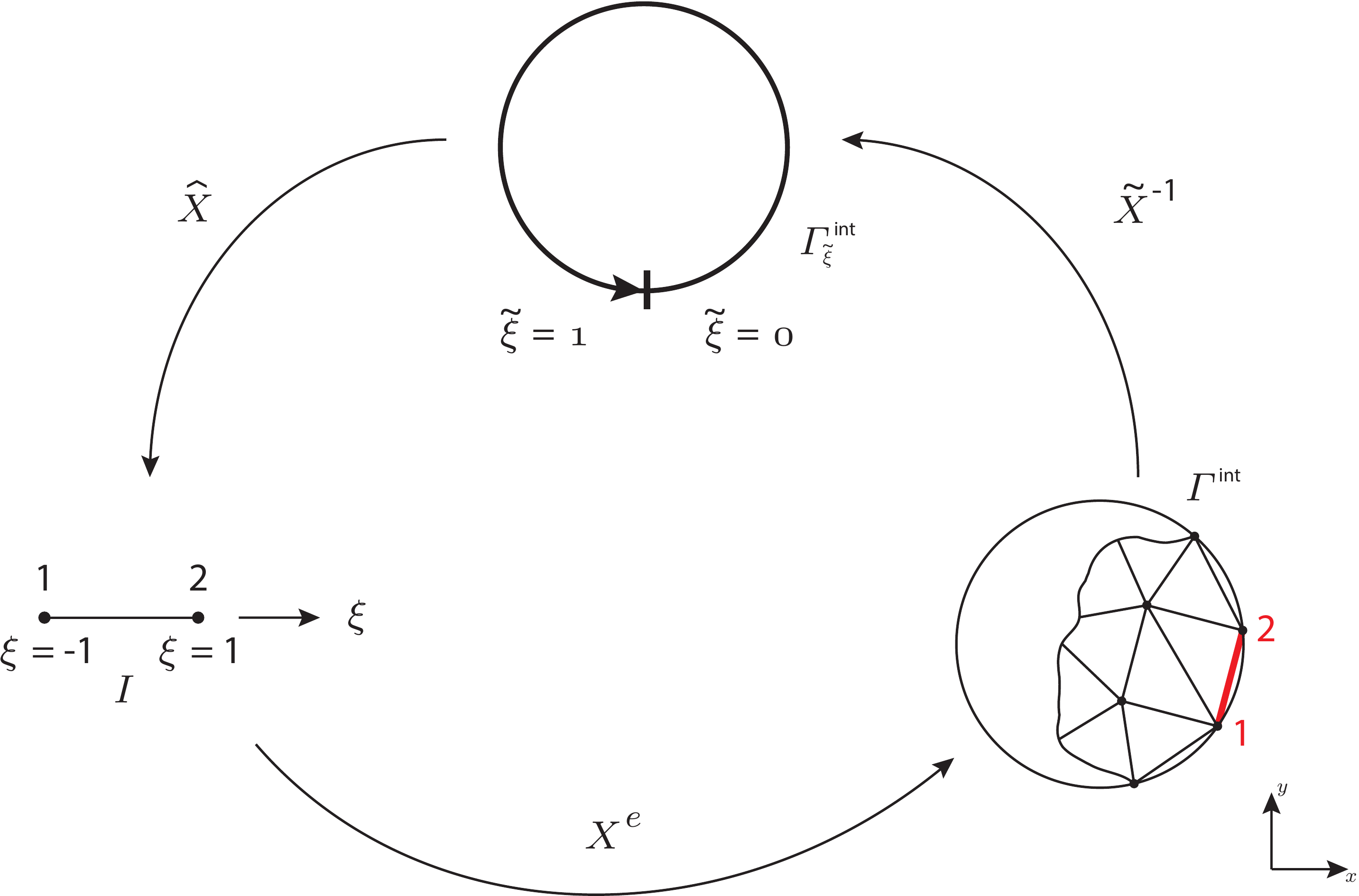}
\caption{The interface $\Gamma^{int}$, which is given in ${\bf x}$-coordinates, is reparametrized with the parameter $\tilde{\xi} \in [0,1]$. This is a smooth embedding and the Laplace-Beltrami technique can be applied. Partial integration can be performed. The integrals are then as usual computed on the reference element in $\xi$ coordinates. The 1D reference element is mapped to a portion of the boundary of a 2D mesh (here, we see an excerpt). }
\label{fig:immersion}
\end{figure}

The steps in Equation~\eqref{st-term} are detailed  in the following: 

\begin{itemize}
\item[(1)] The curvature vector is replaced using the Laplace-Beltrami technique.
\item[(2)] We insert definition~\eqref{eqn:Laplace-Beltrami-expl} for the Laplace-Beltrami operator. The integral is transformed to a parametrization of the interface in $\boldsymbol{\tilde{\xi}}$-coordinates as given in Figure~\ref{fig:immersion}. This parametrization is still smooth.
\item[(3)] Under the assumption that ${\bf w}$ is sufficiently smooth, integration by parts can be performed. 
\item[(4)] The interface is discretized using a boundary conforming finite element mesh. The integral is computed element-wise through a reparametrization into element reference coordinates  $\boldsymbol{\xi}$.
\end{itemize}

For the integration by parts, we assume that the surface is closed and therefore, the boundary term vanishes. This is an assumption regularly made \cite{Ganesan_2007,Gross2010}.

In order to give an example of the actual implementation, we will demonstrate how to evaluate the immersion term in Equation \eqref{st-explicit} in a standard finite element setting. $ {\bf X}^e $ is defined as the mapping from reference coordinates ${\boldsymbol \xi} $ to physical coordinates ${\bf x}$:

\begin{align}
{\bf X}^e: {\boldsymbol \xi} \mapsto {\bf x} \quad {\bf X}^e({\boldsymbol \xi}) = \sum_{k=1}^{nen} N_k({\boldsymbol \xi}) {\bf x}^e_k.
\end{align}

Here, $N_k$ denotes the finite-element shape function for node $k$ and $nen$ indicates the number of nodes per element. Note that we refer to the boundary element as indicated in Figure ~\ref{fig:immersion}.

The derivative of $ {\bf X} $ with respect to the (here only one) reference coordinate $\xi$ evaluated at a specific point $\xi_a$ is then:  

\begin{align}
      \left.\displaystyle\frac{\partial {\bf X}}{\partial \xi}\right|_{\xi_a} =\sum_{k=1}^{nen} \left.\displaystyle\frac{\partial N_k}{\partial \xi} \right|_{\xi_a}  {\bf x}_k \,.
\end{align}

\noindent This concludes the explicit case. 

In the implicit case, the starting point is to express the sum of principle curvatures of the interface $\Gamma^{int}$ by a relation similar to  Equation~\eqref{eqn:kappa-explicit} as \cite{Gross2010,Ganesan_2007}:

\begin{align}
\label{identity}
\kappa {\bf n} = \Delta_{\Gamma} id_{\Gamma^{int}}.
\end{align}

In Equation \eqref{identity}, $id_{\Gamma^{int}}$ signifies the identity mapping on ${\Gamma^{int}}$ defined as $id_{\Gamma^{int}}({\bf x}): {\Gamma^{int}} \mapsto {\Gamma^{int}} := {\bf x} = (x_1, x_2, \dots , x_{nsd})^T$. Note that this definition implies that $\nabla ( id_{\Gamma^{int}})$ are exactly the basis vectors of $\mathbb{R}^{nsd}$. Furthermore, $\Delta_{\Gamma}$  denotes again the Laplace-Beltrami operator. However in this variant, a different definition is employed. This is also the reason, why for now a different symbol has been chosen. Apart from the definition in Equation~\eqref{eqn:Laplace-Beltrami-expl}, it can also be defined using the tangential derivative and the tangential divergence \cite{Gross2010}:

\begin{align}
\Delta_{\Gamma}  = \nabla_{\Gamma} \cdot \nabla_{\Gamma},
\label{eqn:Laplace-Beltrami-implicit}
\end{align}

\noindent where the tangential derivative $\nabla_{\Gamma}$ of an arbitrary function $f$ is defined as:

\begin{align}
\nabla_{\Gamma} f := \nabla f - ({\bf n} \cdot \nabla f) {\bf n} = ({\bf I} - {\bf n}{\bf n}^T)\nabla f, 
\label{eqn:tangential-gradient}
\end{align} 

\noindent and the tangential divergence of ${\bf f}$ is defined as :

\begin{align}
\nabla_{\Gamma} \cdot {\bf f} := \nabla \cdot {\bf f} - {\bf n}^T (\nabla {\bf f}) {\bf n}. 
\label{eqn:tangential-divergence}
\end{align} 

\begin{figure}
\centering
\includegraphics[scale=0.3]{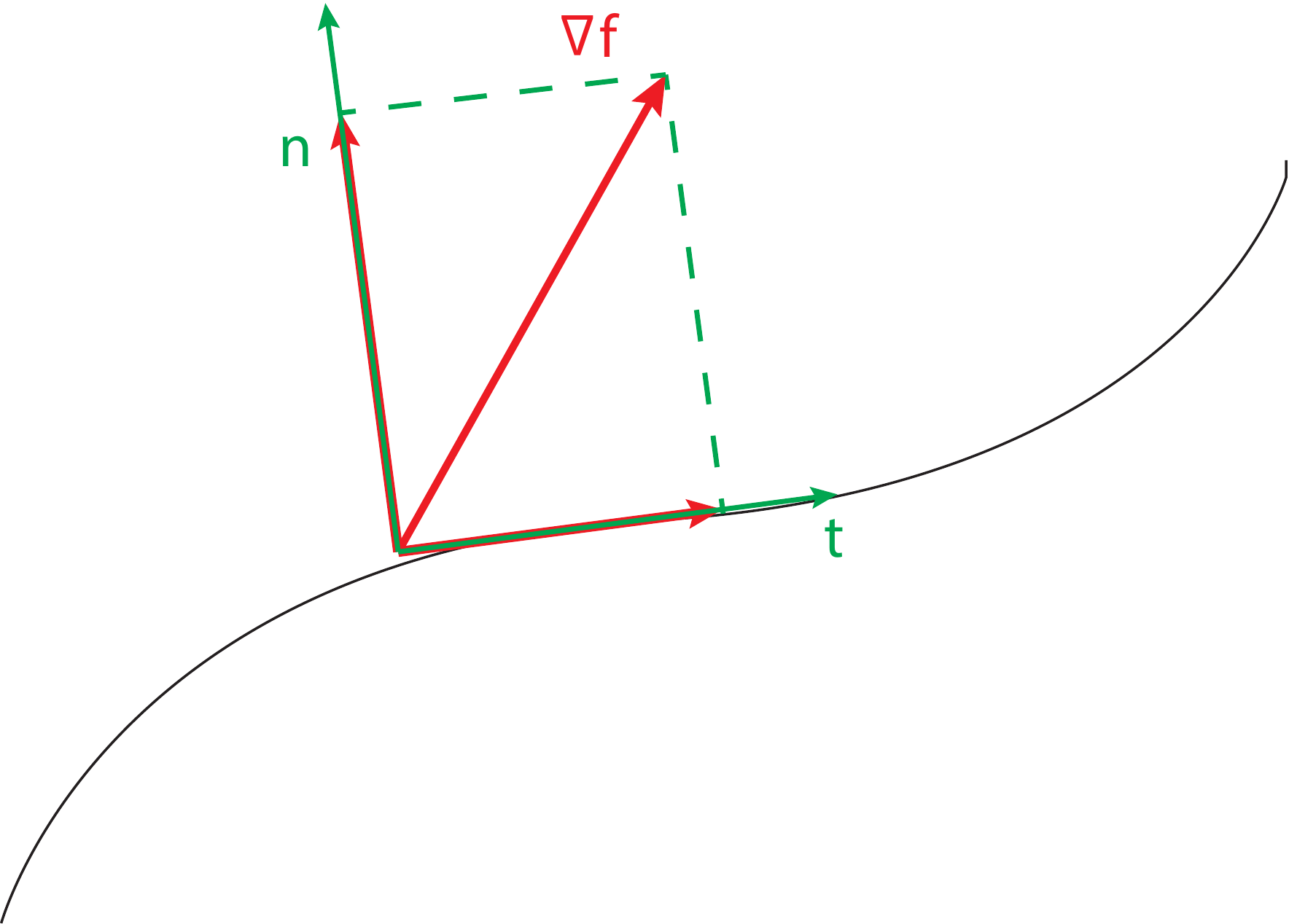}
\caption{The tangential gradient of a function $f$ can be evaluated by computing the full gradient and subsequently subtracting the normal component of the gradient. }
\label{fig:tangential-gradient}
\end{figure}

Using these relations to replace $\kappa {\bf n}$, we can again take advantage of the idea that integration by parts within the weak form can be used in order to reduce the order of derivatives. Again, the integral in equation \eqref{st-term} is restricted to the interface and therefore the integration by parts is not straightforward. In the case of a closed surface $\Gamma^{int}$, where the boundary term vanishes, this leads to \cite{Gross2010}:

\begin{align}
 \gamma \int_{\Gamma^{int}_t} \kappa{\bf w} \cdot {\bf n}\;d{\bf x} = - \gamma \int_{\Gamma^{int}}  \nabla_{\Gamma} {\bf w} : \nabla_{\Gamma} id_{\Gamma^{int}} \;d{\bf x}  = - \gamma \int_{\Gamma^{int}} tr(({\bf I} - {\bf n}{\bf n}^T) \nabla_{\Gamma} {\bf w}) \;d{\bf x}  \,.
\label{eqn:laplace-beltrami}
\end{align}

The assumption of a closed surface is valid in many applications where surface tension is important, such as for example drops and bubbles. Consult \cite{Gross2010}  for derivations for open boundaries and modifications for moving interfaces. Note that as compared to the derivation in the explicit case, the extraction of the tangential gradient using $ ({\bf I} - {\bf n}{\bf n}^T) $ is necessary. This stands contrary to the definition in Equation~\eqref{st-explicit}, where all gradients already point in tangential direction.

The obvious question to pose at this point is how the two presented expressions for the Laplace-Beltrami operator, $\Delta^g$ in Equation~\eqref{eqn:Laplace-Beltrami-expl} and $ \Delta_{\Gamma} $ in Equation~\eqref{eqn:Laplace-Beltrami-implicit}, can be related. In the context of this method, the crucial difference between the explicit and the implicit approach is that in the explicit case, the test function ${\bf w}$ can be expressed directly in terms of the reference coordinates ${\boldsymbol \xi}$; even along the interface. In the implicit case, however, the test function on the interface ${\bf w({\boldsymbol \xi})}$ needs to be expressed in terms of the test function in the bulk domain, which we will denote as ${\bf \tilde{w}}({\bf x})$ :

\begin{align} \label{eqn:implicit-test-function} 
&{\bf w({\boldsymbol \xi})} = {\bf \tilde{w}}({\bf X ( \boldsymbol{\xi}})), \\
&{\bf \tilde{w}}: {\bf X}(D) \subset \mathbb{R}^{nsd} \mapsto \mathbb{R}^{nsd} \,.
\end{align}

Over the course of the next equations, we will demonstrate that despite this difference between the approaches, Equations \eqref{st-explicit} and \eqref{eqn:laplace-beltrami} are to be considered as equivalent. If we insert definition~\eqref{eqn:implicit-test-function} into Equation~\eqref{st-explicit}, which has so far only been used for the explicit case, we obtain 

\begin{align}
\begin{split}
- \gamma \int_{I}  \left( \frac{\partial}{\partial \xi^i} {\bf w}^k \right) g^{ij} \left( \frac{\partial}{\partial {\xi}^j} ({\bf {{X}}}^e)_k \right) \sqrt{\det g}\ d{\boldsymbol{\xi}} \\
= - \gamma \int_{I}  \left( \frac{\partial}{\partial \xi^i} ({\bf X}^e)^l  \frac{\partial}{\partial x^l} {\bf \tilde{w}}^k \right) g^{ij} \left( \frac{\partial}{\partial {\xi}^j} ({\bf {{X}}}^e)_k \right) \sqrt{\det g}\ d{\boldsymbol{\xi}} \\
= - \gamma \int_{I}  \left(  \frac{\partial}{\partial x^l} {\bf \tilde{w}}^k \right) \underbrace{g^{ij} \left( \frac{\partial}{\partial \xi^i} ({\bf X}^e)^l  \frac{\partial}{\partial {\xi}^j} ({\bf {{X}}}^e)_k \right)}_{{(\bf P}(\boldsymbol{\xi}))_k^l} \sqrt{\det g}\ d{\boldsymbol{\xi}} \\
= -\gamma \int_{I} tr\left( \nabla {\bf \tilde{w}} |_{{\bf X}^e(\boldsymbol{\xi})}  {\bf P}(\boldsymbol{\xi}) \right) \sqrt{\det g}\ d{\boldsymbol{\xi}} \,.
\end{split} 
\end{align} 

Based on both the definition of the normal ${\bf n}$ and $ ({g}^{ij}) = ({g}_{ij})^{-1}  $

\begin{align}
{\bf P} {\bf n} = 0 \,, 
\label{eqn:P-def1}
\end{align}
and
\begin{align}
{\bf P} \frac{\partial}{\partial \xi^i} {\bf X}^e = \frac{\partial}{\partial \xi^i} {\bf X}^e  \,,
\label{eqn:P-def2}
\end{align}

\noindent hold. In other words, ${\bf P}$ is a projection onto the tangent space. Since in the 2D case $\{ {\bf n}, \frac{\partial}{\partial \xi^1} {\bf X}^e \}$ and in 3D $\{ {\bf n}, \frac{\partial}{\partial \xi^1} {\bf X}^e, \frac{\partial}{\partial \xi^2} {\bf X}^e \}$ form a basis, Equations~\eqref{eqn:P-def1} -- \eqref{eqn:P-def2} define ${\bf P}$ unambiguously. In particular, this means that 

\begin{align}
{\bf P} = {\bf I} - {\bf n n}^T \,.
\end{align}

With Equation~\eqref{eqn:Laplace-Beltrami-implicit} and $ {\bf P}^2 = {\bf P} $ we obtain

\begin{align}
\begin{split}
- \gamma \int_{I}  \left( \frac{\partial}{\partial \xi^i} {\bf w}^k \right) g^{ij} \left( \frac{\partial}{\partial {\xi}^j} ({\bf {{X}}}^e)^k \right) \sqrt{\det g}\ d{\boldsymbol{\xi}} = \\
- \gamma \int_{I} tr\left(  {\bf P} \nabla_\Gamma {\bf \tilde{w}} \right) \sqrt{\det g}\ d{\boldsymbol{\xi}} = -\gamma \int_{\Gamma^{int}} tr\left(  {\bf P} \nabla_\Gamma {\bf \tilde{w}} \right) d\Gamma \,.
\end{split}
\end{align}

The latter is equivalent to Equation~\eqref{eqn:laplace-beltrami}.

Note that in all cases of both explicit and implicit boundary descriptions, the assumption that Equations \eqref{st-term} and either \eqref{st-explicit}  or \eqref{eqn:laplace-beltrami} are completely equivalent only holds if $\Gamma^{int}$ is sufficiently smooth. This will in general not hold for the discretization of $\Gamma^{int}$~\cite{Gross2010}.

\subsubsection{Continuum Surface Force approach} \label{sec:CSF}
 
Another alternative to Equation~\eqref{st-term} is the Continuum Surface Force approach (CSF) as introduced in \cite{Brackbill92}. This is of particular interest when modelling topologically complex interfaces. The basic idea is to replace the interface coupling condition \eqref{eq:SurfTens} by a localized volume force term in the momentum equation \eqref{NavierStokes1}. The force term has the following form \cite{Gross2007b}:

\begin{equation}
{\bf f}_\Gamma = \gamma \kappa \delta_\Gamma {\bf n}.
\end{equation}

Here, $\delta_\Gamma$ refers to a smoothed Dirac $\delta$-function needed to select the surface force to the proximity of the interface, i.e., it selects a narrow band around the interface. ${\bf f}_\Gamma$ is designed in such a way that it exactly replicates the surface tension force per interfacial area that the surface integral would generate. Away from the interfacial region, the surface force is $0$. Throughout the transition region, all fluid quantities vary continuously.

\subsection{Solving additional equations on the interface}

\begin{figure}
\center
\psfrag{a}{$\Omega_1$}
\psfrag{b}{$\Omega_2$}
\psfrag{u}{${\bf u}$}
\includegraphics[width=5.5cm]{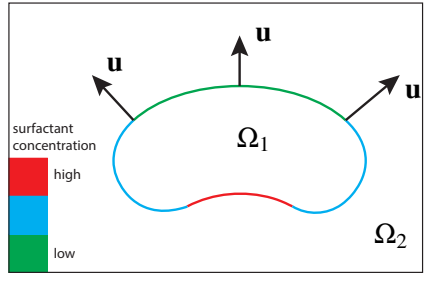}
\caption{Depending on the flow field and the diffusion coefficient, the surfactant concentration on the interface of, e.g., a drop, can be determined in every time step.}
\label{fig:surfactant}
\end{figure}

In some applications, it may be necessary to solve additional partial differential equations on the interface/boundary: an example for the topic of solving partial differential equations on hyper-surfaces (manifolds). One prominent example for such a case is the solution of the transport equation (advection-diffusion equation) as an additional equation to the governing equations of fluid mechanics \cite{Ganesan2009,James2004,Liu2010}, materials science \cite{Deckelnick2001,Fife2001}, and biology \cite{Barreira2011,Leung2003,Lowengrub2007}. The main topic is the consideration of insoluble surfactants (surface active agents), which adhere to the phase interface and influence the fluid flow in an either restraining or stimulating way. One type of surfactant known to all of us in our everyday lives is soap. It is important to note, that the surfactant is located exclusively at the interface/boundary and never in the bulk domain (cf. Figure~\ref{fig:surfactant}). Surfactants appear either as pollutants, which have accidentally entered the flow, or are purposefully added for example to decrease the surface tension and thereby increase the ability of a liquid to wet solid surfaces \cite{Davies1963}.

The two main driving forces for the transport are advection (movement of the interface) and diffusion (molecular diffusion along the interface) \cite{Gross2010}. Even though the governing equations for the scalar transport processes are often known and well understood, solving the respective partial differential equation on arbitrary and even moving manifolds is not straightforward. The effects brought about by the embedding, as already mentioned in Section~\ref{sec:Laplace-Beltrami} have to be accounted for. The solution approaches differ significantly depending on whether an explicit or implicit mesh description has been chosen (cf. Section \ref{sec:ITIC}). In the explicit approach, an interface mesh can be extracted from the bulk mesh thus easily defining the domain on which the transport equation is to be solved. In order to solve the proper equation, the differential operators of the partial differential equation have to be generalized in order to account for the arbitrarily shaped hyper-surface \cite{Dziuk2007,Ganesan2009}. Note that during the integration on the reference domain, it is mandatory to include the metric introduced in Equation \eqref{eqn:integration-on-reference-element}.  Implicit approaches \cite{Adalsteinsson2003,Bertalmio2001,Dziuk2010,Olshanskii2009,Xu2003} typically extend the scalar field away from the interface to the whole domain in order for all differential operators to be defined properly. One possible approach to the extension can be found in \cite{Chen1997,Chessa2003}. This obviously increases the dimension of the problem, but the transport equation can then be solved in the whole fixed computational domain. \cite{Gross2010} provides a good overview of different Eulerian and Lagrangian approaches. 

The solution of the transport equation requires a coupling of the flow solution with the transport. In addition, in many applications we would like to consider the effects that the surfactant has on the interface; i.e., a coupling of the transport solution to the flow equations. This coupling enters through the surface tension coefficient $\gamma$, which can be considered to be a function of the surfactant concentration $G$.  This leads to tangential forces along the interface, the so-called Marangoni convection.  \cite{Ganesan2009,Gross2010}
 
The standard advection-diffusion equation would be given as:

\begin{align}
\frac{\partial G}{\partial t} + \nabla \cdot ({\bf a}  G) - \nabla \cdot ({\bf d} \nabla G) = 0 \,,
\label{eqn:standard-adv-dif}
\end{align}

\noindent where ${\bf a}$ denotes the (possibly time dependent) advection velocity and ${\bf d}$ the diffusion coefficient matrix. In the case, where the above equation needs to be solved on a deformable hypersurface, some modifications need to be taken into account \cite{Ganesan2009,Gross2010}:

\begin{enumerate}
\item In the surfactant context, we will usually assume that ${\bf d}$ is constant and isotropic. Therefore $\nabla \cdot ({\bf d} \nabla G) = {\bf d} \Delta G $.
\item Since the equation is to be solved on an interface embedded into a larger domain, the differential operators in Equation~\eqref{eqn:standard-adv-dif} need to be replaced by the corresponding operators restricted to the interface: the gradient $\nabla$ is replaced by the tangential gradient  $\nabla_{\Gamma}$ (cf. Equation~\eqref{eqn:tangential-gradient}), the divergence $\nabla \cdot$ by the tangential divergence $\nabla_{\Gamma} \cdot$ (cf. Equation~\eqref{eqn:tangential-divergence}) and the Laplace operator  $\Delta$ by the Laplace-Beltrami operator $\Delta_{\Gamma}$ (cf. Equations \eqref{eqn:Laplace-Beltrami-expl} and \eqref{eqn:Laplace-Beltrami-implicit}). 
\item The advection velocity \--- and in the case of surfactant transport the advection velocity will be the fluid velocity ${\bf u}$ \--- needs to be restricted to its tangential component $ {\bf u}_{\Gamma} = ({\bf I} - {\bf n n}^T) {\bf u}$.
\item The term $\nabla_{\Gamma} \cdot ({\bf u}  G)$ can be rewritten as 
\begin{align}
\nabla_{\Gamma} \cdot ({\bf u}  G) = G \nabla_{\Gamma} \cdot {\bf u} +  {\bf u}_{\Gamma} \cdot \nabla_{\Gamma}  G \,.
\end{align}
In the cases of deformable interfaces, the first term ensures conservation of mass even under local changes in the free-surface area. Note that it is important to use $  {\bf u} $ instead of  ${\bf u}_{\Gamma}$. Otherwise, a change in surface area would not be possible. If the interface is fixed, this term is $0$ as ${\bf u} \cdot {\bf n} = 0$. 
\end{enumerate}

This leads to the final equation for the surfactant concentration $G$:

\begin{align}
\frac{\partial G}{\partial t} + ({\bf u}_{\Gamma} \cdot \nabla_{\Gamma}) G - \Delta_{\Gamma} G + G \nabla_{\Gamma} \cdot {\bf u} = 0 \,.
\label{eqn:surfactant-transport}
\end{align}

The above equation may be rewritten in many ways \cite{Adalsteinsson2003,Dziuk2007,Ganesan2009}, making it more accessible to the different numerical schemes.
 \section{Particle Methods} \label{s-particlemethod} 
 
Particle methods utilize mass-less particles distributed in an Eulerian mesh to capture the fluid flow and in particular the interface position. The most prominent particle method is the Marker-And-Cell (MAC) method developed by Harlow and Welch~\cite{Harlow65a}. It is one of the first examples of an interface capturing method; the numerous examples presented in~\cite{Harlow65a} attest to the method's extraordinary flexibility.

The original approach was developed for a single phase \--- i.e., $i=1$ in Equations \eqref{NavierStokes1} -- \eqref{NavierStokes2} \--- thus modelling a fluid surrounded by air. In later versions, the method could be extended to multi-phase flow. In all cases, the particle methods solve the ``one-fluid'' version of the Navier-Stokes equations, where the surface tension force is included directly and the fluid properties, e.g., $\rho$ and $\mu$, vary according to an indicator function. 

\begin{figure}
\centering
\includegraphics[scale=0.3]{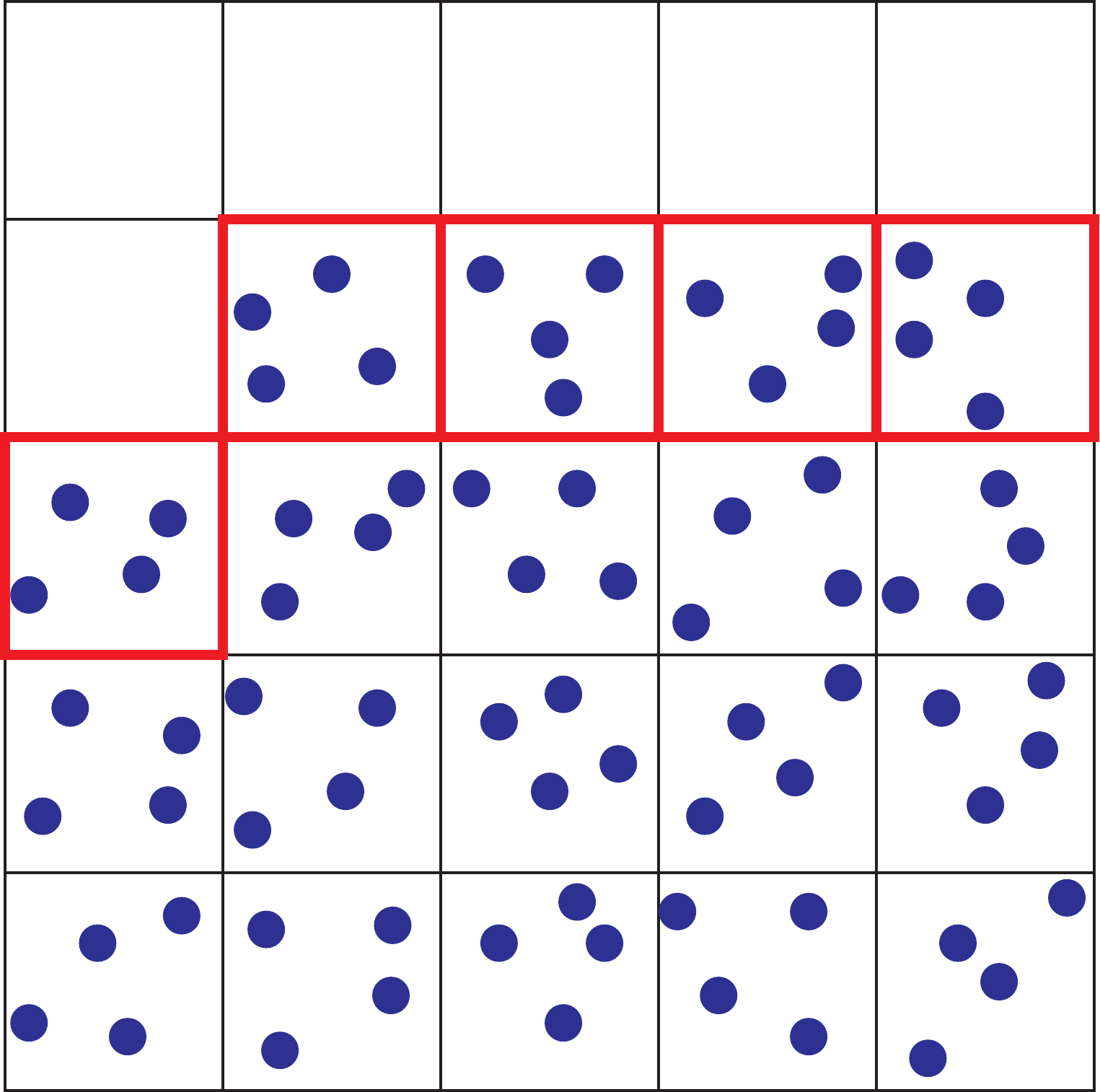}
\caption{In the Marker-and-Cell method by Harlow and Welch~\cite{Harlow65a}, the following convention is used for the identification of the fluid domain and the empty domain: Any cell not containing marker particles (blue dots) is identified as an empty cell. Cells with at least one marker particle and at least one common boundary with an empty cell are the interface cells (marked in red). Cells accommodating at least one marker particle and furthermore surrounded only by other cells containing marker particles are marked as fluid cells. }
\label{particlemethod}
\end{figure}

The method covers with a computational mesh the maximal possible fluid domain and subsequently distributes marker particles throughout the entire domain currently covered with fluid  (cf. Figure~\ref{particlemethod}). Any cell not containing marker particles is identified as an empty cell. Cells with at least one marker particle and at least one common boundary with an empty cell are the interface cells. As a last category, cells accommodating at least one marker particle and furthermore surrounded only by other cells containing marker particles are marked as fluid cells. 

Throughout the simulation, the marker particles are advanced with the local fluid velocity interpolated from the Eulerian grid.

The MAC method is strongly connected to finite-difference discretizations. In principle, however, it can also be combined with other methods, such as the  finite element method in \cite{Girault76}. 

In the MAC method, the velocity and pressure solutions are usually obtained sequentially. First, a preliminary velocity field that we will denote with ${\bf u}^{mom}$ is computed including all effects of the momentum equation (Equation~\eqref{NavierStokes1}) except for the pressure: Of the terms listed in Section~\ref{s-governing}, the inertia term, the friction term and the external forces contribute. The relevant equation reads (adapted from \cite{UnknownA} and \cite{Cline}):

\begin{align}
\frac{\partial^h {\bf u}^{mom}}{\partial^h t} = - ({\bf u}^n \cdot \nabla^h) {\bf u}^n + {\bf f}^n + \nu \Delta^h {\bf u}^n.   
\label{eqn:momentumMAC}
\end{align}

Here, the superscript $h$ on the operators indicate that they are discretized in an appropriate fashion, e.g., with finite-difference schemes. The superscript $n$ denotes the value of the corresponding field at the previous time step. Note that Equation \eqref{eqn:momentumMAC} can either be solved in one step or subdivided into several steps, computing the contributions of the individual terms one by one as done in~\cite{Cline}.

The computed ${\bf u}^{mom}$ will most likely not be divergence free and therefore, Equation~\eqref{NavierStokes2} is not fulfilled yet. However, the contribution of the pressure term to the velocity field has not been considered so far. This contribution is still available to adjust the velocity field by setting the pressure values appropriately. $ {\bf u}^{mom}$ and its correction factor $-\nabla p$ are inserted into Equation~\eqref{NavierStokes2}:

\begin{align}
\nabla^h \cdot (\underbrace{\frac{{\bf u}^{mom}}{\Delta t}}_{known} - \underbrace{\frac{1}{\rho}\nabla^h p}_{unknown}) = 0.
\label{eqn:PoissonPressure}
\end{align}

Again, the superscript $h$ refers to a discretized operator. Equation~\eqref{eqn:PoissonPressure} is a Poisson-type equation for the pressure. By solving this equation, the pressure values $p^{n+1}$ and with those, ${\bf u}^{n+1}$ can be obtained. This procedure is comparable to projection methods \--- or Chorin-Temam method \--- utilized in the context of the finite element method (\cite{Ern2003} and references therein). 

Notable is the development of SMAC (simplified-MAC) \cite{Amsden70} only five years after the original MAC implementation, which reduces the complexity of the original approach particularly when it comes to incorporating different types of boundary conditions. Its secret lies in the computation of the pressure field not based on a preliminary velocity field, but rather a potential function.   

\subsection{Properties}

With regard to the typical applications in fluid flow, MAC has the advantage that the tracer particles are indistinguishable from one another: The joining and separation of fluid parts can be effortlessly simulated. Nevertheless, it suffers from drawbacks especially regarding the description of free surfaces. Stress effects at the free surface, such as surface tension and contact angles are difficult to include. This is particularly true since, if used with finite difference, MAC is entwined with a discretization on staggered grids, i.e., storing velocity on cell faces and the pressure at the cell center, which has been found to provide higher stability \cite{Cline}. Boundary conditions at the free-surface can then for example only be set as pressures at the cell center of all cells identified as interface cells \cite{Easton72}, i.e., not at the actual boundary, or they may require complex surface fitting techniques such as those introduced in~\cite{Daly69a}. In addition, it is difficult to compute normals and curvature from the particles \cite{Li2008} and to control particle drifting \cite{Zheng2010}. In general, the method can suffer from stability problems, which may be difficult to detect, with apparently reasonable particle distributions resulting from gross approximation errors. Some modification were proposed to stabilize the MAC method in~\cite{Chan70a}. McKee \cite{McKee2004} identifies two further consequential restrictions: MAC could not be applied to arbitrary domain shapes and was, for a long time, restricted to two space dimensions. The latter was due to lacking efficiency in  the solution methods for the corrected pressure, preventing the use of a large number of particles.  In addition, care had to be exerted on the choice of time-stepping (no marker is allowed to cross more than one cell per time-step) \cite{McKee2004}. 

\subsection{Recent advances}

Recalling the drawbacks listed in the last section, two main points have recently been addressed in the area of MAC: the restriction to simple shapes in 2D and the topic of applying boundary conditions at walls and on the free surface. 

\begin{figure*}[htbp]
\center
\subfigure[Maximal domain]{
\includegraphics[width=4.6cm]{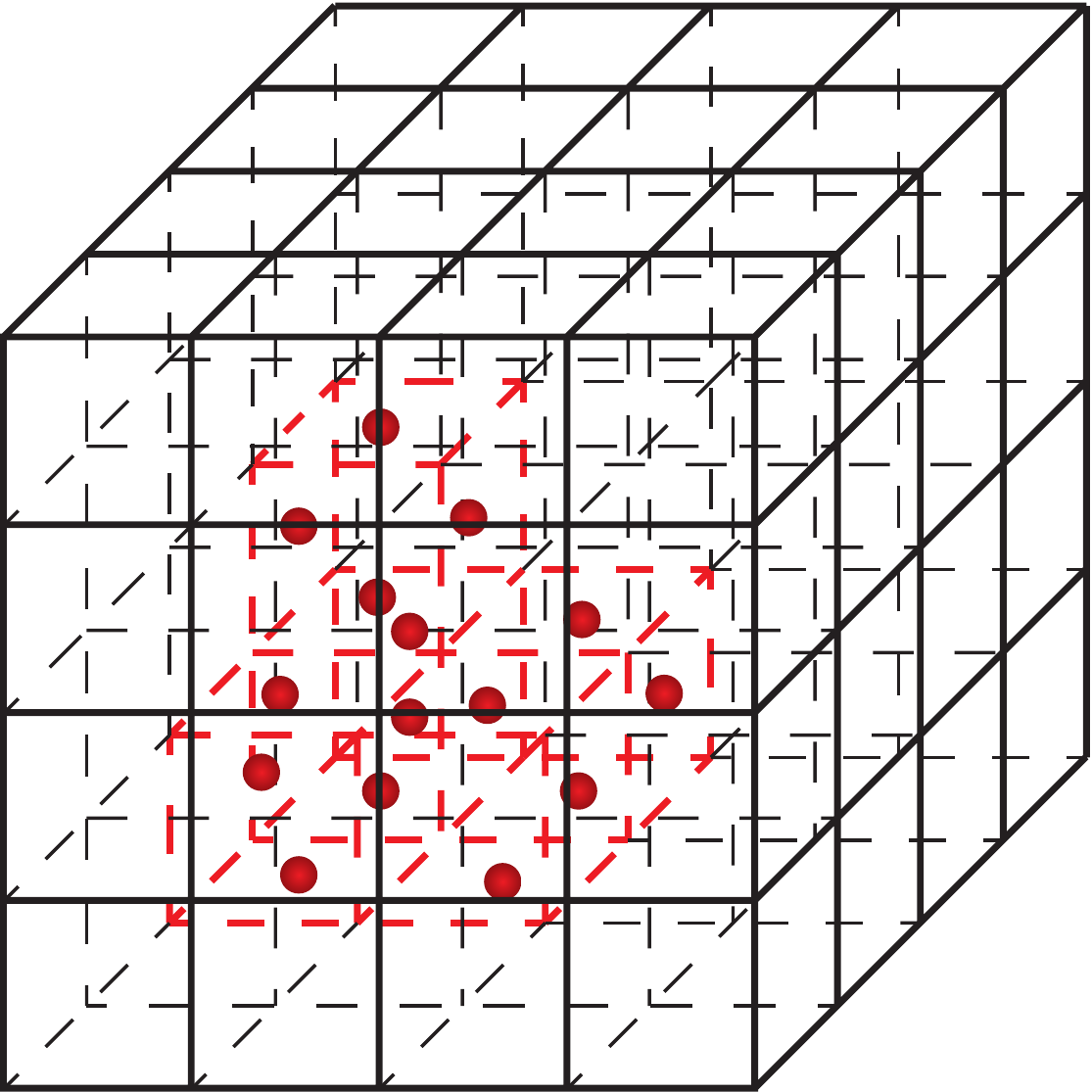}
}
\hfil
\subfigure[Reduced domain]{
\includegraphics[width=3.4cm]{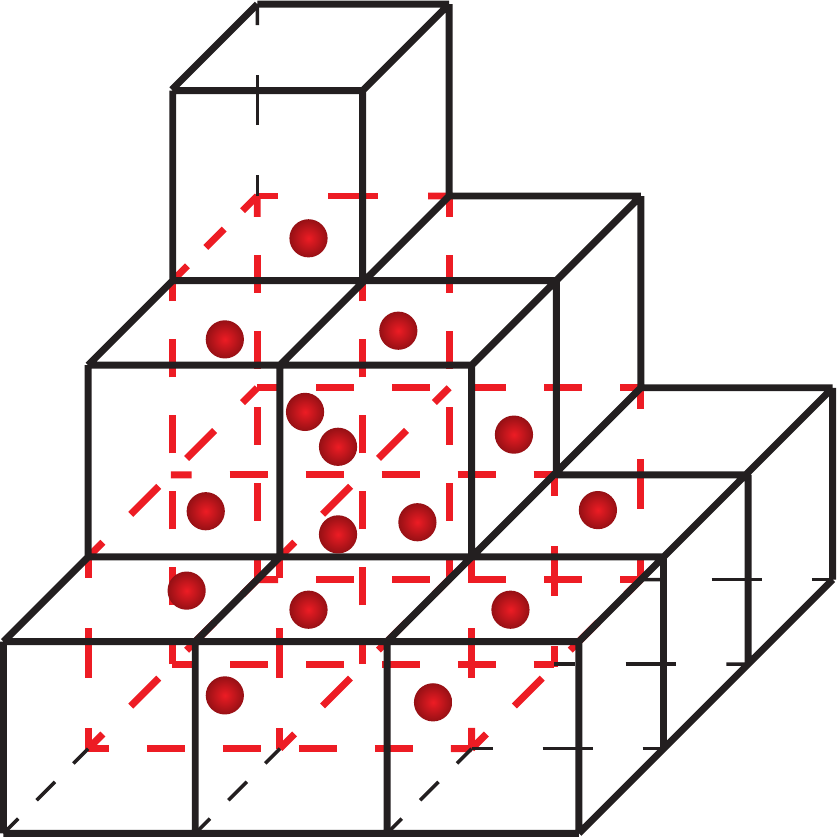}
}
\caption{(a) In the original MAC approach, the maximal possible fluid domain (black cells) was meshed, even if only a minor portion of cells (red cells) contain fluid. (b) In more recent approaches, the mesh consists only of the cells containing fluid plus an additional layer of ghost cells. This is sufficient as there is a time-step requirement, which guarantees that fluid particles cannot move through more than one cell per time step. \cite{Cline} }
\label{fig:MAC-ghost}
\end{figure*}

A development in the area of the former are GENSMAC \cite{Tome94}  and GENSMAC3D \cite{Tome2001}. This extension provides the applicability of MAC to domains of arbitrary shape and dimension as opposed to the straight lines required in the original implementation. This approach has been applied to viscoelastic extrudate swell in \cite{Mompean2011,Tome2008,Tome2012}, i.e., an application where the components of the stress tensor appear as additional unknowns in addition to the only primary variables of velocity and pressure. \cite{Santos2012} proposes an extension to an arbitrary number of phases. In \cite{Cline} a reduction of the mesh size and thus computational effort is proposed. The procedure is illustrated in Figure~\ref{fig:MAC-ghost}.

Boundary conditions on interfaces remained a challenge to MAC for decades. The only way to overcome this challenge was so far the combination with other methods, such as Lagrangian meshes (cf. Section \ref{s-interfacetracking}), level-set methods (cf. Section \ref{s-levelset}) or volume-of-fluid (cf. Section \ref{s-volumeoffluid}) \cite{McKee2004}.
Santos \cite{Santos2012} succeeded by including an additional Lagrangian mesh to track the free surface. \cite{Aubert2006} works with an ordered list of connected surface markers, which are advected along the streamlines of the flow field using a Runge-Kutta integration. Since the flow field will most likely disturb the homogeneity of  the marker distribution (concentration of streamlines or vortical flow), they are subsequently added, deleted and redistributed. \cite{Zheng2010} proposes a hybrid approach in such a sense, that the standard MAC method with its usual Eulerian grid is used for the fluid domain. In addition, special marker particles identify the interface, which is subsequently reconstructed using a Lagrangian mesh. Surface tension was also considered in \cite{Schroeder2012}, where the authors compared an implementation using a Lagrangian mesh with an implementation based on a particle level-set method. In \cite{Li2008}, a level-set approach is combined with Lagrangian marker particles as a means of level-set reinitialization. \cite{Leung2009a,Leung2011} present a method where the interface is tracked using mesh-less particles. The particles are explicitly connected to certain points of an underlying Eulerian reference grid. Connectivity among the particles is not required however. In each step, the interface is resampled based on a local reconstruction of the interface, on the one hand redistributing the particles and on the other hand updating the connecting points in the reference grid. The approach taken in \cite{Wang2012} is based on the point-set method first introduced in \cite{Torres2000}. It uses an indicator field (a function which is either $0$ or $1$ depending on the corresponding fluid domain), which is constructed from massless marker particles. Subsequently, the normals and curvatures can be computed from the indicator function, which has been smoothed using the reproducing kernel particle method (RKPM). A completely different approach has been adopted in \cite{Koshizuka96}, where all gradient terms in the Navier-Stokes equations are computed based on the interaction between all particles within a certain kernel. This method has been extended with more efficient solvers in \cite{Kakuda2013a}. 

The MAC-method served as inspiration for the volume-of-fluid method described in the next section. 

 \section{The Volume of Fluid Method} \label{s-volumeoffluid} 
 
The volume-of-fluid (VOF) method was developed in 1981 by Hirt~\cite{Hirt81a}. At the time, he worked in an environment where there were two options available to represent surfaces: either  marker particles (cf. Section~\ref{s-particlemethod}) or height functions combined with line segments. 
Height functions combined with line segments -- Hirt gives \cite{Hirt75} and \cite{Nichols71} as references -- were the first attempts to represent the actual interface. If the height function is used, the height of the interface as compared to a given reference line is stored for each discrete point in the domain. Problems occur for multiple-valued height functions as they would be seen for drops and when the interface slope exceeds the grid resolution. The height function approach can be generalized to ordered chains of line segments. Although the two given issues with the height function are resolved now, this approach has its drawbacks when it comes to curve intersections or an extension to 3D. 

\begin{figure}
\centering
\includegraphics[scale=0.2]{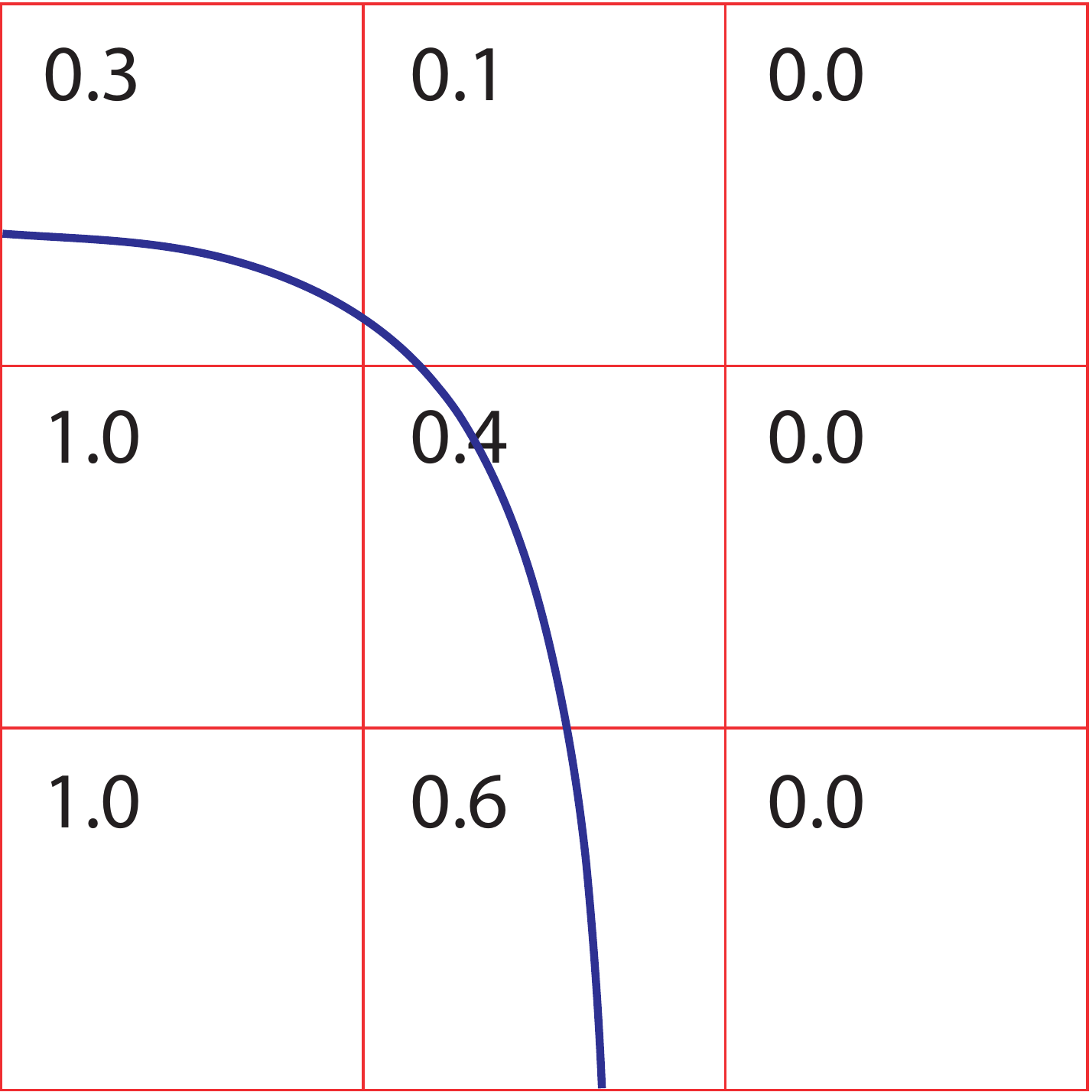}
\caption{Excerpt of a domain containing nine cells: The interface is indicated in blue. The area left of the interface is filled with fluid. In the VOF method, for each cell, the fraction of fluid contained in the cell is stored. Cells that are completely filled with fluid have the value $1$. Cells without any fluid have the value $0$. }
\label{VOF-concept}
\end{figure}

As described in the previous section, the MAC method does not define the interface, but rather the fluid regions, thus leading to no logical problems during interface intersection as well as a straightforward extension to 3D \--- apart from the excessive storage requirements, which were a concern to Hirt. Building on this scenario, he proposed the VOF method. This successor of the MAC concept replaces the discrete marker particles with a global field defining the location of the fluid(s). This global field is called the volume-fraction field $F$, sharply dividing fluid ($F = 1$) and empty ($F = 0$) areas. For each cell, $F$ takes on one specific value $0 \leq F \leq 1$, indicating the amount of fluid contained in each cell (cf. Figure~\ref{VOF-concept}). One way of defining $F$ is the use of a characteristic scalar field $\phi$. In the case of the VOF method, $\phi = \mathcal{H}({\bf{x}},t)$, with $ \mathcal{H} $ as the heaviside function. The volume fraction for each cell $\Omega$ can then be defined as: 
\begin{align}
F(\Omega,t) = \frac{1}{|\Omega|} \int_\Omega  \mathcal{H}({\bf{x}},t) d\Omega.
\end{align}
With only one value per cell, the storage requirements for this method are significantly reduced as compared to the MAC approach. A very important feature of the VOF method is its inherent mass conservation \cite{Caboussat2005}.

The advancement of the interface in time is governed by the following advection-type equation:

\begin{align}
\frac{\partial F}{\partial t} + \nabla \cdot F = 0.
\end{align} 

Due to its discontinuous character, the volume fraction function cannot be advected directly. Instead, it is connected with the mass conversation equation of the Navier-Stokes equations, thus including the volume-fraction field into the flow solution. The resulting equation can be reduced to the computation of fluxes across the element faces. 
In its most basic version, the VOF method has three major difficulties: 

\begin{enumerate}
\item the computation of shape derivatives (as for example needed for curvature evaluation) is obstructed by the use of a non-smooth interface representation,
\item as the volume fraction function is advected, the interface becomes more and more diffuse,
\item imposition of boundary conditions along the interface is impossible \cite{Hirt81a,Aniszewski2014}.
\end{enumerate}

One possible solution to all three problems is the reconstruction of a sharp interface in each time step. Options for this reconstruction will be discussed in the following section. Nevertheless, the VOF method is still utilized today in its original version, as for example in \cite{Choi2013}.

\subsection{Recent Methods for Interface Reconstruction}

Many advances in the VOF method are based on the aspect of reconstruction of the interface line. Over the past decades, numerous approaches have evolved. 
This large number of approaches already hints towards one major problem: As the interface shape is inferred from the volume data it is never unique. The reconstructed shapes are in general either piecewise linear or piecewise constant. \cite{Rider1998}   

\begin{figure}[htbp]
\center
\subfigure[SLIC]{
\includegraphics[width=4.5cm]{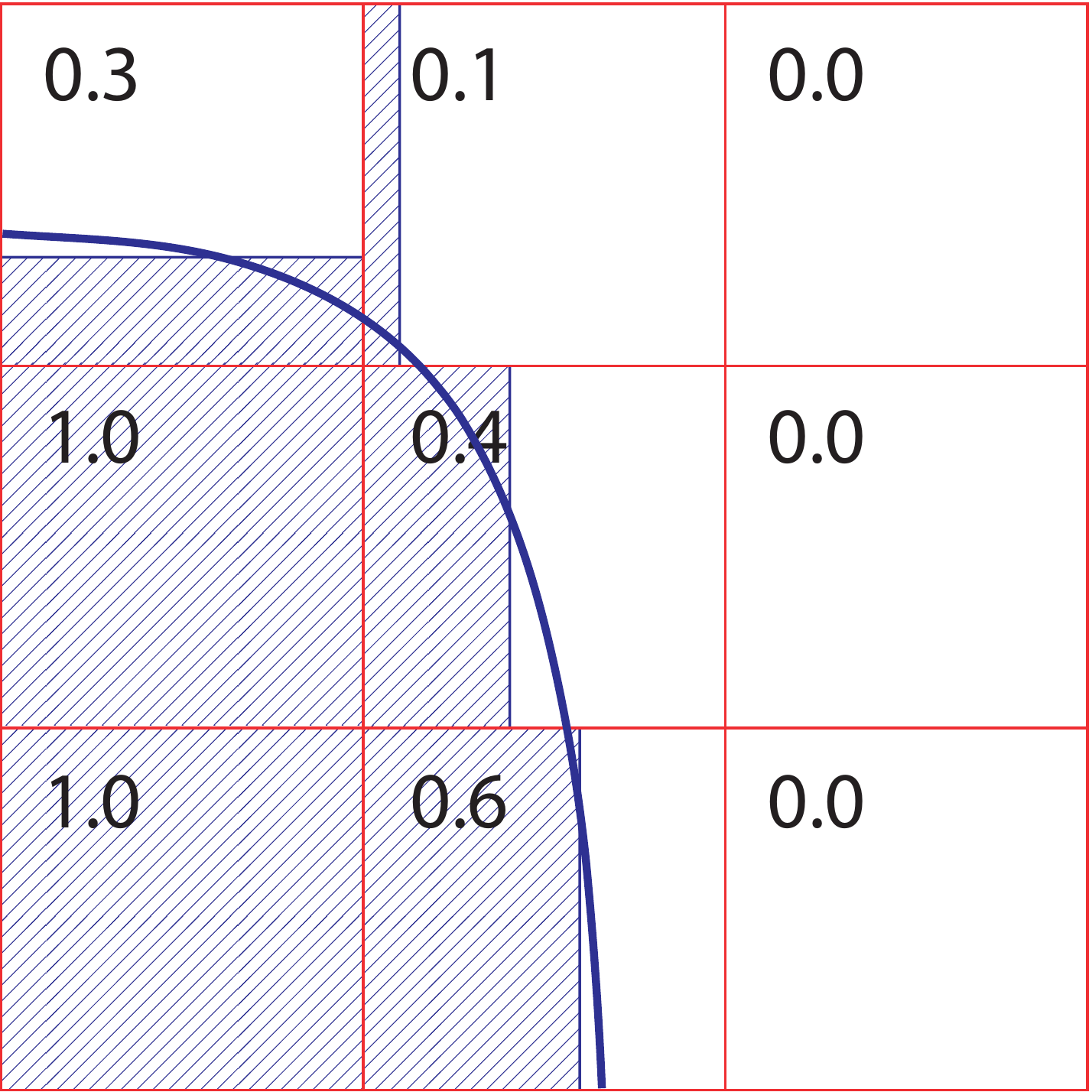}
}
\hfil
\subfigure[PLIC]{
\includegraphics[width=4.5cm]{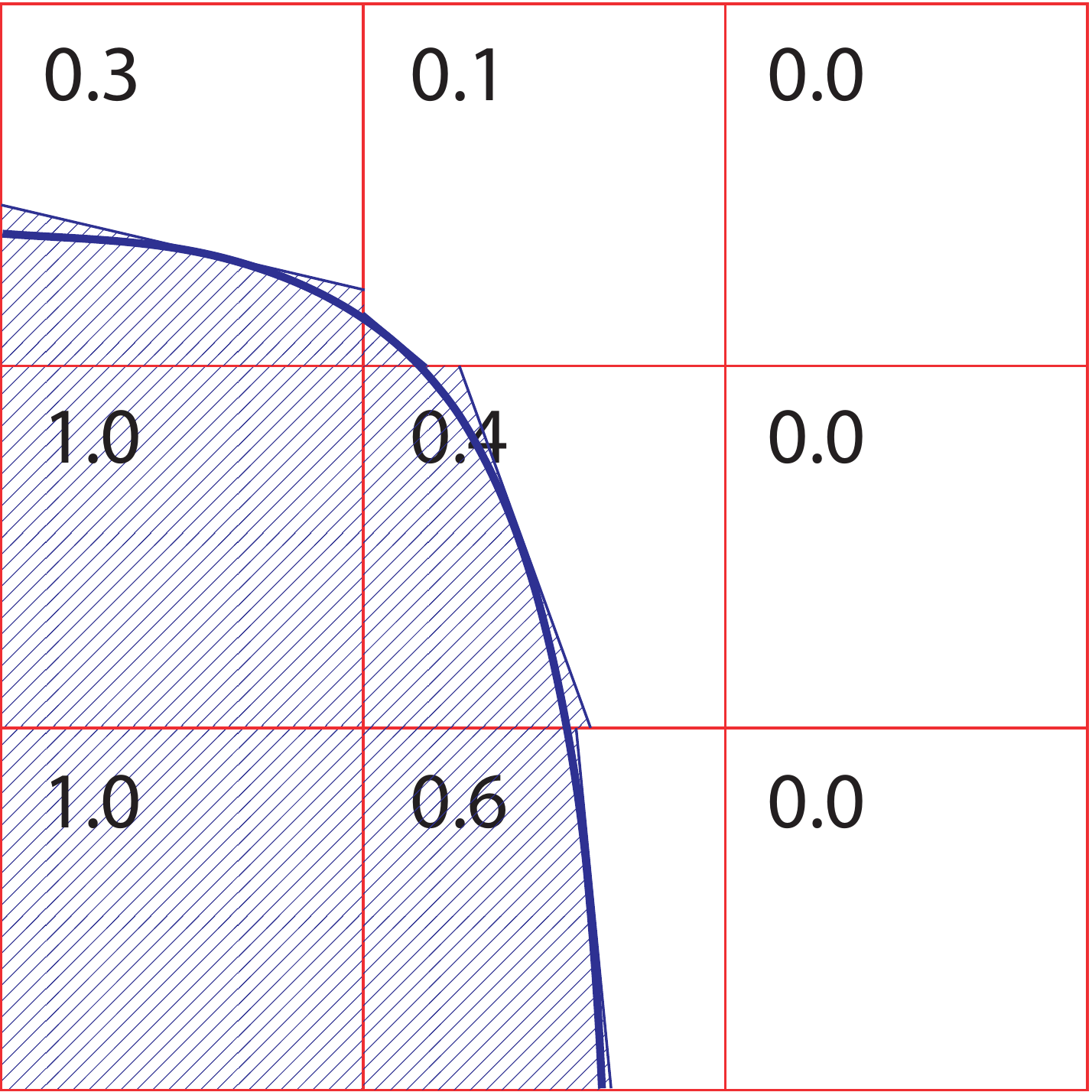}
}
\caption{The two main interface reconstruction approaches for the VOF method \--- SLIC and PLIC \--- are compared. (a) In the SLIC approach, the interface is reconstructed using line segments, which can be either parallel or perpendicular to the major flow axis. (b) In the PLIC approach, the interface is represented by a discontinuous chain of segments.}
\label{fig:VOF-Interface-Reconstruction}
\end{figure}

In the original implementation of VOF, a method that would later become known as  SLIC (Simple Linear Interface Calculation) \cite{Noh76} was proposed, where the interface is reconstructed using line segments, which can be either parallel or perpendicular to the major flow axis. The direction of the line segment is defined through the coordinate direction in which a larger change in the fluid volume occurs (cf. Figure~\ref{fig:VOF-Interface-Reconstruction}).  

A second reconstruction option is PLIC (Piecewise Linear Interface Calculation) -- \cite{Aniszewski2014,Caboussat2005}  and references therein. In this  method, the interface is represented by a discontinuous chain of segments each obeying the definition \cite{Aniszewski2014}:
\begin{align}
{\bf x n} = \alpha.
\end{align}
The parameter $\alpha$ can be determined using the condition of mass conservation. A central question for this scheme is the way, in which the normal is evaluated. The obvious choice, evaluating the normal the gradient of the volume fractions ${\bf n} = \nabla F / | \nabla F|$ approximated through finite differences as proposed in \cite{Youngs84} does not perform well \cite{Mencinger2011}, which brought up a variety of alternative choices \cite{Pilliod2004}. \cite{Mencinger2011} presents a generalization of PLIC to adaptive moving grids, \cite{Ito2013} to arbitrary unstructured grids. Note that, although based on piecewise linear functions, the PLIC approach cannot in general represent linear interfaces. \cite{Rider1998} proposes an appropriate extension.
The definition of the orientation of the line segments can be enhanced by spline interpolation of the interface in every cell \cite{Lopez2004}.
A special challenge is the treatment of more than two fluids, which connect at a triple point within a cell. In equilibrium, this point is connected to the contact angles of the particular fluids. In the instationary case, \cite{Caboussat2005} derives a relation to the normal vectors obtained within a PLIC setting.

All piecewise linear reconstructions, i.e., both SLIC and PLIC, suffer from the unphysical creation of droplets disconnected from the actual surface by pure construction of the numerical method \cite{Luppes2011,Kleefsman2005}.

A similar approach, using a local height function, has been introduced in \cite{Kleefsman2005} and was for example utilized in \cite{Luppes2011} and references therein: The method operates on 3x3 blocks of cells. After the orientation of the surface has been determined based on the values of the volume fraction field of surrounding cells, the surface height is determined row- or column-wise. Then, instead of updating the volume fraction values, the height function is updated. A subsequent adjustment step of the corresponding volume fraction values is required. 

\begin{figure}[h]
\centering
\includegraphics[scale=1.0]{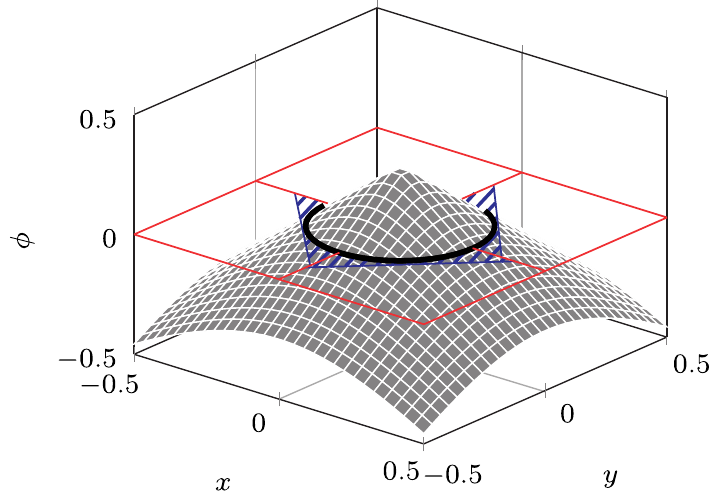}
\caption{Coupled Level-Set Volume-of-Fluid: Both the level-set function and the volume fraction are advected simultaneously, yet individually. The volume fraction interface description is then the basis for a correction of the zero level-set (cf. Section~\ref{s-levelset}).}
\label{fig:CLSVOF}
\end{figure}

A very fruitful alternative is the coupling of VOF with the level-set approach as it will be described in Section~\ref{s-levelset} \--- Coupled Level-Set Volume-of-Fluid (CLSVOF). Commencing with \cite{Sussman2000a}, the idea has been propagated to  combine the level-set approach \--- with its smooth interface \--- and the VOF approach \--- with its inherent mass conservation. Both the level-set function and the volume fraction are advected simultaneously, yet individually. The volume fraction interface description is then the basis for a correction of the zero level-set \cite{Aniszewski2014}. Such approaches, illustrated in Figure~\ref{fig:CLSVOF}, have for example been explored in \cite{Wang2012,Arienti2013}. A similar method is the reconstructed distance function (RDF) in \cite{Cummins2005}, except that here the distance function is deduced from the piecewise linear reconstruction of the volume fractions.

For some applications, simplified approaches that do not rely on a full interface construction can also lead to satisfactory results. For example, \cite{Aniszewski2014} names WLIC (Weighted Linear Interface Calculation) \cite{Yokoi2007,Marek2008} and THINC/SW (Tangent of Hyperbola Interface Capturing with Slope Weighting)  \cite{Xiao2005}.  WLIC is a simplification of PLIC. In a first step it treats each coordinate of the normal vector separately, thus computing fluxes along coordinate directions. In a second step these fluxes are averaged to account for slanted normal vectors.

\subsection{Curvature Computation}

The general definition of the curvature $\kappa$ in the VOF context is:

\begin{align}
\kappa = \nabla \cdot\left( \frac{\nabla F}{|\nabla F|} \right). 
\end{align}

However, being based on the heaviside function, the volume fraction function $F$ is not smooth enough to be differentiated directly. Options can be found in either reconstruction of the interface (as described in the previous section), or some type of smoothing. 

An alternative formulation is:

\begin{align}
\kappa = \nabla \cdot {\bf n}. 
\end{align}

Firstly, we will discuss methods for curvature computation that are based on interface reconstruction. In \cite{Kleefsman2005}, Kleefsman et al. suggest to compute the curvature from the height function $h$ using the relation:
\begin{align}
	\kappa = \frac{\partial}{\partial y} \left( \frac{\partial h/\partial y}{\sqrt{1+(\partial h/ \partial y)^2}} \right)	.
\end{align}
	The derivatives of the height function needed here are approximated using finite differences. Improved versions are among others offered by \cite{Lopez2009}, who introduced a correction scheme based on local error estimation, and \cite{Liovoc2010}, whose scheme involves a collection of possible height functions and corresponding curvatures, which are then sampled using quality statistics. In the context of PLIC, \cite{Meier2002} uses a local approximation of the curvature with polynomials. The polynomial coefficients are obtained from a least-squares fit. Raessi addresses the problem from a different perspective. In \cite{Raessi2010}, unit normals are advected, which are the principal constituent for both the reconstruction of the interface and the curvature calculation.

Secondly, we will address the topic of smoothing. In \cite{Cummins2005}, three different approaches to curvature approximation are compared: (1) In convolved VOF (developed in \cite{Williams2000}), $F$ is replaced by a smoothed version $\tilde{F}$, which has been mollified by a smoothing kernel $K$. $K$ could for example be defined as $K(r) = (1 - r^2)^4$ for $r < 1$.  The danger here is that the smoothing can be so strong that the curvature effects are diminished. (2) An improvement with respect to error, but not with respect to convergence behavior, was found using RDF. (3) Superior to both method was however a height function approach by Sussman \cite{Sussman2003}. 
Another commonly used approach to circumventing interface reconstruction is to employ the continuum surface force (CSF) model as described in Section~\ref{sec:CSF}. Very recently, Baltussen has suggested the tensile force method \cite{Baltussen2014}. This is based on a substantially different surface tension model: the tensile force model introduced in \cite{Tryggvason2001}. In this model, the surface tension force is determined as a sum of tensile forces, i.e., the force neighbouring cells exert on each interface element. This method does not require any curvature information. In \cite{Denner2014} a method the authors name CELESTE is presented. It revolves around a second-order Taylor expansion of $F$ from any given cell to its neighbouring cell. It is constructed from an overdetermined system of equations, utilizing a least-squares fit to determine both the normals and the curvature. 

An interface capturing method which, contrary to VOF, provides a smooth interface is the level-set method, which will be subject of the next section.

 \section{The Level-Set Method} \label{s-levelset} 
 
The level-set method is an interface capturing strategy. Level set methods~\cite{Osher88a,Chang96a} mitigate difficulties associated with a discontinuous volume-fraction function by using instead a smooth level set function $\phi$, separating areas filled with fluid A ($\phi < 0$) from those filled with fluid B ($\phi > 0$). The interface $\Gamma^{int}$ is then given by:

\begin{align}
\Gamma^{int} = \{ {\bf x} \in \Omega | \phi({\bf x},t) =0 \}.
\end{align}

 In most cases, the level-set function $\phi$ is defined as a signed distance function \cite{Osher2003}: 

\begin{align}
\phi({\bf x},t) = \pm \min_{{\bf x^*} \in \Gamma^{int}_t } \parallel {\bf x - x^*} \parallel, \quad \forall {\bf x} \in \Omega.
\end{align}

\begin{figure}[htbp]
\center
\subfigure[Implicit circular interface]{
\includegraphics[width=3.0cm]{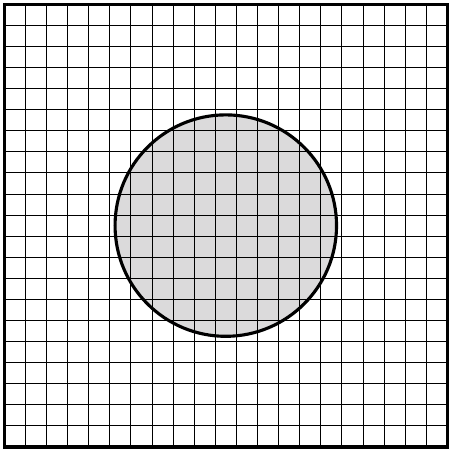}
}
\hfil
\subfigure[Level-set field of a circular interface]{
\label{resistance}
\includegraphics[width=5cm]{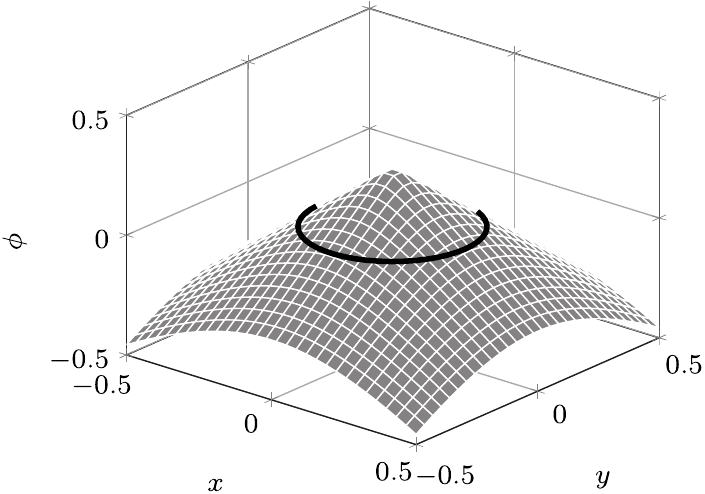}
}
\caption{A 2D domain containing a circular interface: In (a), an implicit definition of the interface is shown. Note that the elements of the mesh do not conform with the interface. The implicit description is realized through a signed distance level-set function outlined in (b).}
\label{fig:level-set}
\end{figure}

At each point ${\bf x}$ and $t$, the level-set function $\phi$ stores the shortest distance from the point to the current interface. Figure~\ref{fig:level-set} illustrates the concept by example of a circle. 

The method requires the definition of an initial level-set field from the original interface position:

\begin{align}
\phi({\bf x}, 0) = \phi^0({\bf x})  \quad \in \Omega. 
\end{align}

\begin{figure}[htbp]
\center
\subfigure[Variables jump at the (approximated) interface]{
\includegraphics[width=4.5cm]{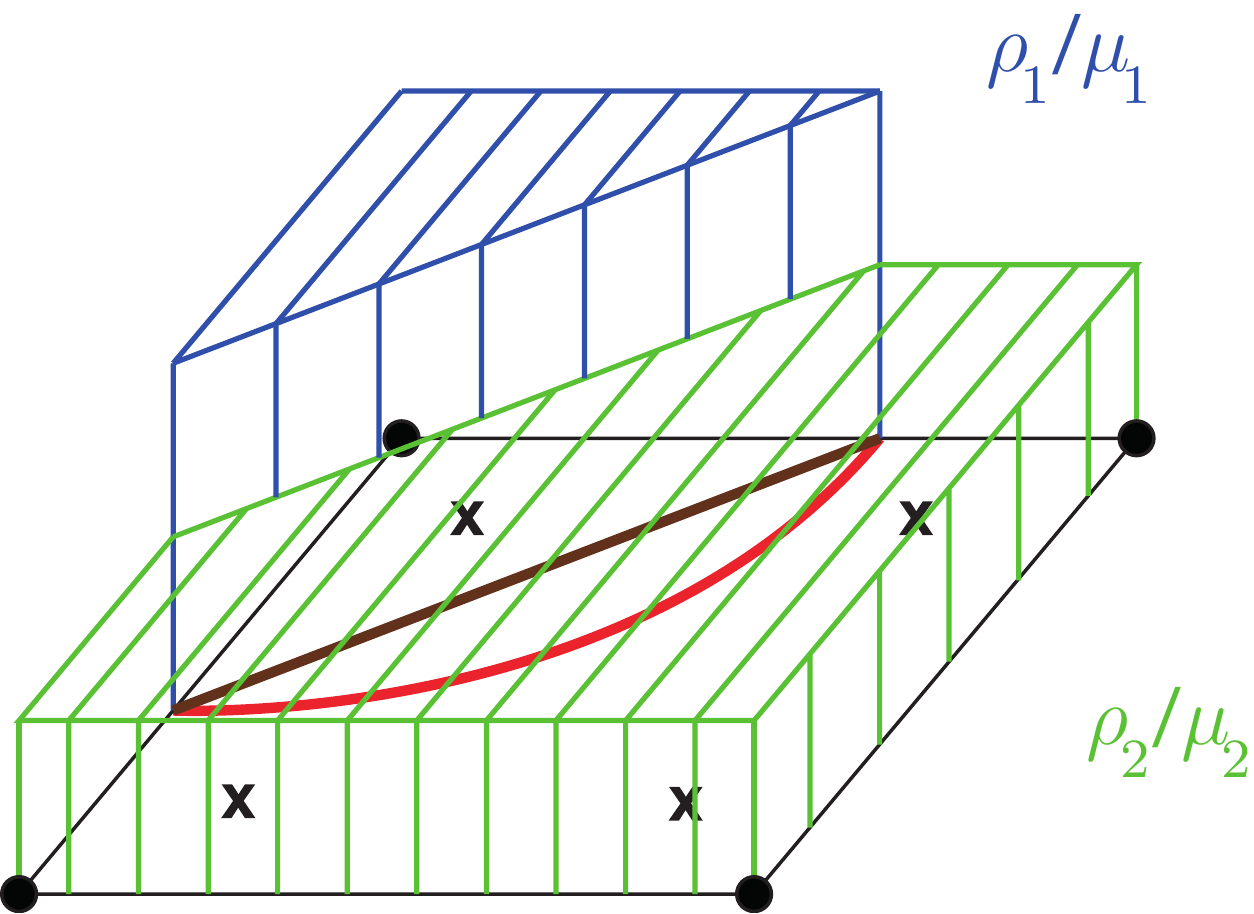}
}
\hfil
\subfigure[Variables are smoothed across the interface]{
\includegraphics[width=4.5cm]{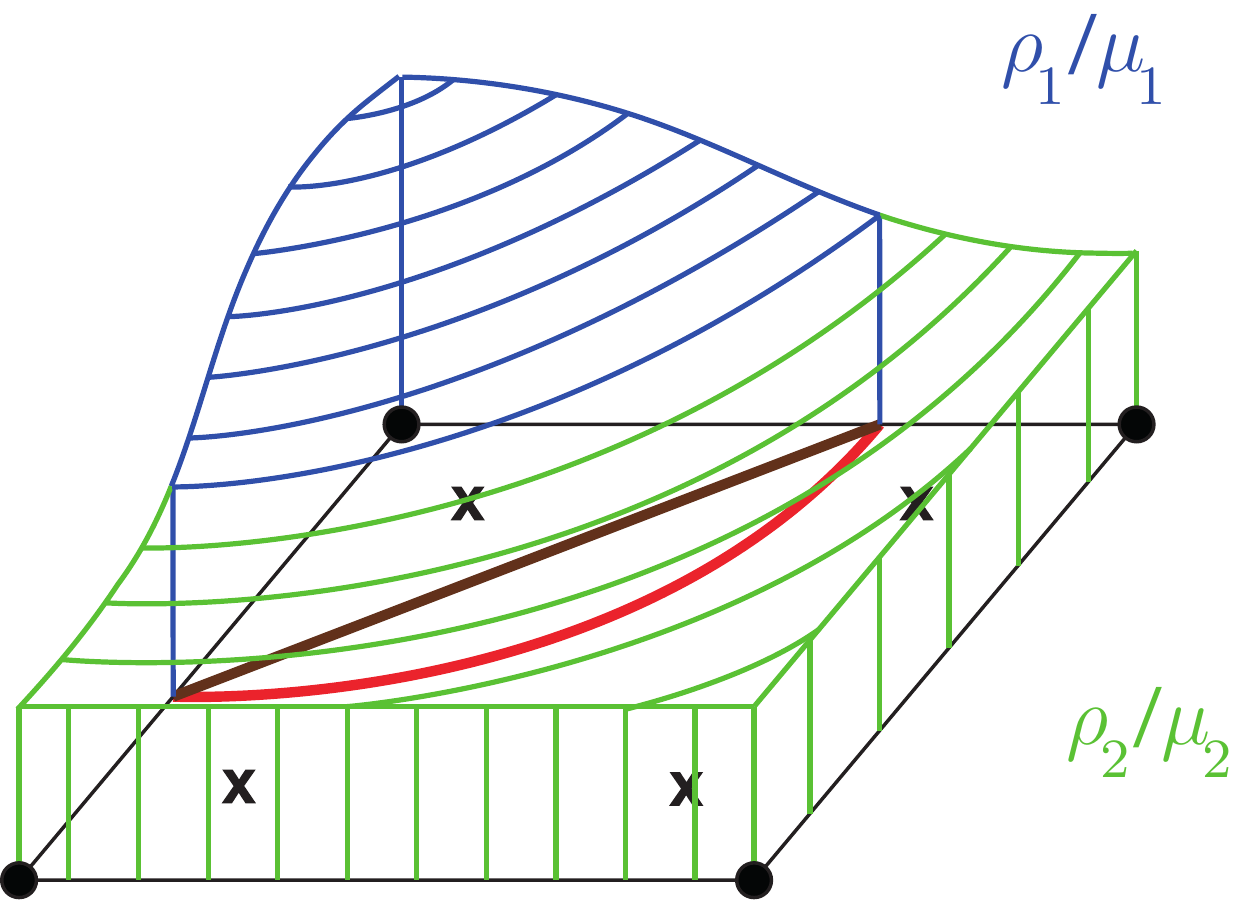}
}
\caption{A quadrilateral element with four nodes is intersected by the interface (red): The interface is approximated linearly (brown). (a) During integration, all Gauss points {\bf x} are associated with either the properties of fluid 1 ($\rho_1, \mu_1$) of of fluid 2 ($\rho_2, \mu_2$). (b) As an alternative, smooth transition of the properties can be generated throughout the element. Note that in both cases, an integration point might lie in the area between the interface as indicated by the level-set function and the linearized interface. One consequence of the approximation is then that this integration point will be associated with the wrong material domain.}
\label{fig:CutElement}
\end{figure}

Furthermore, at inflow boundaries, the level-set field has to be defined throughout the entire time domain, thus assuring that the level-set function is well-defined at boundary points, where information is transported into the spatial domain. During the progress of the simulation, the initial $\phi$ is then transported with the fluid velocity ${\bf u}({\bf x},t)$ using the following transport equation:

\begin{align}
\label{ls-transport}
\frac{\partial \phi}{\partial t} + {\bf u} \cdot {\bf \nabla} \phi = 0 \quad \mathrm{in} \ \Omega,\ t \in [0,T].
\end{align}
  
The value of the level-set function indicates, which density and viscosity values $\rho_i$ and $\mu_i$ are associated with which element, or, along the interface, only part of an element. To illustrate the latter case, Figure \ref{fig:CutElement} shows a sample element, which is cut by the interface. The interface is approximated (e.g., linearly) through the intersection points with the element edges. Depending on the side of the interface where an integration point is located, the material property associated with it is chosen:

\begin{align}
\rho, \mu (\phi) = \begin{cases} \rho_1,\mu_1, & \phi({\bf x}, t) < 0; \\ \rho_2, \mu_2, & \phi({\bf x}, t) > 0. \end{cases} 
\end{align}

If standard discretization methods are used, this procedure may lead to stability problems. It is therefore suggested to create a smooth transition of the material properties through the element. This can for example be achieved by including a smoothed Heaviside function into the formulation \cite{Sussmann1999}. The price to pay is, however, a  smeared out interface where it was once  infinitesimally thin. 

Already in its original set-up, the level-set method attains the main advantage of the interface-capturing methods: the topological flexibility. Furthermore, as opposed to other interface capturing methods, the evaluation of geometrical quantities is straightforward. As $\phi$ is a smooth function, it can be used directly to compute the interface curvature, normals, etc.  What is left is to cope with the challenges that this approach contains: treatment of discontinuities across the interface, ensuring mass conservation, applying boundary conditions on the interface, and evaluating equations on the interface. 

Note that the volume-of-fluid method of the previous section can be interpreted using the same scheme, except that $\phi= Heaviside$.

\subsection{Mass conservation}
\label{sec-ls-massconservation}

\subsubsection{Reinitialization}

\begin{figure}[htbp]
\center
\subfigure[Degraded signed distance function]{
\includegraphics[width=4.5cm]{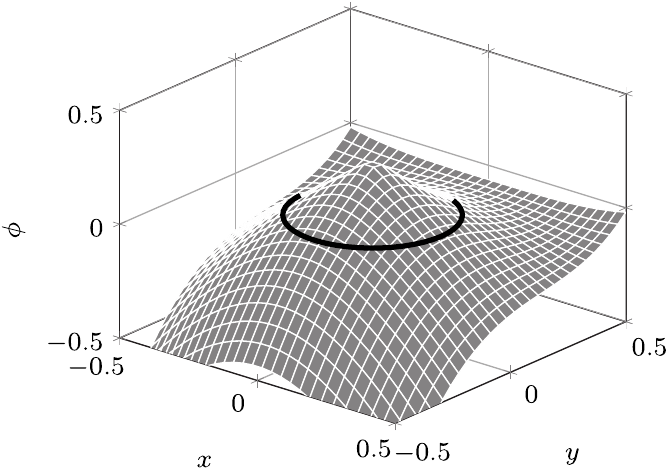}
}
\hfil
\subfigure[Locally reinitialized signed distance function]{
\label{resistance}
\includegraphics[width=4.5cm]{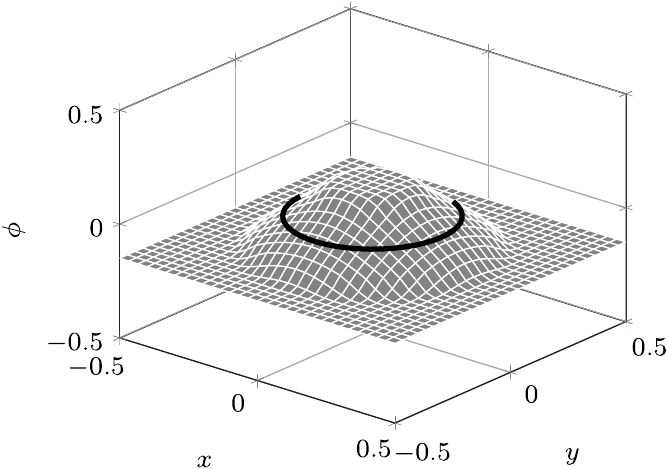}
}
\caption{To handle the degraded (as compared to Figure~\ref{fig:level-set}) level-set function (a), a reinitialization restricted to a band around the interface is performed (b).}
\label{fig:reinitialization}
\end{figure}

Since the velocity field entering the transport equation \eqref{ls-transport} is usually non-uniform, it will over time degrade the signed-distance property of the level-set function. Even though the interface is still represented correctly, this has severe influence on the quality of the computation of geometrical properties from the level-set function. Consequently, frequent reinitialization, i.e., restoring of the signed-distance character of $\phi$, is recommended. \\
Different approaches for reinitialization can be found in literature \cite{Hysing2005}. In direct reinitialization approaches, the closest distance to the interface is determined  individually for each node. The approach is based on identifying the intersections between the interface and element edges. Out of the collection of these intersection points, the Euclidean distance to the nearest point is computed for each node. This value is stored as the new level-set value for the node, keeping the sign intact. Unfortunately, this approach can quickly lead to efficiency problems. The efficiency can be improved by restricting the reinitialization to a narrow band around the interface \cite{Cho2010,Peng1999,Vaikuntam2008} (cf. Figure~\ref{fig:reinitialization}). 
As an alternative, fast search tree methods can be employed to organize the distance evaluations \cite{Marchandise2007}. Fast marching methods compute the distances sequentially from low to high, each computation building on the previously computed distances. \cite{Sauerland2013} features timings for the different approaches. 

\subsubsection{Mass correction}

One very undesirable side-effect of the reinitialization procedure is the mass loss: The reinitialization does not preserve the location of the interface. This effect is enhanced through the choice of temporal and spatial discretization of the transport equation \eqref{ls-transport} \cite{Reusken2011}. Investigated solution approaches are manifold. Some approaches aim at the minimization of the displacement during reinitialization \cite{Hartmann2008,Mut2006,Sussmann1999}. In \cite{Croce2010} an iterative mass correction is described, which provides global mass conservation. As mentioned in the Sections \ref{s-particlemethod} and \ref{s-volumeoffluid}, another possibility lies in the combination with particle or volume-of-fluid methods.

\subsection{Obtaining sharp interfaces}
\label{sec-ls-sharp-interfaces}

With the implicit definition through the level-set function, the discontinuity conditions along the interface will usually occur inside of elements, where, as described before, a smearing of the bulk properties is recommended. One possibility that allows to account for a sharp interface within the element is the extended finite element method (XFEM). XFEM is then often connected to the other improvement option: adaptive mesh refinement. 

\subsubsection{XFEM \-- the extended finite element method}

In the finite element method, the unknown function is usually interpolated with polynomial basis functions that are continuous within the elements and $C^0$ continuous across element faces. This approach is very much suitable for smooth unknown functions, but leads to problems if jumps occur. The XFEM circumvents this deficiency by local and well-targeted enrichment of the basis function space. The origin of XFEM lies in fracture mechanics \cite{Belytschko1999,Moes1999}. Overviews of the method can be found in \cite{Belytschko2009,Fries2010}. 
In XFEM, a sample function $g$ is approximated on a fixed mesh as:

\begin{align}
g^h({\bf x},t ) = \underbrace{\sum_{i\in I} N_i({\bf x}) g_i}_{standard\ FE} + \underbrace{\sum_{i\in I^*} N^*_i({\bf x}) \psi({\bf x},t) a_i}_{enrichment}.
\end{align}   

\begin{figure}[htbp]
\center
\includegraphics[width=6.0cm]{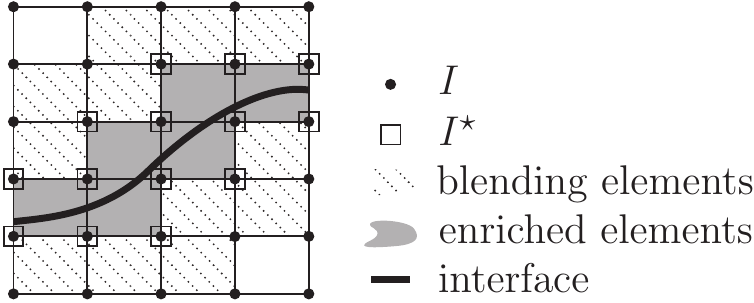}
\caption{Sample domain with an interface. The three different element types of the XFEM: regular, enriched, and blending elements, are indicated.}
\label{XFEM}
\end{figure}

$ N_i({\bf x}) $ and $ N^*_i({\bf x}) $ are the standard finite element shape functions for node $i$ with the nodal unknown $g_i$. The set $I$ is the set of all nodes in the domain, whereas $I^*$ contains all enrichment nodes. They are enriched using the enrichment functions $ \psi({\bf x},t) $, which are connected to additional XFEM unknowns $a_i$. Refer to Figure \ref{XFEM} for a sample mesh.

An enrichment function that is typically chosen for strong discontinuities is the sign-enrichment:

\begin{align}
\psi_{sign}( {\bf x},t) = sign(\phi( {\bf x},t)) = \begin{cases} -1\,, & \phi( {\bf x},t) < 0 \,, \\ 0 \,, & \phi( {\bf x},t) = 0 \,, \\ 1 \,, & \phi( {\bf x},t) > 0 \,. \end{cases}
\end{align} 

In the fluid flow problems considered in this paper, the scenario is as follows: Across the interface, a jump in density and viscosity, a kink in the velocity,  as well as both a jump and kink in the pressure need to be considered (cf. Figure \ref{scenario}).

Chessa and Belytschko \cite{Chessa2003a,Chessa2003b} use a kink enrichment for only the velocity approximation. Depending on whether surface tension is considered, a kink or jump enrichment of the pressure is added \cite{Fries2008,Minev2003,Rasthofer2011,Zlotnik2009}. K\"olke \cite{Koelke2005} enriches both spaces with a jump.  Even so it seems like a natural enrichment to account for the kink in the velocity field by appropriately enriching the basis, \cite{Coppola2009} states that this does not lead to an improved overall quality of the solution. \cite{Sauerland2010} even reports stability problems if this is included. As a consequence, several authors apply an enrichment for only the pressure \cite{Cheng2010,Coppola2005,Gross2007}.

\begin{figure}[htbp]
\center
\subfigure[Original interface]{
\includegraphics[width=2.5cm]{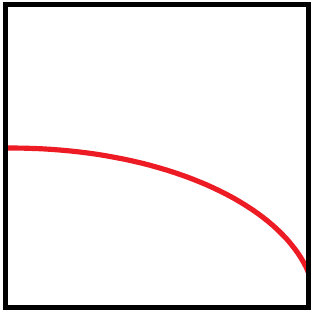}
}
\hfil
\subfigure[Linearized interface]{
\label{resistance}
\includegraphics[width=2.5cm]{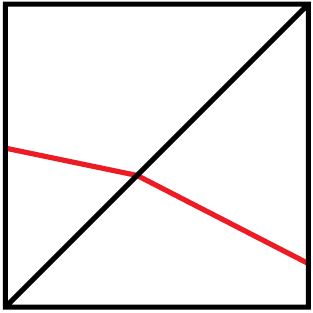}
}
\hfil
\subfigure[Sub-cells and Gauss points]{
\label{resistance}
\includegraphics[width=2.5cm]{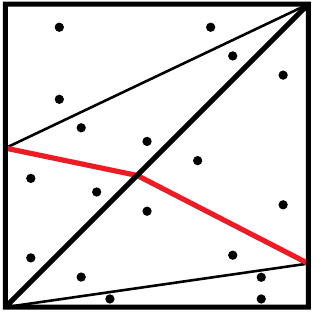}
}
\caption{Linearization of a curved interface in a quadrilateral reference element by decomposing the quadrilateral into linear triangular sub-cells, where finally the integration points are placed.}
\label{XFEM-quad}
\end{figure}

It is to be noted that the enrichment also entails special quadrature procedures. In order to be able to invoke standard Gauss rules, the interface cells are subdivided into triangles as indicated in Figure~\ref{XFEM-quad} \cite{Belytschko1999,Fries2010,Moes1999}.

\subsubsection{Adaptive mesh refinement}

\begin{figure}[htbp]
\center
\subfigure[Initial mesh]{
\label{AMR-1}
\includegraphics[width=3.5cm]{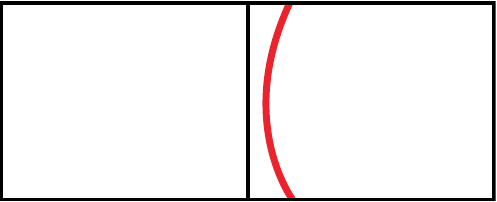}
}
\hfil
\subfigure[One level of refinement in interface elements]{
\label{AMR-2}
\includegraphics[width=3.5cm]{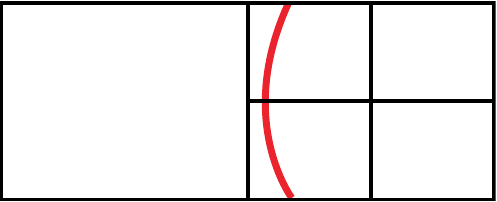}
}
\hfil
\subfigure[Two levels of refinement in interface elements]{
\label{AMR-3}
\includegraphics[width=3.5cm]{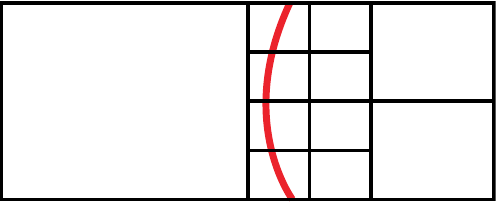}
}
\hfil
\subfigure[Additional refinement in neighbouring elements]{
\label{AMR-4}
\includegraphics[width=3.5cm]{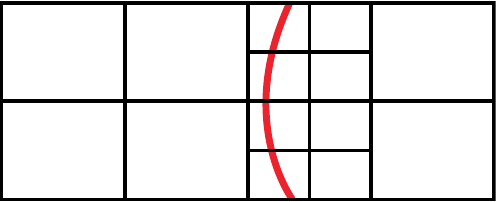}
}
\hfil
\subfigure[Level-set indicator field (interface) defined on a background mesh]{
\label{AMR-5}
\includegraphics[width=3.5cm]{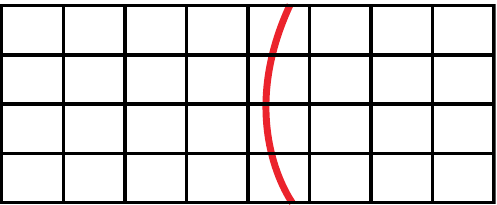}
}
\caption{Adaptive mesh refinement with respect to the interface (bold red line).}
\label{AMR}
\end{figure}

An alternative or an addition to XFEM can be found in adaptive mesh refinement (AMR). AMR can be used to increase the solution quality in the vicinity of the interface: the most interesting part of the computational domain. While global mesh refinement usually leads to prohibitively large mesh sizes, AMR provides an alternative at reasonable computational cost \cite{Babuska1983,Plewa2005}. In \cite{Cheng2010,Fries2011}, an adaptive refinement approach tailored to XFEM is introduced. Here, the level-set functions serves as an indicator for the refinement area. The interface elements are successively subdivided up to a desired level. 
Figure \ref{AMR} depicts the general refinement procedure for a simple setting. Here, two levels of refinement are applied. In a first step, the elements of the initial mesh, Figure \ref{AMR-1}, which are in contact with the interface, are subdivided once (cf. Figure \ref{AMR-2}). Those elements from the first refinement, which are cut by the interface are refined a second time, see Figure \ref{AMR-3}. Subsequently, further elements are subdivided step by step such that neighbouring elements differ by at most one level of refinement (cf. Figure \ref{AMR-4}). This assures a smooth variation of the element sizes. The level-set indicator field is defined on a fine background mesh with a mesh density corresponding to the highest refinement level (cf. Figure \ref{AMR-5}). Thereby, no interpolation of the level-set values at nodes created during refinement is required and the accuracy of the interface does not degrade.

\subsection{Computing geometrical quantities from the interface}

Since the level-set function $\phi$ is smooth, it can be used to compute geometric entities. As an alternative, the (linear) approximation of the level-set function within the elements could be employed. In terms of convergence order, \cite{Gross2010,Gross2007b} suggest the second option. 

\subsubsection{Normal vectors}

The normal vector is defined as:

\begin{align}
\label{normals-levelset}
{\bf n} = \left.\displaystyle\frac{{\bf \nabla} \phi( {\bf x},t)}{\parallel {\bf \nabla} \phi( {\bf x},t) \parallel} \right|_{\Gamma^{int}}.
\end{align}

Due to the signed distance property of $\phi$, the denominator in \eqref{normals-levelset} could essentially be dropped. However, since during the time-stepping, the level-set function is only an approximation of a signed distance function (cf. Section \ref{sec-ls-massconservation}), it is advisable to keep the normalization. 

\subsubsection{Curvature}

The curvature can be computed from the level-set function as:

\begin{align}
\kappa = {\bf \nabla} \cdot \left.\displaystyle \frac{{\bf \nabla} \phi( {\bf x},t)}{\parallel {\bf \nabla} \phi( {\bf x},t) \parallel} \right|_{\Gamma^{int}}.
\end{align}

This formulation is restricted to approaches, where the level-set function is not approximated linearly. An example is given in \cite{Olsson2007}, where a quadratic approximation in combination with a filtering technique is employed. However, even when theoretically possible, Marchandise points out in \cite{Marchandise2007} that it is in general not advisable to use the thereby computed curvature to evaluate the surface tension. One popular option is to resort to the Laplace-Beltrami technique already described in Section \ref{sec:Laplace-Beltrami}.

Particularly in combination with the XFEM, level-set provides very sharp interface descriptions. A contrasting approach is the phase-field approach described in the next section, which deliberately introduces diffuse interfaces.

 \section{Phase-Field Method} \label{s-phasefield} 
 
One of the oldest numerical methods in multi-phase problems is the phase-field method. According to \cite{Emmerich2002}, as early as 1873, the work of Gibbs on thermodynamics already served as a foundation~\cite{Gibbs1873}. 

The main difference compared to the previously introduced interface-capturing methods is that the phase-field method works with diffuse interfaces --- i.e., the transition layer between the phases has a finite size. There is no tracking mechanism for the interface, but the phase state is included implicitly in the governing equations. The interface is associated with a smooth, but highly localized variation of the so-called phase-field variable $\phi$. In two-phase problems, $\phi$ is a scalar value; if more phases are present, it can become vector-valued. $\phi$ is assigned a distinct value, e.g., $-1$, in phase $A$  and another distinct value, e.g., $1$, in the other phase $B$.  The interface could then be assumed to be located at $\phi=0$ \cite{Boettinger2002}; however, the phase-field equations never require the knowledge of the exact interface location.  $\phi$ will vary with space and time. Typically, we observe small variations in the bulk phases and rapid variations close to the interface \cite{Emmerich2002}. 

\begin{figure}
\centering
\includegraphics[scale=0.3]{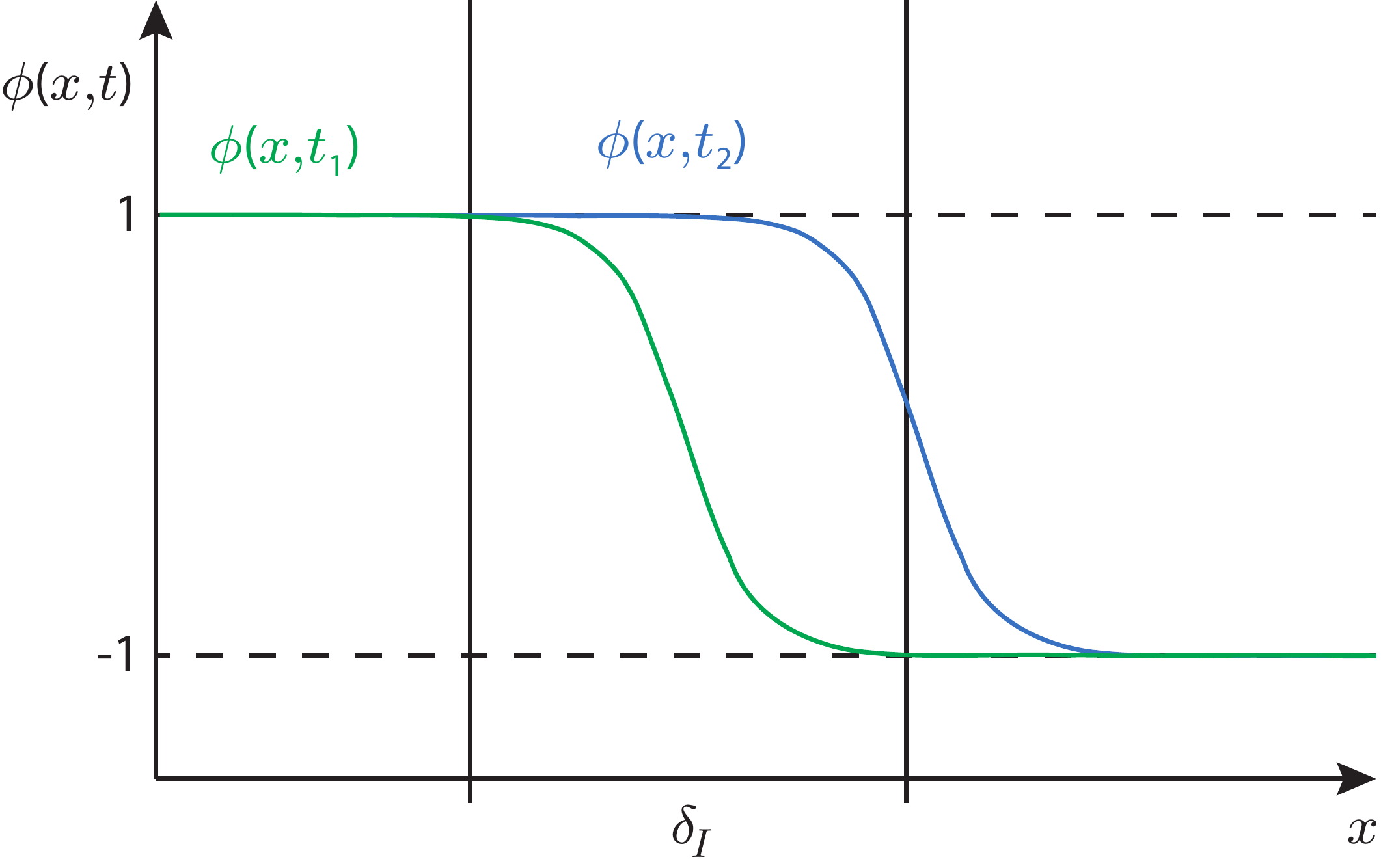}
\caption{Illustration of the phase-field variable: $\phi(x,t)$ takes on one value (e.g., $1$ in this case) in one phase and another value (in this case $-1$) in the second phase. The transition is rapid, but not sharp. We speak of a diffuse interface with interface width $\delta_I$.}
\label{fig:phase-field-variable}
\end{figure}

Comparing Figure~\ref{fig:phase-field-variable} and Figure~\ref{fig:level-set} gives the notion, that phase-field resembles level-set in many ways. Indeed, they have been applied for identical applications, however, there are some fundamental differences that should be considered in the choice of application. Major differences are: (a) The phase-field approach avoids setting continuum boundary conditions at the interface, which would usually be required to match the individual solutions of the bulk phases. Note also that one of these boundary conditions usually defines the advection velocity of the interface. What is particularly difficult about these boundary conditions is that they have to be set at a boundary whose shape and position is also part of the solution. Instead, these boundary conditions are substituted by the inclusion of the phase-field variable into the system. Consequently, the phase-field method is particularly useful in cases, where the continuum equations for a sharp interface do not  (yet) exist. Sharp-interface descriptions can be interpreted as a phase-field approach with infinitely small interface width. However, it is not considered to be very useful to use phase-field in such a sense \cite{Emmerich2002}. (b) In the phase-field method, the interface is not advected specifically as would be the case for level-set; moreover, the interface location is never computed explicitly \cite{Emmerich2002,Boettinger2002}.  

Phase-field can be seen as a computational method, in which the presence of a non-sharp interface dampens instabilities, or it can be seen as a physically motivated approach, where the governing equations do not impose sharp interface conditions \cite{Emmerich2002}. 

The derivation of the governing equations for a system described by a phase-field approach is strongly tied to the notion of thermodynamic equilibrium. In engineering applications, we usually deal with open or closed systems; and hardly with isolated systems. In open or closed systems, the thermodynamic equilibrium is not defined through the state of maximal entropy, as it would be for isolated systems, but through the extrema of some other thermodynamic state variable. In principle, we distinguish between four variations of the equilibrium by thermodynamic forces: mechanical, thermal, material, and chemical. \cite{Lucas2007} 

\begin{itemize}
\item {\it{Mechanical}}: Two systems under different pressures are inclined to reach mechanical equilibrium: mechanical work is performed. In such a case, the pressure is interpreted as a thermodynamic potential.
\item {\it{Thermal}}: Two systems with different temperature are inclined to reach thermal equilibrium: heat flow is induced. In such a case, the temperature is interpreted as a thermodynamic potential.
\item {\it{Material}}: A system containing either different material states of one material, e.g., solid and liquid phase, or two chemically inert materials, will tend to the material equilibrium. In such a case, the chemical potential is interpreted as the thermodynamic potential. 
\item {\it{Chemical}}: A system containing substances, which are chemically reactive, will develop a specific equilibrium between educts and products. Again, the chemical potential is interpreted as the thermodynamic potential. 
\end{itemize}

The state of equilibrium is characterized through the minimum of the thermodynamic potential $P$. The condition for equilibrium is defined as $dP = 0$. 

So far, we have stated that in order to describe a thermodynamic equilibrium of any kind, the definition of a suitable thermodynamic potential is essential. In realistic problems, several potentials interact. Therefore, the need for finding the extremum of one specific potential is now replaced by minimizing the free energy functional of the system. The free energy is the portion of the energy that is available to perform thermodynamic work. After performing such work, it is irreversibly lost \cite{Emmerich2002}. Depending on the applications, different definitions of a free-energy functional can be used: the Gibbs free energy formulation is for example applicable to systems where the temperature and the pressure remain constant, whereas the Helmholtz free energy formulation can be employed for situations where the temperature and the volume remain constant. In general we can write \cite{Emmerich2002}:

\begin{equation}
\label{eq:potential1}
{\mathtt{P}}(X_1,X_2, \dots, X_n) = \int_V P(X_1,X_2, \dots, X_n) dV,
\end{equation}

\noindent where ${\mathtt{P}}$ refers to the appropriate thermodynamic potential with $P$ signifying the respective density function. $X_1,X_2, \dots, X_n$ refer to the governing variables of the system, i.e., the independent state variable. Depending on the system, these could for example be temperature, pressure, the concentration of the phases, etc..  The free-energy functional together with the equations of motion leads to a full description of the system in the sense that the value of all state variables can be determined from the free-energy via extremal principles if we assume thermodynamic equilibrium. In detail, this means that through differentiation with respect to the individual variables, equations for all those variables can be obtained (cf. \cite{Emmerich2002} and references  \cite{Onsager1931,deGroot1984,Gyrmati1970,Prigogine1966} therein). 
In the phase-field method, an additional, in a sense artificial, dependence on the phase-field variable is generated. Equation \eqref{eq:potential1} now reads  \cite{Emmerich2002}:

\begin{equation}
\label{eq:potential1}
{\mathtt{P}}(X_1,X_2, \dots, X_n, \phi) = \int_V P(X_1,X_2, \dots, X_n, \phi) dV,
\end{equation}

\noindent where $\phi$ is the phase-field variable. 

In order to include $\phi$ into the formulation, a contribution of $\phi$ to the density function $P$ has to be defined: the potential $V(\phi)$. $V(\phi)$ shall be constructed in such a way that it features two minima \--- one for each of the two phases. This ensures that exactly these two phases (e.g., liquid and solid) are the most stable states. Furthermore, the potential needs to be symmetric. If we expand $V(\phi)$ in a Taylor series, and break off after the first two terms, we will find the definition  \cite{Emmerich2002}:

\begin{equation}
V(\phi) = \frac{V_0}{4} (\phi^2 - 1 )^2.
\end{equation} 
 
Here, it is assumed that the two minima have an energy level of $0$. In order to be able to also represent inhomogeneous states, the potential also has to depend on the gradient of $\phi$. The lowest order terms invariant under both rotation and translation are either $\nabla^2 \phi$ or $(\nabla \phi)^2$. The contribution of $\phi$ to the energy can then be expressed as:

\begin{equation}
P = \dots + \frac{\epsilon}{2} (\nabla \phi)^2 + V(\phi) \,,
\end{equation}

\noindent with $\epsilon$ as a constant. This formulation is termed the Cahn-Hilliard potential. 

\begin{figure}
\centering
\includegraphics[scale=0.3]{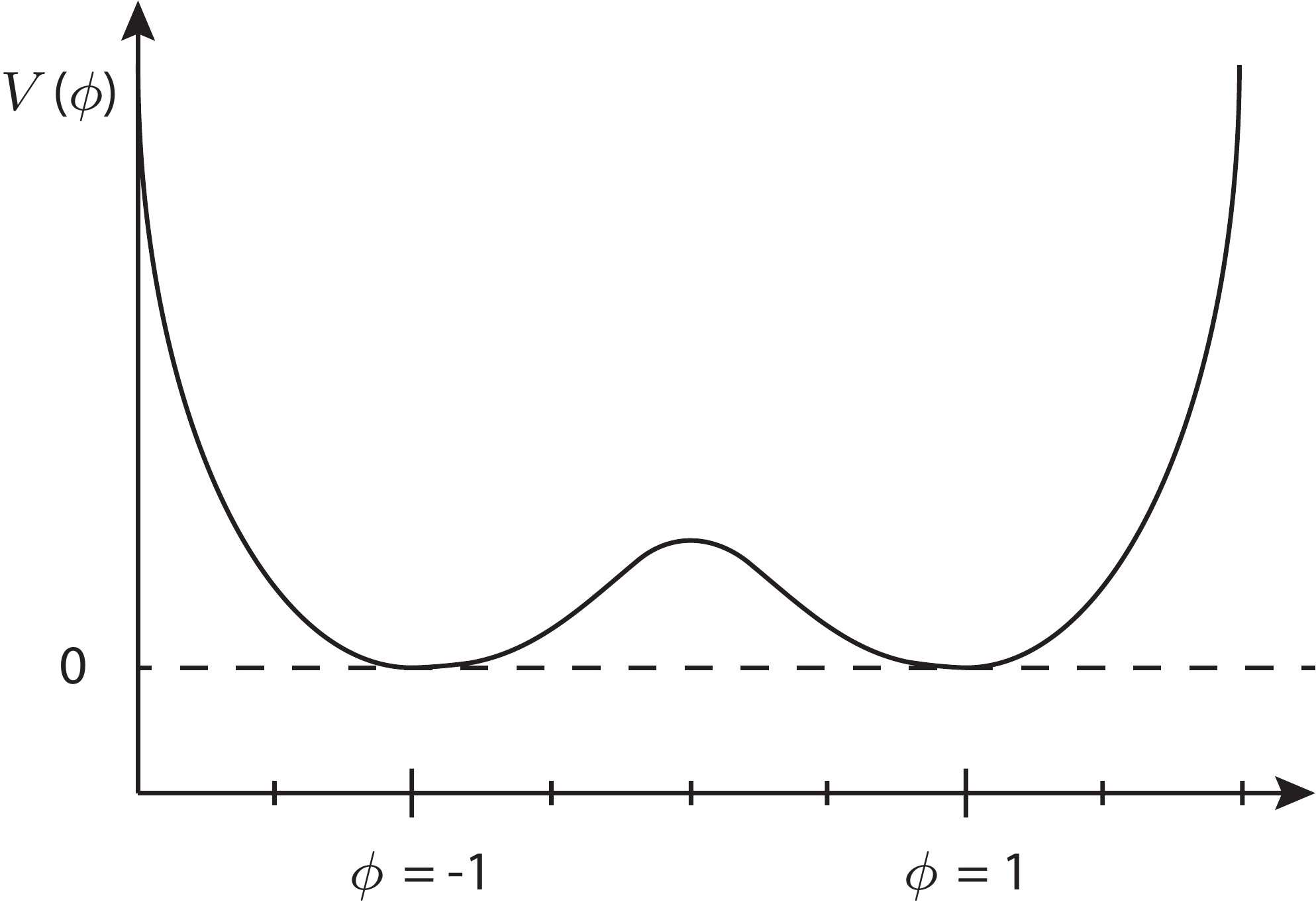}
\caption{Double-well potential: A potential depending on the phase-field variable $\phi$ is defined, which fulfills the basic criteria of symmetry and two minima.}
\label{fig:double-well}
\end{figure}

Again using variational principles, evolution equations can be derived for the phase field.

In non-equilibrium state, two important equations that can be developed from the Gibbs free-energy functional are the Cahn-Hilliard and the Allen-Cahn equation. The latter is a reaction-diffusion equation and in general not of importance for flow problems. The Cahn-Hilliard equation reads \cite{Boettinger2002}

\begin{align} \label{eqn:Cahn-Hillard}
\frac{\partial C}{\partial t} = \nabla \cdot \left( M_C \nabla \left( \frac{\partial f}{\partial C} - \varepsilon_C^2 \Delta C \right) \right) \,.
\end{align}

In the above equation, the unknown quantity is the concentration field $C$ of one of the two quantities present in the simulation (e.g., fluid A and fluid B or fluid and solid). C can be interpreted as equivalent to $\phi$. $f(C, \dots)$ denotes a free energy. It depends on the concentration, but also on other independent state variables such as possibly the temperature. $M_C$ is the mobility related to the solute diffusion coefficient and $\varepsilon_C$ is a constant related to the length scale of the diffuse interface. Note that the Cahn-Hilliard equation contains a fourth derivative of the concentration field; a challenge to many discretization methods.

The most important areas of application of the phase-field method are solidification processes, especially including phase segregation in binary alloys, and dentritic growth \cite{Echebarria2004,Gonzalez2013,Singer2008,Siquieri2011,Steinbach2009}. It has also been applied to other scenarios, such as spinodal decomposition \--- the unmixing of a mixture of liquids or solids \--- \cite{Axelsson2013}, surfactant transport \cite{Aland2010,Engblom2013}, or even topology optimization \cite{Dede2012}.  

In \cite{Liu2010}, isogeometric analysis (cf. Section~\ref{sec:IGA}) is utilized to analyze liquid-vapor phase transition phenomena modelled through the Navier-Stokes-Korteweg equations, which can be categorized as a phase-field model; in \cite{Gomez2008} the same approach is used for the solution of the Cahn-Hilliard equations in the context of spinodal decomposition. For both applications, the isogeometric analysis approach profits from the arbitrarily high continuity of the finite-element basis functions across element interfaces. For the first time, this allows for a standard use of the finite element method not augmented by mixed methods or discontinuous or continuous/discontinuous Galerkin methods \cite{Gomez2008}. 

The discussion of the phase-field method concludes the topic of interface-capturing approaches. In the next section, interface tracking will be introduced.

\section{The Interface-Tracking Approach} \label{s-interfacetracking} 
 
In interface tracking, the boundary or interface is explicitly resolved by the computational mesh. The necessity of accurately following the deforming boundary of the domain requires some degree of Lagrangian description to be present in the formulation, at least in the vicinity of the boundary. Since the use of a Lagrangian viewpoint is unavoidable, some researchers considered fully Lagrangian formulation for the entire domain~\cite{Hayashi1991,Okamoto90a,Radovitzky1998,Hirt70a,Okamoto90a}. The simplicity of such solutions is offset by several difficulties. The flow field in the interior of the domain may exhibit a large amount of circulation or shear, which are often not directly related to the motion of the boundary. The fully Lagrangian approach then results in unnecessary mesh distortion and invalid meshes in such interior regions, even for moderate or null displacements of the boundary. However, internal circulation and shear are handled without difficulty by an Eulerian description.\\
This provided a strong motivation for the development of methods capable of combining the Lagrangian and Eulerian approaches in the same domain. The ALE (Arbitrary Lagrangian-Eulerian)  formulation was initially stated in the finite difference context by Hirt~\cite{Hirt74a}, and later adopted also in the finite element community---see~\cite{Donea2003,Belytschko1980,Hughes81a,Huerta88a} and references contained therein. The basic idea of ALE is to work with a deforming mesh, as in the Lagrangian method, but to decouple the mesh deformation from the displacement of the fluid particles \--- or in other words, the mesh velocity is no longer equal to the fluid velocity, but can be chosen arbitrarily. The usual choice in ALE is to track boundary and interface nodes in a Lagrangian manner (or a manner that comes very close), while all other mesh nodes are adjusted solely with the intention of preserving mesh quality as best as possible. The ALE approach requires that the mesh velocity enters the momentum equation~\eqref{NavierStokes1} explicitly, thus modifying the original flow equations, specifically the convection term. Furthermore, this equation now needs to be written over a reference domain that may or may not coincide with either material or spatial domains. This added complexity is one of the main drawbacks of ALE. Nevertheless, it is used actively in the area of fluid flow computation, e.g., in~\cite{Chippada1994,Ganesan2006,Ganesan2012,Wang2010}. The ALE philosophy found another expression in the family of space-time finite element methods, which are described further below.

\subsection{Deforming spatial domain space-time finite element methods}

The classical choice of space and time discretization within the finite element community is to use finite elements (FE) in the spatial domain, but finite differences (FD) in the time direction. One alternative are the so called space-time finite element methods, which use finite element approximations for both space and time. Leaving the spatial discretization unaltered as compared to the classical approach, an additional coordinate direction \--- the time $t$ \--- is introduced for both the computational domain as a whole, but also for each finite element. As a consequence, all volume integrals in Equation~\eqref{eqn:weakform} are now integrals not only over the spatial domain $\Omega$, but over the space-time domain $Q$ of dimension $nsd+1$. The shape functions $N_i$, which interpolate the unknown functions, are equipped with an additional dependency on the time:   $N_i({\bf x},t)$. 
 
\begin{figure}[h]
\centering
\includegraphics[scale=0.35]{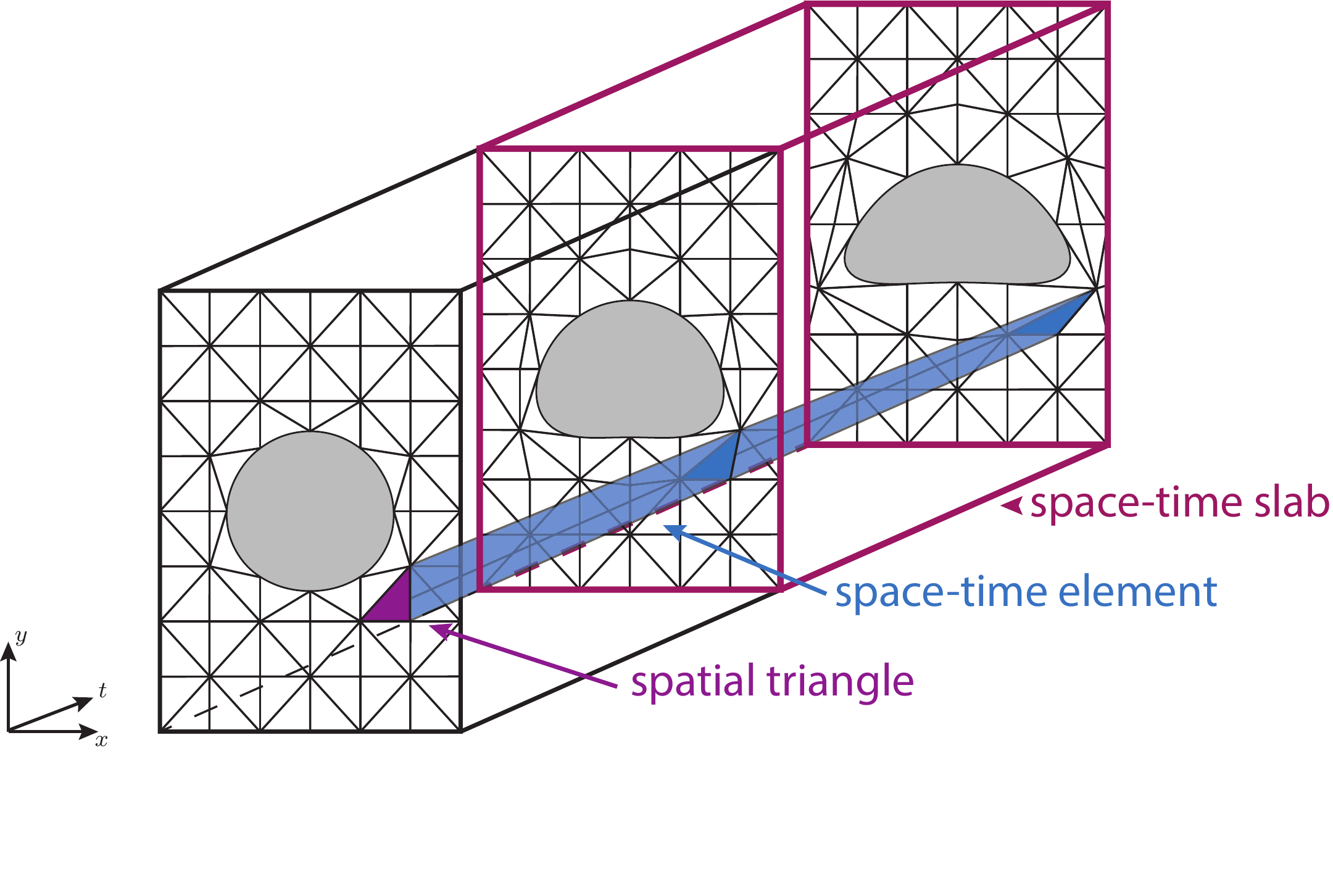}
\caption{Illustration of two successive space-time slabs $\hat{\Omega}$. From the lower time level to the upper time level, the mesh can deform. In blue, one space-time element is exemplified. It is triangular in space and forms a prism in space-time. }
\label{fig:spacetimeslab}
\end{figure}

Although it is in principle possible to employ the space-time discretization with an interpolation, which is continuous-in-time, the modern space-time approach usually relies on discontinuous in time discretization in the spirit of the \textit{Discontinuous Galerkin} method \cite{Donea2003,Hughes88a,Hansbo90a}. In order to avoid unmanageable number of degrees of freedom that must be solved for at any given time, the space-time approach is typically applied to subsets of the temporal domain called space-time slabs (cf. Figure~\ref{fig:spacetimeslab}), which are roughly comparable to time steps in the FE/FD approach. Due to the discontinuous in time interpolation, each node stores two values of every unknown per point in time. In order to enforce a relation between these two values, the continuity of the solution in temporal direction is imposed weakly by adding a jump term to the variational form \eqref{eqn:weakform}.      

The space-time method offers the unique opportunity to write the variational form directly over a deforming domain. This concept was pioneered by Jamet and Bonnerot~\cite{Jamet75a,Bonnerot77a,Frederiksen81a,Lynch80a,Jamet78a}. It was later extended concurrently by Tezduyar et al.~\cite{Tezduyar92a} in the form of the Deformable-Spatial-Domain/Stabilized Space-Time (DSD/SST) method, also applied in \cite{Behr92b,Mittal91a}, and  by Hansbo~\cite{Hansbo92a} with a particular mesh moving scheme conforming to the characteristic streamline diffusion ideology. Note that in the deforming domain space-time approach, there is no need to explicitly include the mesh velocity into the variational form, as it would be mandatory for ALE.  

In most space-time implementations to date, the meshes for the space-time slabs are simply extruded in the temporal direction from a spatial mesh, resulting in reference element domains that are always Cartesian products of spatial and temporal reference domains. Such an approach is best described as semi-structured (unstructured in space, structured in time) and does not leave the option of increasing temporal refinement in portions of the domain. In such a case, the space-time slab exactly corresponds to a time step of a semi-discrete procedure. In fact, many stencils of the FE/FD methods may be re-derived by using the semi-structured space-time approach with appropriate weighting and interpolation functions.

Note that the applicability of the space-time approach is of course not limited to the deforming domains in an interface-tracking context. Extensions of the XFEM concept to space-time~\cite{Chessa2004a,Pasenow2012a,Lehrenfeld2014a}, the monolithic space-time fluid-structure-interaction solvers~\cite{Huebner2004a}, or shape optimization~\cite{Elgeti2012d,Pauli2012} are just a few examples.

\subsection{The mesh deformation}

In interface tracking, the mesh is continuously adapted to the flow. Procedures for mesh deformation are described in the following, categorized into boundary deformation and displacement of the inner nodes. 

\subsubsection{Displacement of the boundaries} \label{sec:boundarydisplacement} \label{sec:boundarydeformation}

Section \ref{sec:interfacetracking1} describes the standard choice for the motion of the domain interface/boundary in the interface tracking context: a full Lagrangian deformation as specified by Equation \eqref{eqn:Lagrangian}. This translates to a deformation, where the mesh velocity ${\bf v}$ is chosen exactly equal to the fluid velocity ${\bf u}$ :

\begin{align}
  {\bf v}({\bf x}) = {\bf u} ({\bf x}) \,.
\label{eqn:fullLagrange}
\end{align}

However, this is not the only possible choice. Physically, the deformation of the interface/boundary mimics a no-penetration boundary condition, meaning that no fluid is allowed to cross the interface \--- this condition will always be fulfilled if the interface moves in accordance with the fluid as in Equation~\eqref{eqn:fullLagrange}. However, it will also be fulfilled by any choice where the mass flux through the boundary, $\dot{m}$, is $0$. The mass flux $\dot{m}$ through the interface can be computed as \cite{Ferziger1999}: 

\begin{align}
\dot{m} = \int_{\Gamma^{int}} \rho( {\bf u} - {\bf v} ) {\bf n}\;d{\bf x} \overset{!}{=} 0 \,.
\end{align}

\begin{figure}
\center
\psfrag{a}[rb][lt]{${\bf v}_1$}
\psfrag{b}[rb][lt]{${\bf v}_2$}
\psfrag{c}{${\bf n}$}
\psfrag{d}{${\bf u}$}
\psfrag{e}{${\bf v}_3$}
\includegraphics[width=5.5cm]{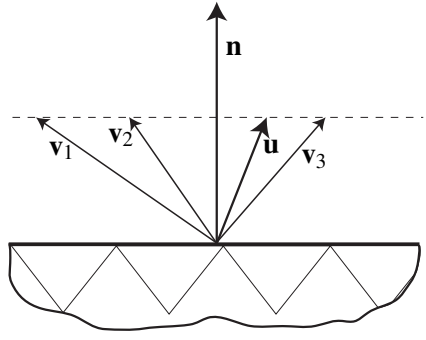}
\caption{Possible displacement directions satisfying the kinematic boundary condition.}
\label{kbc}
\end{figure}

Consequently, all choices for ${\bf v}$ that comply with the kinematic boundary condition

\begin{align}
{\bf v}({\bf x}) \cdot {\bf n}({\bf x}) = {\bf u}({\bf x}) \cdot {\bf n}({\bf x})
\label{eqn:kbc}
\end{align}

\noindent are valid (cf. Figure \ref{kbc}). It can be deduced from Equation~\eqref{eqn:kbc}, that tangential node movement remains irrelevant to the interface/boundary shape, but may significantly influence the mesh quality. Particularly applications with large tangential velocities \--- such as die swell in extrusion processes or rising droplets \---  profit tremendously if the tangential velocity component is modified or even suppressed when computing the interface/boundary velocity ${\bf v}$.

Note, however, that the assumptions made by the kinematic boundary condition are only valid strictly for the analytical geometry of the interface/boundary. This is particularly due to the computation of the normal vector, which is usually ambiguous at the  finite element nodes (usual finite element discretizations employ shape functions, which offer only $C^0$ continuity across element boundaries). Still, the finite element nodes are exactly the points of the interface/boundary where the mesh velocity needs to be defined. It is therefore inherent to the nature of the discretization that Equation~\eqref{eqn:kbc} can only be fulfilled approximately by any choice other than Equation~\eqref{eqn:fullLagrange}. \cite{Stavrev2015}

 Nevertheless, Behr \cite{Behr92} indicates several possibilities for boundary displacement resulting from Equation~\eqref{eqn:kbc}: displacement with the normal velocity component with $ {\bf v} = ({\bf u} \cdot {\bf n}){\bf n} $ and displacement only in a specific coordinate direction ( e.g., $y$-direction), i.e., $ {\bf v} = \frac{({\bf u} \cdot {\bf n}) {\bf e}_y}{ {\bf n} \cdot {\bf e}_y} $. Note here that a direct projection onto   the unit vector $ {\bf e}_y$ in the coordinate direction perpendicular to the main flow direction  is not possible, as this would result in a displacement in the correct direction, but not with the correct magnitude to satisfy the kinematic boundary condition.\\
In cases with very extensive deformations, even this method can lead to strongly distorted elements. As a remedy, it is possible to allow the nodes to slip along the tangential direction. The amount of tangential movement is then not connected to the fluid velocity, but determined according to the same method responsible for the inner nodes, which will be described in the next section.  

\subsubsection{Retaining mesh quality: mesh update for inner nodes }

As described in the previous section, both the ALE approach and the deforming domain space-time approach prescribe only the mesh displacement on (or possibly in the vicinity of) the boundary/interface. All other nodes have to be adapted to this displacement in a way that retains good mesh quality. Wall \cite{Wall1999} identifies two possible categories for the treatment of the inner nodes: methods that rely on explicit functions, which have been defined a priori, and methods that require the solution of their own system of equations.  

Examples of the first category can be found in \cite{Probst2013,Kjellgren1998}, where the boundary deformation is distributed smoothly throughout the entire domain using an interpolation kernel. 

In the second category, probably the most popular option is the Elastic Mesh Update Method (EMUM) \cite{Johnson94a} and comparable approaches (e.g., \cite{Masud97a,Onate2001a}). In this method, the computational mesh is treated as an elastic body reacting to the boundary deformation applied to it. For this purpose, the linear elasticity equation is solved for the mesh displacement $\bm{\upsilon}$, which relates to the mesh velocity ${\bf v}$ as $\bm{\upsilon} = {\bf v} \Delta t$:

\begin{align}
\label{eq:sopt-emum}
\nabla \cdot \mbox{${\bm{\sigma}}$}_{\mathrm{mesh}} &= 0 \,, \\
\mbox{${\bm{\sigma}}$}_{\mathrm{mesh}}({\bm{\upsilon}}) &= \lambda_{mesh} \left( {\mathrm{tr}}  {\bm{\varepsilon}}_{\mathrm{mesh}} ({\bm{\upsilon}})\right){\bf I} + 2 \mu_{mesh} {\bm{\varepsilon}}_{\mathrm{mesh}}({\bm{\upsilon}}), \\
{\bm{\varepsilon}}_{\mathrm{mesh}} ({\bm{\upsilon}}) &= \frac{1}{2} \left(\nabla {\bm{\upsilon}} + (\nabla {\bm{\upsilon}})^T\right).
\end{align}
In structural mechanics, \( \lambda_{mesh} \) and \( \mu_{mesh} \) are the Lam\'e-parameters that are needed to define the stiffness tensor for isotropic and linear elastic materials. In the context of the elastic mesh update, these parameters are prescribed for each element in order to control its stiffness. This is used to increase the stiffness of the smaller elements compared to the larger elements in order to allow for larger deformation before the mesh becomes invalid. The choice of Lam\'e-parameters can differ significantly between the approaches. 

\begin{figure}[htbp]
\center
\subfigure[Undeformed mesh]{
\includegraphics[width=4.5cm]{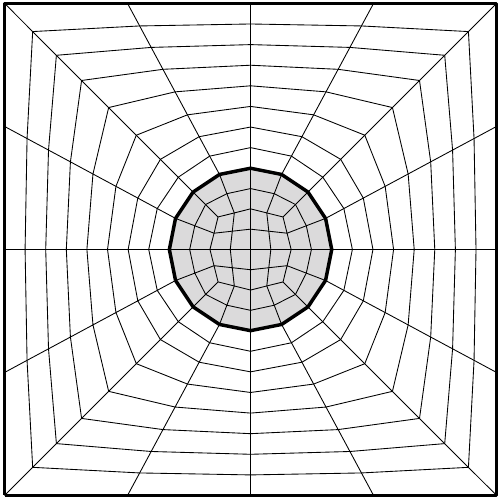}
}
\hfil
\subfigure[Elastically deformed mesh]{
\includegraphics[width=4.5cm]{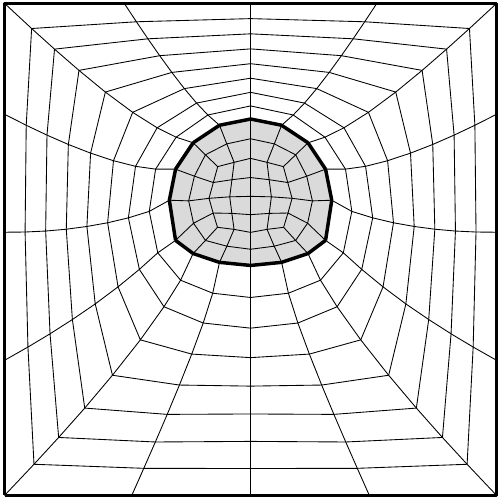}
}
\caption{Elastic Mesh Update Method: An initially circular interface (a) is subject to a non-uniform motion. (b) shows the ideal elastic response of the entire mesh with respect to the interface movement.}
\label{fig:EMUM}
\end{figure}

As boundary conditions, the displacement of the boundary/interface nodes, which is prescribed due to the rules given in Section~\ref{sec:boundarydeformation}, is imposed:

\begin{align}
\bm{\upsilon} = \hat{\bm{\upsilon}} \quad {\mathrm{on}} \ \Gamma_{\bm{\upsilon}}(t) \,.
\end{align} 

Another possible boundary condition, which, in the context of mesh deformation, is usually used in addition to the condition above, would be:

\begin{align}
{\bf n} \cdot \mbox{${\bm{\sigma}}$}_{\mathrm{mesh}} = {\bf h}_{\mathrm{mesh}} \quad {\mathrm{on}}\ \Gamma_{{\bf h}_{\mathrm{mesh}}}(t) \,.
\end{align}

This condition lets the nodes slip along certain boundaries. $\Gamma_{\bm{\upsilon}}(t)$ and $ \Gamma_{{\bf h}_{\mathrm{mesh}}}(t) $ are no longer distinct subsets of $\Gamma^{int}(t)$, but can be combined for different degrees of freedom. The calculated displacement $\bm{\upsilon}$ can now be used to update the nodal coordinates in all phases:

\begin{align}
{\bf x}^{n+1} = {\bf x}^n + {\bm{\upsilon}} \,.
\end{align} 

As a possible alternative, a simple Laplace equation can be solved for the displacement instead of Equation \eqref{eq:sopt-emum}, which diffuses the displacement throughout the domain, but gives less control over the individual elements. 

Also in the second category, but from a completely different view point, a method where a functional associated with the mesh distortion is minimized is proposed in \cite{Lopez2008}. The approach presented in \cite{Knupp2002} is based on a series of optimization steps with regard to the mesh quality. Loub\`ere et al. \cite{Loubere2010} apply a method, where all nodes are displaced in a fully Lagrangian fashion, but each time-step contains a reconnection step, where a new connectivity is established based on the existing nodal positions. 

As a final option, remeshing at every time-step can be applied, as for example in \cite{Onate2004} in the context of the particle finite element method (PFEM). Strictly speaking PFEM is a particle method, with marker particles carrying the information of the physical properties. What relates this method to the interface-tracking approach however is that the particles are connected by a Lagrangian mesh, where the marker particles simultaneously act as finite element nodes, on which the flow equations are evaluated.

\subsection{Normal vector and curvature approximation}

In principle, the explicit description of the interface/boundary allows the immediate evaluation of normals and curvature within the elements, as long as the shape functions can be differentiated sufficiently often. If the values for the normal vector or the curvature are needed directly at element nodes, they can be obtained by a weighted average of the values within neighbouring elements. One example for defining this average can be found in \cite{Engelman82a}.

In cases, where the geometric interpolation cannot provide the necessary geometric information \--- i.e., mainly in the case of linear interpolation \---, one may either resort to the Laplace-Beltrami technique (cf. Section \ref{sec:Laplace-Beltrami}) or introduce an additional geometrical entity from which the desired values can be obtained. Here, Mier-Torrecilla et al. \cite{Mier2011} propose the use of the osculating circle to the curve at each grid point. In the implementation, the osculating circle is identified as the circle with shoulder points at the respective node and its two neighbours. Using the radius $r$ of the osculating circle, the term $\kappa {\bf n}$ can be evaluated as $\frac{r}{|r|^2}$. The use of spline interpolations for this purpose is proposed in \cite{Ganesan2006} with the use of piecewise defined cubic splines, and in \cite{Elgeti2010} with the use of a single NURBS curve or surface to represent the boundary. 

The utilization of splines to obtain geometric quantities in interface tracking naturally leads to methods, which integrate Computer-Aided-Design (CAD) information into the finite element analysis.

\subsection{Recent advancement in CAD-integration for the finite element method} \label{s-NURBS} 

Geometries for engineering applications are generated using Computer-Aided-Design (CAD) systems: the CAD model is what we assume to be the real geometry. As of now, all major CAD systems share one common basis for geometry representation: Non-Uniform Rational B-Splines (NURBS). They provide a  a common standard to exchange  geometry information. Similar to the close relation between level-set methods and XFEM, an interconnection between CAD-related numerical methods and interface tracking can be identified. (Note however that in a few cases, the level-set method has also been combined with spline descriptions, as for example in \cite{Bernard2008,Fuchs2007}.) 

In the classical discretization methods, the exact NURBS geometry is lost as a finite element mesh, which is in general only an  approximation of the geometry, is generated. In this section, we will explain the structure of NURBS and then introduce two methods that make use of the NURBS format in order to integrate the exact geometry into the finite element method, thus avoiding the approximation character of the finite element mesh. 

\subsubsection{Geometry representation using splines}

We will start out with the most basic NURBS structure: the NURBS curve. From the curve definition, it will later be straightforward to derive the definitions of higher-dimensional objects. This section is based on \cite{Rogers2001,Piegl97}.
NURBS are an example of a parametric geometry representation. In this type of representation, a local parameter (e.g., $\theta$)  is defined along the curve. Usually, we use $\theta=0$ at one end of the curve and $\theta=1$ at the other end. The global coordinates of the curve $(x,y,z)$ are then each defined as functions of $\theta$:

\begin{align}
x = x(\theta), \quad y = y(\theta), \quad z = z(\theta) \,.
\end{align}

Important specifications during the historical development of NURBS were to generate a curve 

\begin{itemize}
\item[(1)] whose smoothness is completely under user control,
\item[(2)] which occupies a predefined space,
\item[(3)] which allows for local shape control, and
\item[(4)] which is able to represent both free-forms and analytical shapes.
\end{itemize}

\begin{figure}[htbp]
\centering
\includegraphics[scale=1.0]{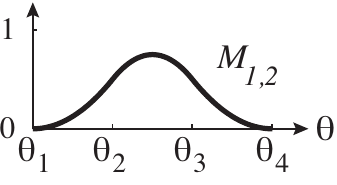}
\caption{Example of a B-spline basis function with polynomial degree $2$. $M_{1,2}$ is non-zero between $\theta_1$ and $\theta_4$. }
\label{fig:basisexample}
\end{figure}

In the NURBS context, the space occupied by the curve is defined using control points, which are connected to a control polygon. It is a property of the NURBS that it will always remain within the convex hull of the control points (Requirement (2)). By modifying the position of the individual control points, the shape of the curve will also be modified, as the control points guide the curve. The NURBS definition allows giving certain control points a larger influence on the curve compared to others, as each control point is assigned a weight. Another important feature of a NURBS is that each control point is only responsible for the shape of a minor portion of the curve (Requirement (3)). In the NURBS definition, this effect is realized through the B-spline basis functions $M$, the building blocks of the NURBS basis which will be introduced below. The B-spline basis functions are non-zero only in a sub-domain of the parameter space. Each control point $i$ is associated with exactly one basis function $M_{i,p}$ of polynomial degree $p$. It is the non-zero domain of this basis function that indicates which portion of the curve is influenced by the specific control point. This indication is given in terms of the local parameter $\theta$. A typical basis function for a NURBS is depicted in Figure~\ref{fig:basisexample}. The control point associated with this basis function will only influence curve points associated with $\theta_1 < \theta < \theta_4$. The values of $\theta$, where the influence of a basis function begins and ends, are referred to as {\it knots}. 

The curve definition ${\bf C}$ of a NURBS curve is then

\begin{align}
{\bf C}(\theta) = \sum_{i=1}^n R_{i,p}(\theta) {\bf P}_i,
\label{eqn:NURBScurve}
\end{align}

where $R_i$ denote the NURBS basis functions and $ {\bf P}_i$ are the coordinates of the control point $i$ \--- which are always given in the global coordinate system $(x,y,z)$. $n$ indicates the total number of control points. As we can see, ${\bf C}$ is a mapping from the parametric coordinate $\theta$ to the global coordinate system $(x,y,z)$. Notice that a NURBS curve has only one local parameter, no matter if it exists in $\mathbb{R}^1$ or whether it is embedded into the $\mathbb{R}^2$ or the $\mathbb{R}^3$. In this respect, NURBS are an example of an immersion.

\begin{figure*}[Htbp]
\centering
\includegraphics[scale=1.0]{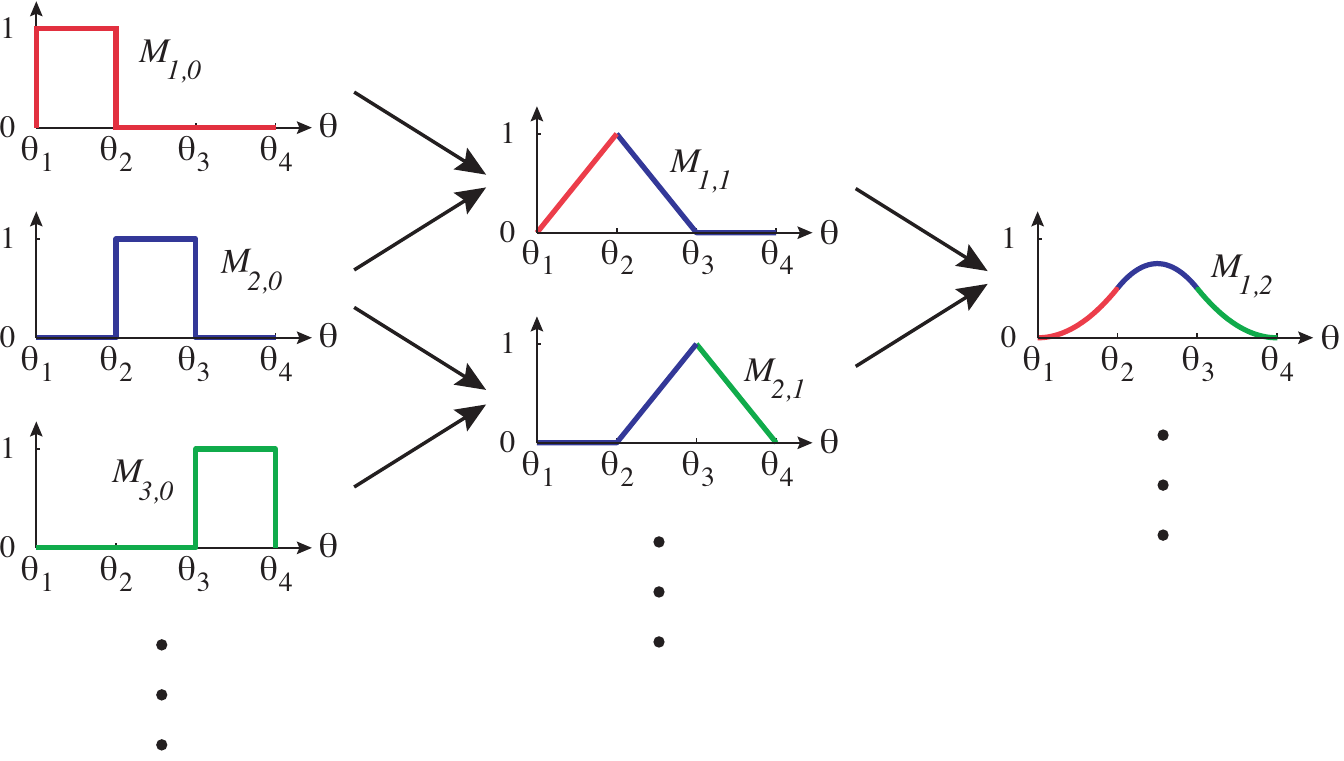}
\caption{Recursive definition of the B-spline basis functions: For each degree $p$, the total number of basis functions reduces by one. Each basis function of the higher degrees is built-up of $p+1$ constant basis functions, meaning that it is non-zero exactly between $\theta_i$ and $\theta_{i+p+1}$. In any interval $\theta_i < \theta < \theta_{i+1}$, exactly $p+1$ basis functions will be non-zero. }
\label{fig:basisstructure}
\end{figure*}

The NURBS basis functions $R$ emanate from the B-spline basis functions $M$. These in turn are usually defined recursively (cf. Figure~\ref{fig:basisstructure}), starting from the constant basis functions, using the Cox-deBoor relation:

\begin{align}
M_{i,0}(\theta) = \begin{cases} 1 \,, & \theta_i \leq \theta < \theta_{i+1} \,, \\ 0  \,,& otherwise \,,\end{cases}
\end{align}

\begin{align}
M_{i,p}(\theta)= \frac{\theta-\theta_i}{\theta_{i+p}-\theta_i}M_{i,p-1}(\theta) + \frac{\theta_{i+p+1}-\theta}{\theta_{i+p+1}-\theta_{i+1}}M_{i+1,p-1}(\theta) \,.
\label{eqn:Bsplinebasis}
\end{align}

$R$ is then defined as:

\begin{align}
R_{i,p}(\theta) = \frac{M_{i,p}(\theta) w_i}{\sum_{{\hat{i}}=1}^n M_{{\hat{i}},p}(\theta) w_{\hat{i}}}.
\label{NURBSbasis}
\end{align}

The rational structure of the basis is responsible for the fact that NURBS are able to represent all analytical shapes, including conic sections (Requirement (4)). The values $\theta_i, \theta_{i+p}$, etc. are the knot values, i.e., the indication of the influence domain of the basis functions. The knot values are collected in the so-called knot vector:  

\begin{align}
\Theta = [\underbrace{\theta_1, \theta_2}_{M_{1,0} \neq 0}, \underbrace{\theta_3,\theta_4}_{M_{3,0} \neq 0}, \dots] \,.
\end{align}

For the proper definition of the basis, it is essential that the knots are in a non-decreasing sequence. However, it is specifically possible to repeat knot values. With each repetition, the continuity of the basis is decreased. With regard to the curve, there are then two measures that influence the continuity: the continuity of the basis and the topology of the control points. If control points are placed on top of each other, the curve continuity is decreased, e.g., to represent sharp corners. If control points are aligned, the continuity can be increased. Both effects are a direct consequence of property that the NURBS curve is always contained within the convex hull of its control polygon. These two measures give complete control over the curve continuity (Requirement (1)). In general, the basis functions, and therefore the curve, are infinitely differentiable at all points, and $p-1$ times differentiable at the knots. 

The derivatives and the curvature of NURBS curves can be calculated analytically. From differential geometry, the curvature of a parametric curve $\bf{C}(\theta)$ can be calculated with \cite{Gray_1998}:
\begin{align}\label{eq:curvature}
	\kappa(\theta) = \frac{|\bm{C}'(\theta)\times\bm{C}''(\theta)|}{|\bm{C}'(\theta)|^{3}}.
\end{align}
$\bf{C}'(\theta)$ and $\bf{C}''(\theta)$ denote the first and second derivative of the curve with respect 
to parameter $\theta$. For the sake of brevity, $M_{i,p}(\theta)$ is written as $M_{i,p}$ and $\sum_{i=1}^{n}$ 
as $\sum_{i}$ in the following equations.
\begin{align}
	\bm{C}'(\theta) =& \frac{\sum_{i}M'_{i,p}w_{i}\boldsymbol{P}_{i}\sum_{i}M_{i,p}w_{i} - 
			\sum_{i}M_{i,p}w_{i}\boldsymbol{P}_{i} \sum_{i}M'_{i,p}w_{i}}
			{\left(\sum_{i}M_{i,p}w_{i}\right)^2},\\
	\begin{split}
	\bm{C}''(\theta) =& \frac{\sum_{i}M''_{i,K}w_{i}\bm{P}_{i}\sum_{i}M_{i,p}w_{i}
				-\sum_{i}M_{i,p}w_{i}\bm{P}_{i}\sum_{i}M''_{i,K}w_{i}}
				{\left(\sum_{i}M_{i,p}w_{i}\right)^2}\\
			&- \frac{2\sum_{i}M'_{i,p}w_{i}\left(\sum_{i}M'_{i,p}w_{i}\bm{P}_{i}
				\sum_{i}M_{i,p}w_{i}\right)}{\left(\sum_{i}M_{i,p}w_{i}\right)^3} \\
				&-\frac{2\sum_{i}M'_{i,p}w_{i}\left(\sum_{i}M_{i,p}w_{i}\bm{P}_{i}
				\sum_{i}M'_{i,p}w_{i}\right)}{\left(\sum_{i}M_{i,p}w_{i}\right)^3}.
	\end{split}
\end{align}
$M'_{i,p}$ and $M''_{i,p}$ are the first and second derivative of the basis functions with respect 
to the parameter $\theta$.

\begin{figure}[h]
\centering
\psfrag{a}{$\theta$}
\psfrag{b}{${\bf P}_6 = (-1,1,0)$}
\psfrag{c}{${\bf P}_7 = (-1,0,0)$}
\psfrag{d}{${\bf P}_8 = (-1,-1,0)$}
\psfrag{e}{${\bf P}_1 = {\bf P}_9 = (0,-1,0)$}
\psfrag{f}{${\bf P}_2 = (1,-1,0)$}
\psfrag{g}{${\bf P}_3 = (1,0,0)$}
\psfrag{h}{${\bf P}_4 = (1,1,0)$}
\psfrag{i}{${\bf P}_5 = (0,1,0)$}
\includegraphics[scale=0.3]{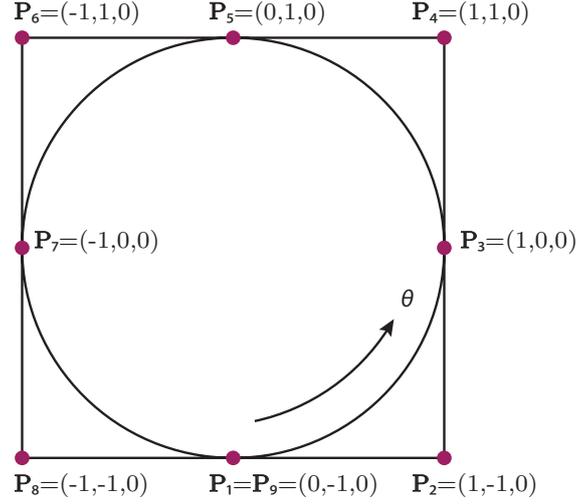}
\caption{Example of a NURBS curve: with 9 control points and a quadratic basis, it is possible to represent a full circle. The corresponding knot vector is 
$\Theta = [0\ 0\ 0\ 1/4\ 1/4\ 1/2\ 1/2\ 3/4\ 3/4\ 1\ 1\ 1 ] $. The weights have to be set to $w_1 = w_3 = w_5 = w_7 = w_9 = 1$ and $w_2=w_4=w_6=w_8= \frac{\sqrt{2}}{2}$.}
\label{fig:CurveExample}
\end{figure}

Figure~\ref{fig:CurveExample} gives an example of a NURBS curve representing a circle. A closed curve is obtained by overlapping the first and the last control point. Higher continuity at the connection point of the curve can be reached if $p$ control points at the ends are overlapped (${\bf P}_{n-p+1} = {\bf P}_1, \dots , {\bf P}_n = {\bf P}_{p}$).

If surfaces (or volumes) are to be described through NURBS, a second (or third) local parameter is introduced. The basis functions for the higher-dimensional objects are computed as a Cartesian product of the univariate basis functions $M_{i,p}(\theta)$ and $ L_{j,q}(\iota) $  defined in Equation~\eqref{eqn:Bsplinebasis}. For each local coordinate direction, an individual knot vector is required. Note that it is possible to choose different polynomial degrees of the basis for each local coordinate. The NURBS surface ${\bf{S}}$ is described as

\begin{align}
{\bf{S}}(\theta, \iota) = \sum_{i=1}^n \sum_{j=1}^m R_{i,j}^{p,q}(\theta, \iota) {\bf P}_{i,j} \,,
\label{eqn:NURBSsurface}
\end{align}

\noindent with 

\begin{align}
R_{i,j}^{p,q}(\theta, \iota) = \frac{M_{i,p}(\theta)L_{j,q}(\iota)w_{i,j}}{\sum_{\hat{i}=1}^n \sum_{\hat{j}=1}^m M_{\hat{i},p}(\theta)L_{\hat{j},q}(\iota)w_{\hat{i},\hat{j}}}.
\end{align}

\subsubsection{Isogeometric analysis} \label{sec:IGA}

Recall that the finite element method invokes the \textit{isoparametric concept}, as described in \cite{Hughes2000} or many other works on finite elements. The isoparametric concept signifies that the  element nodes are interpolated with the same shape functions as the unknown function, i.e., the solution of the partial differential equation. This approach instantly presents a mapping from reference coordinates to global coordinates, thus offering a simple means for evaluating all  components of the transformation theorem  necessary to compute the integrals of the weak form on a reference element instead of the global elements. In comparison to an integral evaluation in global coordinates, the approach using a reference element makes easy use of the local support of the shape functions and furthermore avoids the more tedious and computationally inefficient definition of the shape functions in global coordinates.

Cottrell, Bazilevs and Hughes used the exact same underlying approach when devising the concept of isogeometric analysis  \cite{Hughes2005,Cottrell2009,Bazilevs2010}. In the spirit of the isoparametric concept, it utilizes the NURBS basis functions in order to represent both the geometry and the unknown solution. For the geometry, the already presented definition \eqref{eqn:NURBScurve} or \eqref{eqn:NURBSsurface} is employed. The unknown function $u$ is then approximated in exactly the same fashion as: 

\begin{align}
u^h = \sum_{i=1}^n R_{i,p}(\theta) u_i \,.
\end{align}

\begin{figure}
\psfrag{a}{$u_1$}
\psfrag{b}{$u_2$}
\psfrag{c}{$u_3$}
\psfrag{d}{$u_4$}
\psfrag{e}{$u^h$}
\psfrag{f}{$x$}
\psfrag{g}{$u^h(x)$}
\psfrag{h}{${\bf P}_1$}
\psfrag{i}{${\bf P}_2$}
\psfrag{j}{${\bf P}_3$}
\psfrag{k}{${\bf P}_4$}
\psfrag{l}{$\theta$}
\psfrag{x}{$x$}
\centering
\includegraphics[scale=1.0]{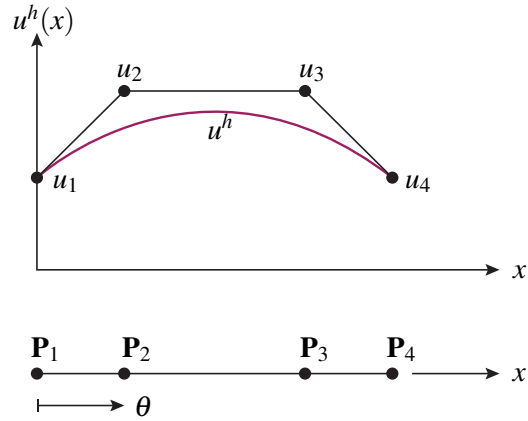}
\caption{Function representation in IGA: Consider a quadratic NURBS curve whose control points are aligned in such a way that it represents a line (lower part of the picture). The local parameter along the curve is $\theta$, and it is embedded into $\mathbb{R}^1$ with coordinate $x$. On this curve, a function $u$ is approximated ($u^h$). The function is represented through the control variables $u_1$ through $u_4$ which are interpolated using the NURBS basis function $R_{i,p}$. }
\label{fig:iga-function-representation}
\end{figure}

The values $u_i$ are the control variables. These are the unknown values solved for in the finite element code. Note that just like the control points for the geometry, the values of the control variables are in general not coinciding with the solution; see Figure~\ref{fig:iga-function-representation}.

Not only does this idea essentially provide the capability to perform finite element analysis directly on CAD models, but it also profits from superior approximation properties of NURBS as compared to the polynomial shape functions used in the standard finite element method.

\subsubsection{The NURBS-enhanced finite element method}

One major challenge still remaining for the isogeometric analysis described in the previous section is the generation of closed volume splines describing complex geometries. In the traditional CAD/CAM systems, complex geometries are usually described as a large number of individual surface patches, connected neither geometrically nor parametrically. The reason for this is that traditionally, the CAD system tries to  mimic classical manufacturing approaches (such as turning, drilling, or milling) starting out from a form similar to a semimanufactured product. This makes the handling of the CAD tool easily accessible to engineers. However, this also makes it very complicated to use for isogeometric analysis.

Sevilla, Fernandez-Mendes, and Huerta have proposed an alternative on middle ground between isogeometric analysis and standard finite elements: the NURBS-enhanced finite element method (NEFEM) \cite{Sevilla2008a,Sevilla2011e,Sevilla2011a,Stavrev2015}. In NEFEM the boundary of the geometry is represented using NURBS, thus leading to an exact representation, whereas for the interior, a standard finite element mesh, with all benefits of existing meshing algorithms, is utilized. Note that this leads to two different kinds of elements: (1) standard finite elements in the interior and (2) elements with one NURBS edge alongside the boundary. Usually, elements of category (1) will be in the vast majority, keeping the computation very efficient. Through the elements of category (2), the exact geometry is made available during the process of evaluating the integrals of the appropriate weak form (e.g., of the governing equations Equations~\eqref{NavierStokes1}--\eqref{NavierStokes2}): The integration domain $\Omega^h$ is no longer only an approximation of the real domain $\Omega$.  The availability of the exact geometry is however not resorted to for the interpolation of the unknown solution. The latter is approximated using standard Lagrangian shape functions both along the boundary and in the interior. As a consequence, NEFEM is no longer an isoparametric method. This is in contrast to IGA. 

\begin{figure}[htbp]
\center
\subfigure[Position of quadrature points]{
\includegraphics[width=3.0cm]{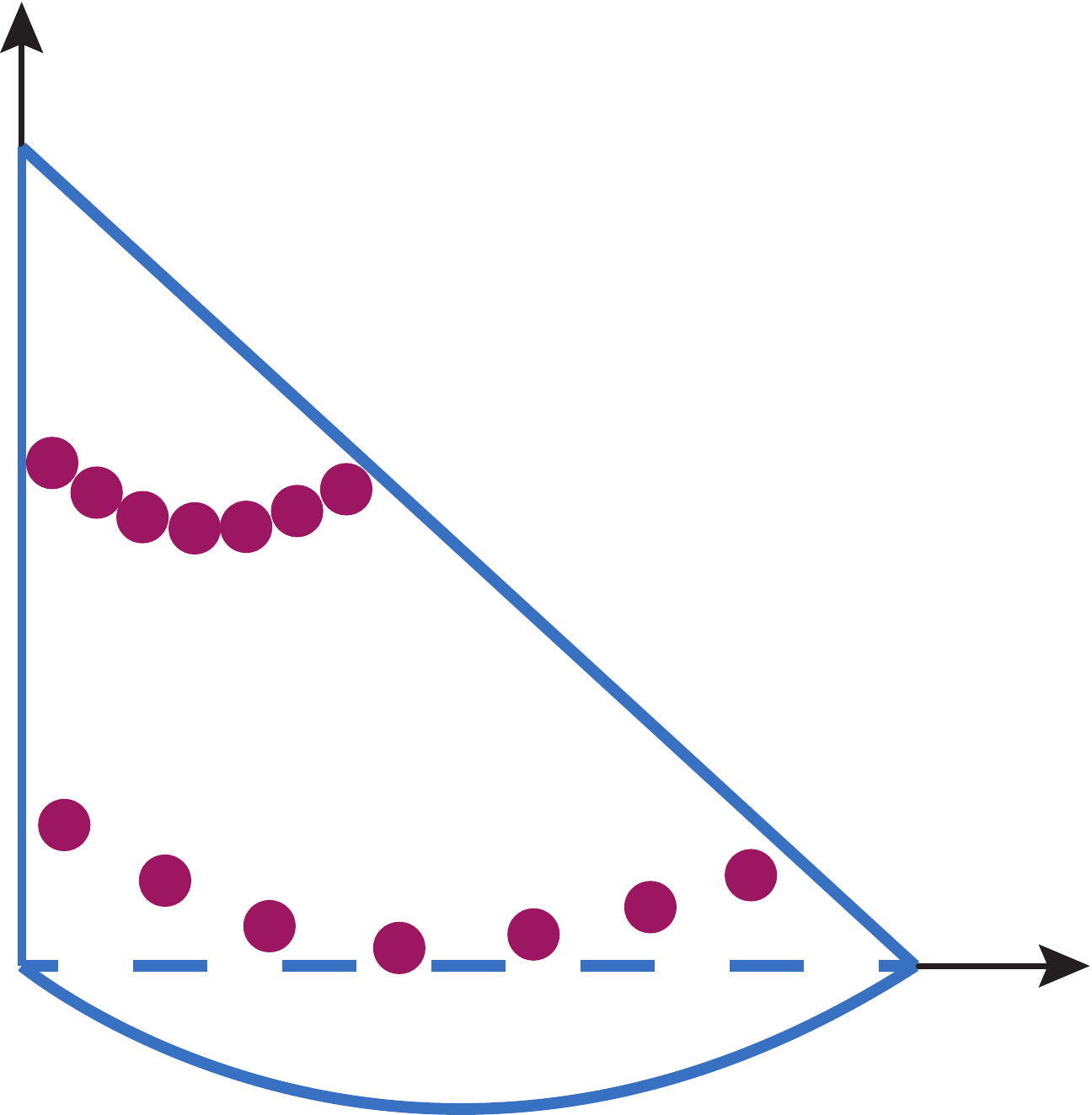}
}
\hfil
\subfigure[Example of a linear shape function]{
\label{resistance}
\includegraphics[width=3.2cm]{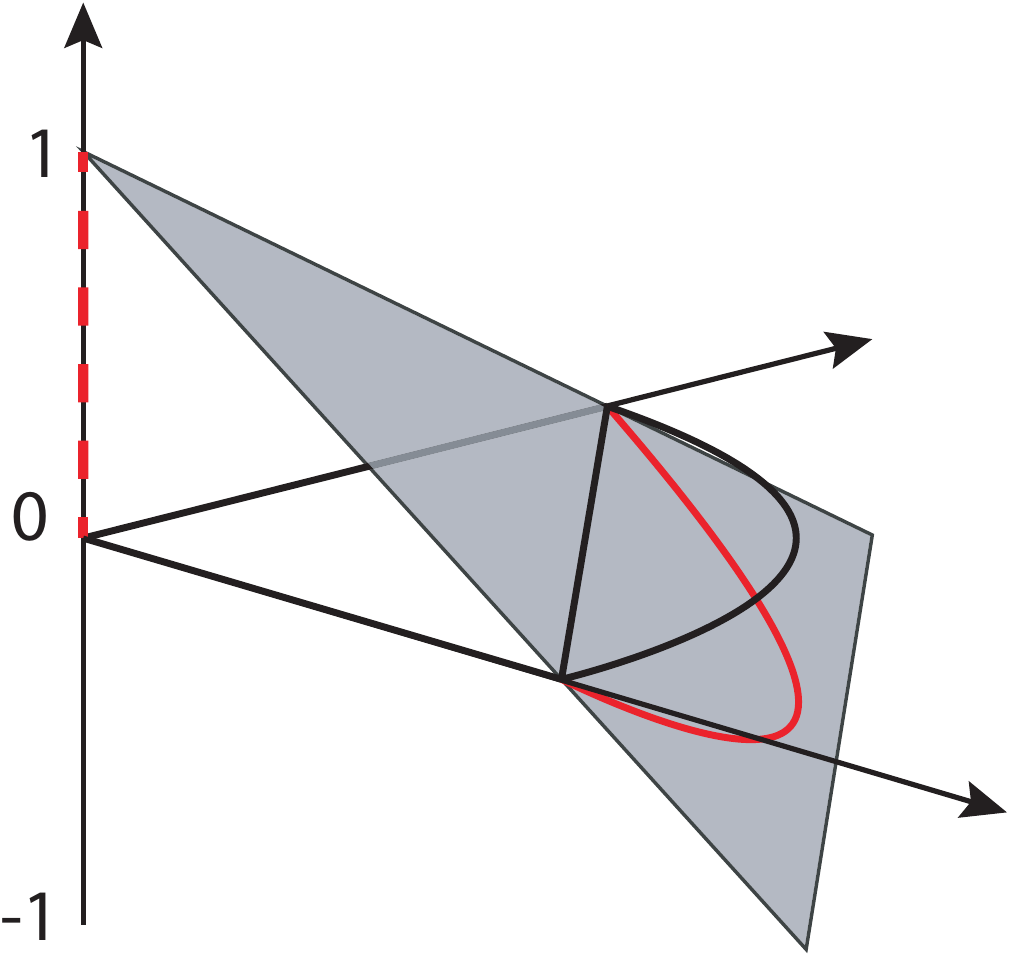}
}
\caption{(a) The NEFEM quadrature points are adapted to the curved triangle shape. (b) Example of a linear shape function in the NEFEM context. Negative values and values larger than $1$ may occur.}
\label{fig:NEFEM}
\end{figure}

The main advantage of NEFEM over standard finite elements is that the position of the integration points needed in the finite element method are determined from the curved NURBS geometry and not from the approximated geometry (cf. Figure~\ref{fig:NEFEM}). If we consider boundary integral, such as the integral arising from the surface tension, the integration points are evenly distributed on the curved geometry. Due to the continued approximation of the unknown function through Lagrange polynomials an oddity arises: negative shape function values may occur along the boundary (cf. Figure~\ref{fig:NEFEM}). Furthermore, the shape function connected to all nodes, also those not located on the NURBS edge, are non-zero along the NURBS edge. This has the consequence, that interior nodes contribute to boundary integrals \--- a highly unusual situation in the finite element method.     

This section concludes the introduction of the numerical methods. The next sections will concentrate on standard applications.
 \section{Bubbles and drops} \label{s-drops} 

The simulation of drops (liquid) or bubbles (gas) in a surrounding bulk phase is one of the classical benchmark cases in multi-phase flow. It addresses the issues of surface tension effects, moving interfaces, as well as the question of break-up and coalescence of phase domains (in this case, the drops).  The test case is motivated by the application of liquid-liquid or liquid-gas extraction columns \cite{Bertakis2010}. The relevant effects can however also be extended for example to the simulation of the alveoli in the lung.
 
 \subsection{Static drop inside a rectangular domain}
 
\begin{figure}[htbp]
\center
\subfigure[Computational domain]{
\includegraphics[width=3.0cm]{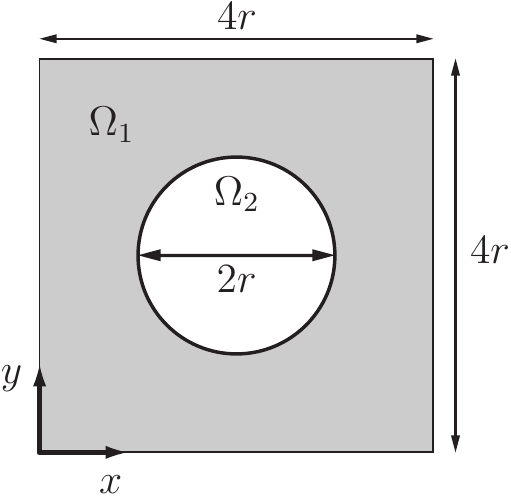}
\label{fig:staticdrop1a}
}
\hfil
\subfigure[Exact pressure solution]{
\includegraphics[width=4.3cm]{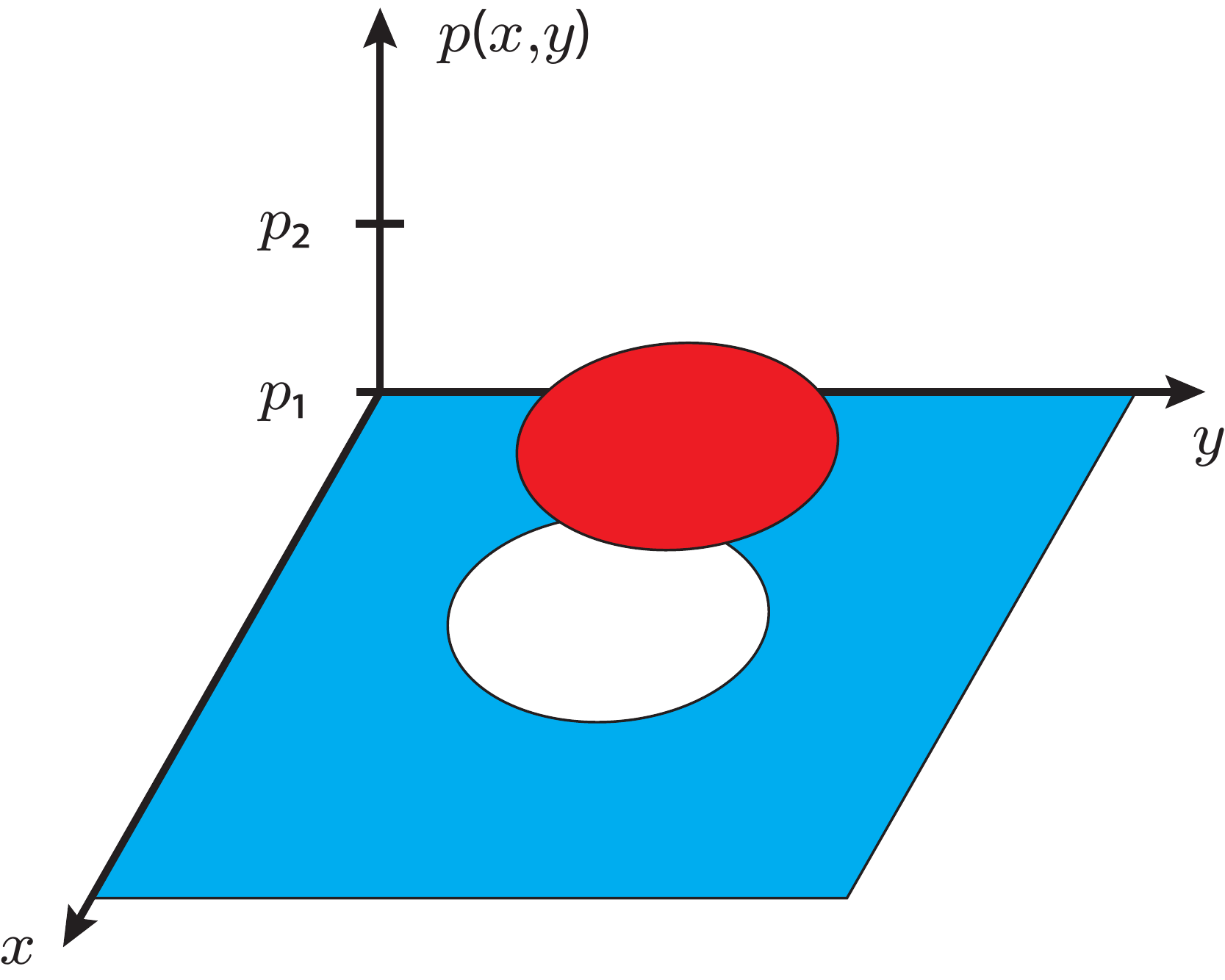}
\label{fig:staticdrop1b}
}
\caption{Static test case with circular interface. In (b), a zero reference pressure is assumed for $\Omega_1$. Furthermore, the parameters have been chosen to $\rho_1$=$\rho_2$= $1.0 kg/m^3$, $\mu_1$=$\mu_2$= $1.0 kg/m/s$, $\gamma = 1.0 kg/s^2$ and $r=0.5m$.}
\label{fig:staticdrop1}
\end{figure}

 The first test case regularly used in the context of drop simulations is a circular drop inside a rectangular domain (cf. Figure~\ref{fig:staticdrop1a}). The domain $\Omega_2$ is occupied by fluid 2, which is surrounded by another fluid 1 residing in the domain $\Omega_1$. The dimensions of the domains differ within the literature. What remains in common is that  the steady Stokes equations without external forces (i.e., in particular without gravity) are solved on the two domains (cf. Equation~\eqref{StatStokes_1}--\eqref{StatStokes_2}). Density and viscosity of the two fluids are usually chosen to be equal. However, the presence of surface tension with a surface-tension coefficient of $\gamma$ is assumed. Since the analytical value of the sum of principal curvatures is known for circles ($\kappa = \frac{1}{r}$) and spheres ($\kappa = \frac{2}{r}$), according to the Laplace-Young equations, the analytical solution is (cf. Section~\ref{sec:bc} ):

 \begin{align}
 {\bf u} ( {\bf x} ) &= {\bf 0} \,, \\
 p( {\bf x} ) &= \begin{cases} 0, & {\bf x} \in \Omega_1 \,,\\ \gamma \kappa, & {\bf x} \in \Omega_2 \,.\end{cases}
 \end{align}
  
With the analytical solution at hand, this test case is particularly valuable to demonstrate the capabilities of a chosen method to handle the jump of quantities across the interface as well as the evaluation of geometrical quantities of the interface. This is even more the case as Ganesan et al. \cite{Ganesan06a} point out that this solution is in general not represented well by a discretization scheme: the discrete velocities are not equal to ${\bf 0}$ \--- we see spurious velocities. \cite{Ganesan06a} identifies two causes: inadequate representation of the pressure solution (particularly the pressure jump) and an insufficient approximation of the curvature. They give the error bound on ${\bf u}^h$ as (in our own notation): 

\begin{align}
|u^h|_1 &\leq C \Big( \inf_{q^h \in Q^h} \parallel p - q^h \parallel_0 + \sup_{w^h \in V^h} \frac{|(\kappa^h, w^h \cdot {\bf n}) - (\kappa, w^h \cdot {\bf n}) | }{|w^h|_1} \Big) \,.
\end{align} 

\begin{figure}[h]
\center
\subfigure[Continuous interpolation.]{
\includegraphics[width=7cm]{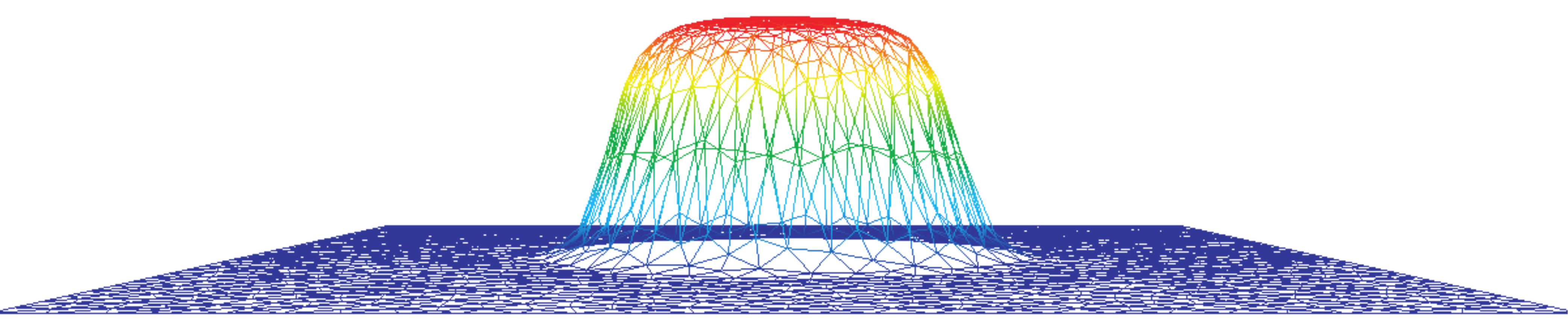}
\label{pjcont}
}
\hfill
\subfigure[Doubled interface nodes.]{
\includegraphics[width=7cm]{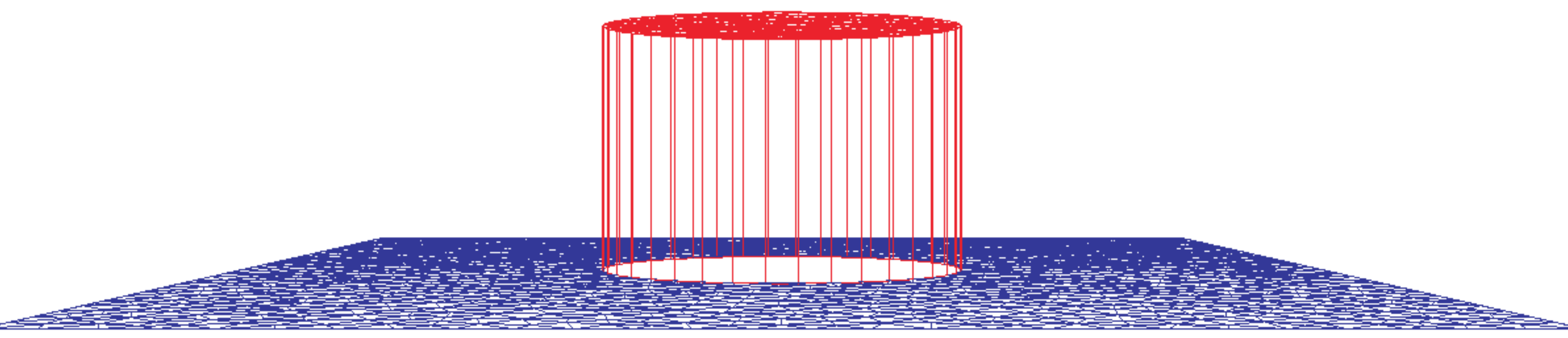}
\label{pfdiscont}
}
\vspace{-6pt}
\caption{Comparison of the pressure jump using continuous interpolation functions and interface nodes with doubled pressure degrees of freedom \cite{Elgeti2010}. \copyright Wiley-Blackwell}
\label{fig:pressurejump}
\end{figure}

From theoretical considerations, Ganesan et al. \cite{Ganesan06a} conclude that the use of continuous approximation functions in the case of a discontinuous solution is to always be avoided -- this holds true even if the discontinuities are resolved by the mesh. The jump is then spread over several elements, thus distorting the solution (cf. Figure \ref{pjcont}). Here, and also in Gross and Reusken \cite{Gross2010}, it is shown that in general, the error in the pressure approximation $\inf_{q^h \in Q^h} \parallel p - q^h \parallel_0$ is bounded by $c\sqrt{h}$, leading to very slow convergence. 

These theoretical results brought about a whole set of ideas aiming at numerically obtaining both an instant pressure jump and an improved curvature approximation.

\emph{Pressure jump}: In the area of level-set methods, the work concentrated on developing problem-suitable finite element pressure spaces: XFEM evolved as the main means of incorporating the pressure jump (cf. Section \ref{sec-ls-sharp-interfaces}). In  \cite{Gross2010,Gross2007,Reusken2008}, a Heaviside enrichment is employed for this purpose. In the case of continuous basis functions, convergence orders for the pressure of only $O(h^{3/2})$ in the $L^2$-norm were reached. For the velocity, the reported convergence order is $O(h^{1/2})$ in the $H^1$-norm. However, if the Heaviside-enrichment was employed in combination with an appropriate curvature approximation, a convergence order of $O(h^\alpha)$ with $\alpha \geq 1$ could be recovered.
Also in the area of level-set approximations, Ausas et al. \cite{Ausas2009} present an approach where the pressure space is locally enriched at those elements, which contain the interface. The corresponding elements are subdivided along the interface, and each part of the domain is then only influenced by pressure degrees of freedom on one side of the interface. The method does not introduce new degrees of freedom for the pressure. In \cite{Sousa2012}, the stability of this method is compared to the method of Coppola-Owen \cite{Coppola2005} already mentioned in Section~\ref{s-levelset}. 

Along similar lines, but in the context of the PFEM method, the pressure degrees of freedom at the interface are duplicated, once associated to $\Omega_1$ and once to $\Omega_2$, in \cite{Mier2011}. This leads to the possibility of an exact representation of the pressure jump even on coarse meshes. In the context of boundary-conforming meshes, the same idea has been employed in \cite{Elgeti2010}. The velocity values of the doubled nodes are coupled in the context of the boundary condition in Equation~\eqref{eq:contvel}, but the pressure is allowed to take on different values. A sample result is depicted in Figure~\ref{fig:pressurejump}.
The advantage of the local approaches is that they take features of the individual problem into account, e.g., the knowledge that the pressure jump will only ever occur at the interface. This makes the enrichment less general, but usually much more efficient and less detrimental to the numerical scheme (such as, for example, the ill-conditioned matrices XFEM can induce).

\emph{Curvature approximation}: Most approaches never compute the curvature of the interface, but resort to the Laplace-Beltrami technique or the continuum surface force technique as described in Section~\ref{s-challenges}  (among many others \cite{Barrett2013}). Among the rare explicit approaches, several researchers have attempted the use of splines as an interface representation, from which the curvature can then be derived analytically. This includes \cite{Ganesan_2007} where piecewise interpolated cubic splines provide an additional description of the discretized interface. In  \cite{Elgeti2010,Stavrev2015}, a single NURBS (cf. Section \ref{s-NURBS} ) represents the entire circle, thus providing an excellent means for evaluating the curvature at any given point of the interface. \cite{Stavrev2015} furthermore employs a space-time version of the NEFEM (cf. Section \ref{s-NURBS}). Here, a combination of P1P1 finite elements and a quadratic NURBS for the geometry description is employed. With this combination, the analytical solution can be obtained exactly even on a coarse mesh with only 8 FE nodes on the interface.

\subsection{Rising drops or bubbles}

The test case described in the previous section can be extended when independent material parameters ($\mu_{1,2}$ and $\rho_{1,2}$) as well as gravity are taken into account. Depending on the density relation, the drop in $\Omega_2$ will then either rise or fall in $\Omega_1$ due to buoyancy. In comparison to the previous test case, the set-up of this section needs to additionally consider the velocity kink at the interface, the jump in the pressure gradient, and most of all, a scheme for deforming domains.

Experimentally, the phenomenon of a rising drop is for example described in \cite{Clift1978}. The book defines three dimensionless numbers, any two of which completely characterize the flow regime. They involve the gravity $g$, the drop diameter $d$, viscosity $\mu$ and density $\rho$, as well as the surface-tension coefficient $\gamma$:

\begin{itemize}
\item The Reynolds number, relating inertia to friction forces,
\begin{align}
Re = \frac{\rho_2 \sqrt{g d} d}{\mu_2},
\end{align}
\item the E\"otv\"os number, relating buoyancy forces and surface tension,
\begin{align}
Eo = \frac{\rho_2 d^2g}{\gamma},
\end{align}
\item and the Morton number, relating viscous forces and the surface tension,
\begin{align}
Mo = \frac{g \mu_2^4}{\rho_2 \gamma^3} = \frac{Eo^3}{Re^4}.
\end{align}
\end{itemize}

Depending on these characteristic numbers, three regimes can be distinguished: spherical, ellipsoidal, and spherical cap. Spherical drops occur at low $Re$ and $Eo$. The ellipsoidal regime is located around relatively high $Re$ and intermediate $Eo$, while both high $Re$ and $Eo$ yield a drop in the spherical cap regime. 

Hysing et al. published benchmark computations for a single rising bubble \cite{Hysing2009}. It contains two benchmark cases, one in the ellipsoidal regime and one in the spherical cap regime, featuring thin filaments and break-up. The results of three different numerical methods are described; two use the finite element method in conjunction with level-set and one finite elements in an ALE setting. 

\begin{figure}
\centering
\includegraphics[scale=0.8]{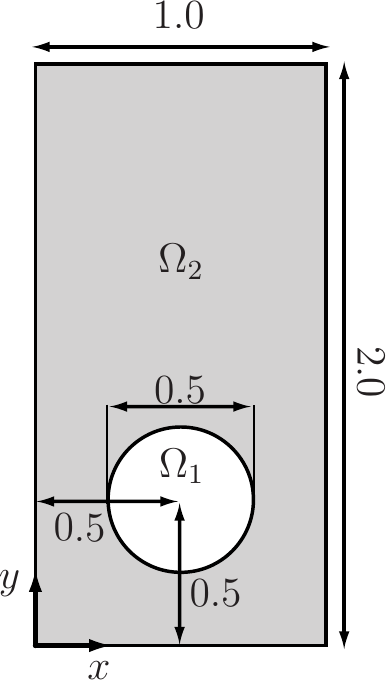}
\caption{Rising bubble: Initial computational domain.}
\label{fig:risingdropsetup}
\end{figure}

In \cite{Sauerland2013}, the first benchmark of \cite{Hysing2009} was repeated. Here, the surface tension effects dominate the bubble shape and no break-up should occur. The rectangular domain has the size of $1.0 \times 2.0 m$ with an initially circular bubble with diameter $d= 0.5 m$, as shown in Figure~\ref{fig:risingdropsetup}.

 The properties of the fluids are $\rho_1 = 100 kg/m^3$, $\rho_2=1000 kg/m^3$, $\mu_1 = 1 kg/m/s$, $\mu_2 = 10 kg/m/s$, $f_y=-g=-0.98 m/s^2$, and surface tension coefficient $\gamma = 24.5 kg/s^2$. Under these circumstances, the characteristic Reynolds and E\"otv\"os numbers can then be evaluated to:
 
 \begin{align}
Re = \frac{\rho_2 \sqrt{gd} d }{\mu_2} = 35,
 \end{align}
 
 \begin{align}
Eo = \frac{g \rho_2 d^2}{\gamma} = 10.
 \end{align}
 
\begin{figure}[htbp]
\center
\subfigure[Bubble shape]{
\label{fig:bubble-shape}
\includegraphics[width=5.0cm]{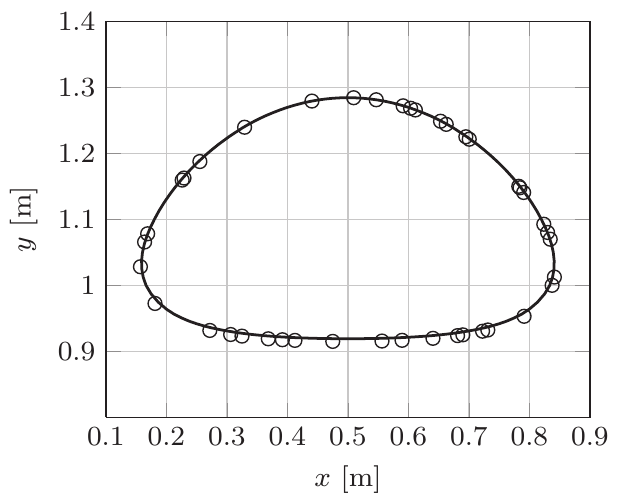}
}
\hfil
\subfigure[Rise velocity]{
\label{fig:rise-velocity}
\includegraphics[width=5cm]{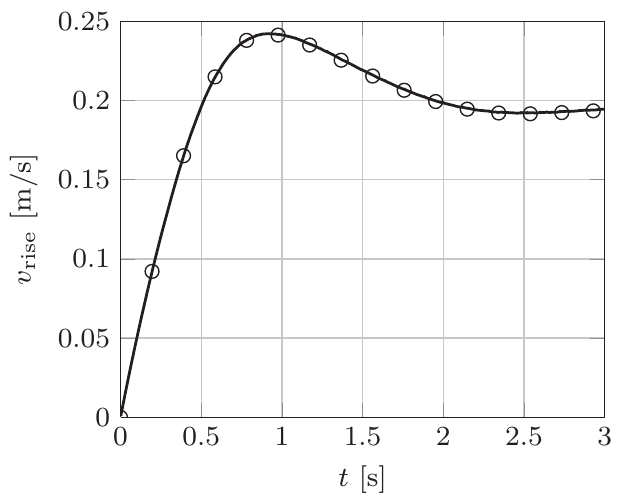}
}
\hfil
\subfigure{
\includegraphics[width=4.0cm]{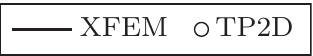}
}
\caption{Rising bubble: Comparison of the bubble shape and rise velocity obtained with the XFEM level-set approach in \cite{Sauerland2013} and simulation data from \cite{Hysing2009}.}
\label{fig:rising-drop-comparison}
\end{figure}

 The spatial resolution is $80 \times 160$ elements and $\Delta t = 0.002 s$. No-slip boundary conditions are assumed at the top and bottom boundary, slip boundary conditions are used along the vertical walls, and zero pressure is specified at the upper boundary. As an initial condition, the velocity field is set to ${\bf 0}$. The results \--- Figure~\ref{fig:bubble-shape} reports the bubble shape over time \---   give further support to the conclusions in \cite{Hysing2009}. In Figure~\ref{fig:bubble-shape}, the bubble shape at $t=3.0s$ is compared with the result from \cite{Hysing2009}, showing a very good agreement. In order to compare a time-evolving quantity, the rise velocity of the bubble, defined by the authors of the benchmark as:
 
\begin{align}
v_{rise} = \frac{\int_{\Omega_1} v({\bf x},t) d{\bf x}}{\int_{\Omega_1} 1 d{\bf x}}
\end{align} 

\noindent is depicted in Figure~\ref{fig:rise-velocity}.

Other publications where this or comparable benchmarks have been studied are, e.g., the following. In \cite{Mier2011}, a particle finite element method (PFEM) is employed to evaluate the second test case. In this method, the marker particles are treated in a Lagrangian manner. The interaction between the particles is portrayed through a connecting mesh, which is renewed after each displacement update. In addition, the article discusses the difficulties in modeling break-up and coalescence, which are always mesh-depended. Aland et al. consider both test cases with a phase field method \cite{Aland2012}. The governing equations for the mixture of incompressible, immiscible fluids lead to a coupling of the momentum equation \eqref{NavierStokes1} with a Cahn-Hillard type phase field equation. These are solved using an adaptive finite element scheme. Also with the phase-field method, rising bubble examples similar to \cite{Hysing2009} are investigated in \cite{Gruen2012}. A thermodynamically consistent set of governing equations is presented. The numerical scheme combines finite element discretizations for the momentum, phase-field, and chemical potential equations with a finite volume approach for all convective terms. Very large deformations as well as break-up and coalescence are inherent for these methods. In \cite{Klostermann2013} the authors were able to reproduce both benchmark cases (with and without breakup) with the finite volume method. Qualitatively, the results show the same features as in  \cite{Hysing2009}; the quantified results, however, present significant deviations from the ones obtained with the finite element codes. Doyeux et al. utilize the benchmark to validate their finite element method with interface capturing using level-set before they move on to the simulation of vesicles (a structure which serves as an imitation of the  mechanical behavior of red blood cells).  Again, both benchmark cases are presented \cite{Doyeux2012}. In \cite{Adelsberger2014}, the benchmark case is extended to 3D. The results of three different codes: XFEM with level-set, finite differences with level-set, and finite volumes with VOF are compared.  
 \section{Sloshing tank} \label{s-tank} 
 
A second common application in free-surface flows are sloshing tanks. The numerical simulation of sloshing tank problems is motivated by, e.g., the analysis of storage tanks under seismic loading, or partially filled tanks on ships such as fuel tanks and maritime transport of fluids or liquid gas. The aim of the simulation is to gain information about the forces acting on the surrounding structures, e.g., the tank walls or the ship, in order to investigate possible failure mechanisms. When considering such effects, the phenomenon of sloshing is in general not negligible when investigating the responses of the tank structure to the load due to the large mass of the liquid in combination with the low stiffness of the tank walls \cite{Ozdemir2009}.  If the problem is to be considered to its full extent, it falls into the category of fluid-structure interaction problems, where the thin structure of the tank wall reacts to the forces generated by the motion of the fluid and vice-versa. In this work, we will concentrate on the fluid component of such simulations. 

The amplitude of the sloshing depends on the amplitude and frequency of the seismic load as well as the fluid fill level, the fluid properties, and the tank geometry \cite{Ozdemir2009}. In the next section, we will focus on the last point, including arbitrary tank geometries. Usually, the considered tanks are either rectangular or cylindrical. With regard  to the fluid-structure-interaction problems, methods for arbitrarily shaped tanks need to be developed, however, in order to incorporate the different failure mechanisms of the tank wall into the sloshing phenomenon. Known shapes under failure of the tank structure are elephant footing, diamond-shaped buckling, as well as random convex and concave bumps \cite{Klinkel2012}. 

Within the overall context of the fluid-structure-interaction problem, concessions are often made with regard to the accuracy of the fluid flow model. Many authors decide to model the liquid in the tank as inviscid, incompressible and irrotational. This enables the use of a simple Laplace equation, which describes the velocity potential. If the full Navier-Stokes equations~\eqref{NavierStokes1}--\eqref{NavierStokes2} are chosen as governing equations, the sloshing tank is an application where the body forces are of great importance. In general, a gravitational force will be considered, since otherwise, any upward sloshing could never be reversed. In cases of rotating tanks, Coriolis forces are added to the general body forces. If seismic loading is to be considered, an alternating acceleration (usually orthogonal to gravity) is imposed. 

Another factor, which categorizes different scenarios, is the regularity of the flow conditions. While moderate deformations can be accurately and efficiently handled by interface tracking methods, wave breaking, air entrapment, and touching of the tank ceiling is more easily conducted with interface-capturing methods. Benchmark problems in this area are in general restricted to rectangular tanks with either an initially flat surface, which is subsequently excited (e.g., \cite{Huerta88a,Tezduyar92b}), or a free surface with a prescribed original deformation, which is then allowed to settle into equilibrium (e.g., in \cite{Fries2008,Koelke2005}). Even though this is a two-phase scenario, with a liquid phase in the container secluded by a gaseous phase on the top, in this application, the modelling is usually restricted to the liquid phase, as the density of the gaseous phase is several orders of magnitude lower \cite{Gomezgoni2013}. Exceptions need to be made, e.g., if the flow conditions are so violent that gas entrapments occurs. Ansari et al. \cite{Ansari2011} use the modal method, a reduced-order method, to model tank sloshing based on the flow potential (Laplace equation). The authors investigate the question under which conditions the top fluid can be neglected during the simulation. The results were later further analyzed in \cite{Gomezgoni2013}. In \cite{Gavrilyuk2013}, the scope of the modal method is extended from vertical walls to tapered walls. Numerical simulation and experimental validation of sloshing in a 2D rectangular tank is presented in \cite{Cruchaga2013}. The full Navier-Stokes equations are solved and the interface is captured using marker particles, which are updated in a Lagrangian manner. Subsequently, a global mass correction scheme is applied in order to ensure conservation of mass. Sloshing liquid motion in combination with a floating roof, which underlies the non-linear equation of motion, is the topic of \cite{Matsui2013}. The flow is modelled through the analytical solution modes of the velocity potential, and only the interface to the floating roof requires finite-element discretizations. In \cite{Ozdemir2009}, the authors decide to solve the compressible Navier-Stokes equations with an ALE finite element method.

In a similar application, Lee et al. \cite{Lee2012} simulate non-linear free-surface flows around ship bows including wave breaking phenomena. Due to the wave breaking, interface-capturing methods are more appropriate then interface tracking. In the referenced paper, a modified marker-density method is employed, circumventing the oversimplifications a traditional MAC-method would involve. As a solution method, finite differences with a sequential solution of velocity and pressure are used. One important challenge in this type of simulation is the conjunction point between the ship hull and the free surface, namely imposing the appropriate boundary conditions.

\subsection{The boundary conditions} 
 
\begin{figure}
\centering
\includegraphics[scale=0.9]{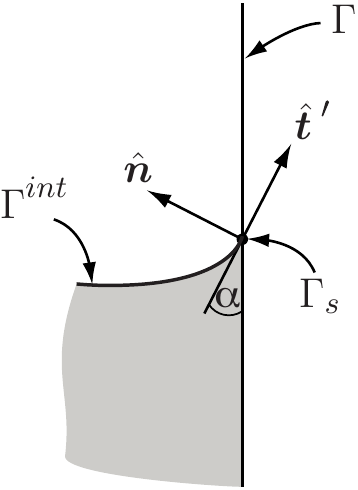}
\caption{Close-up of the contact line $\Gamma_s$ between an interface $\Gamma^{int}$ and a domain boundary $\Gamma$. The tangent vector of the interface at $\Gamma_s$, $\hat{{\bf n}}'$, and the boundary $\Gamma$ span the contact angle $\alpha$. If wetting models are not considered, $\alpha = 90^\circ$.}
\label{fig:contactline}
\end{figure}

The free surface requires a no-penetration boundary condition (cf. Equation~\eqref{eqn:fullLagrange}), indicating that the shape of the free surface is dictated through the fluid velocity. If one would like to consider surface tension effects in the pressure solution, Equation \eqref{eq:SurfTens} needs to be incorporated. Note that any technique relying on partial integration now needs to take the boundary integral into account (in contrast to the situation described in Section~\ref{s-challenges}), as the free surface in a sloshing tank is not closed. An exception can be made if we assume that the contact angle $\theta$, i.e., the angle between the tangent vector along $\Gamma^{int}$, $\hat{{\bf n}}'$, and the boundary $\Gamma$ is assumed to be constant at $90^\circ$ (cf. Figure~\ref{fig:contactline}). Then the boundary term is again $0$ as shown in \cite{Ganesan2009b}. This generalization is acceptable if wetting effects, as described in \cite{Kistler1993}, are assumed to not play an influential role. This holds true for the sloshing tank: As Behr points out in \cite{Behr2003a}, the capillary effects are negligible in the case of a sloshing tank, thus precluding the need for a specific wetting model. Instead, a global slip condition is found appropriate to allow the free surface to slide freely along the tank walls. The slip boundary condition can be expressed in the following way \cite{Behr2003a}:

\begin{align}
\label{eqn:slip-bc}
{\bf u} \cdot {\bf n} &= 0 \quad \mathrm{on} \ \Gamma_{slip} \nonumber \\
{\bf t}_1^T \bm{\sigma} ({\bf u},p) {\bf n}  &=0\quad \mathrm{on} \ \Gamma_{slip} \\
{\bf t}_2^T \bm{\sigma} ({\bf u},p)  {\bf n}  &=0\quad \mathrm{on} \ \Gamma_{slip} \,. \nonumber
\end{align}
 
As an alternative, a Navier slip condition \cite{Silliman78a} could be employed to incorporate at least some amount of wall friction. 

The tangential and normal vectors can either be computed analytically, if an appropriate description is available, or from the finite element mesh.

\subsection{Tanks with walls of arbitrary shape}

If interface-capturing methods are used, the shape of the tank wall is insignificant, as long as the computed velocity field complies with the slip boundary condition. Since the marker particles, the level-set field, or the VOF field are advected with the fluid velocity, they will obey the shape of the boundary. In an interface-tracking context, tanks with walls of arbitrary shape have so far posed a significant challenge. The reasons for this are two-fold. One reason lies in the mesh deformation, which needs to comply with the analytic boundary at all times. This means that in addition to the discretized boundary, sufficient information on the analytic boundary needs to be stored, e.g., in form of a CAD model. The second reason is due to the discrepancy between the discrete slip condition and the analytic slip condition. In \cite{Behr2003a}, it was found that this remains true even if a conforming normal vector is utilized. 

The general procedure for incorporating a slip boundary condition into an interface-tracking context is to first assemble the full system matrix ${\bf A}$, irrespective of any boundary conditions that need to be applied. Subsequently, all equations are rotated in such a way, that they are now formulated in a coordinate system spanned by the tangent and the normal vector of the respective boundary. Only then are the Dirichlet conditions, e.g., zero velocity in normal direction, applied. This again requires two steps: (1) At any finite element node (or degree of freedom) with a Dirichlet boundary condition, the corresponding weighting function is equal to zero. Therefore, the corresponding equation is deleted from the matrix  ${\bf A}$. (2) Furthermore, all occurrences of the specified degree of freedom are moved to the right-hand-side. Note that, usually, we encounter homogeneous Dirichlet conditions, where nothing needs to be done in this second step.
The described procedure is straightforward if the tangential coordinate system is aligned with the main coordinate axes. In all other cases, the tangential coordinate system together with the corresponding rotation matrix remains to be determined.

Consider the original linear system of equations with the nodal velocities in the $(x,y,z)$-coordinate system: 

\begin{align}
\underbrace{\begin{pmatrix}
{\bf A}_{11} & {\bf A}_{12}  & \dots \\
{\bf A}_{21} & {\bf A}_{22}  & \dots \\
\vdots & \vdots & \ddots & 
\end{pmatrix}}_{\bf A} \underbrace{\begin{pmatrix}
{\bf u}_{1} \\ {\bf u}_{2}  \\ \vdots
\end{pmatrix}}_{\bf u} = \underbrace{\begin{pmatrix}
{\bf b}_1 \\ {\bf b}_2 \\ \vdots
\end{pmatrix}}_{\bf b} \,.
\label{eqn:matrixsystem}
\end{align} 

For each node $i$ with degrees of freedom ${\bf u}_i = \begin{pmatrix}
u_{i,x} & u_{i,y} & u_{i,z},
\end{pmatrix}^T $ a rotation matrix $ {\bf O}_i $ is determined, which transforms from the $(x,y,z)$-coordinate system to a given tangential coordinate system. Rotation matrices are  square, orthogonal matrices with determinant 1 ( otherwise $ {\bf O}^{-1} = {\bf O}^{T} $ would not hold and the rotation could not be easily reversed). They can be obtained by placing the maps of the basis vectors of $\mathbb{R}^{nsd}$ into the columns of the rotation matrix. The local rotation matrix then takes the form \cite{Collins2004}:

\begin{align}
\label{eqn:rotationmatrix}
{\bf O}_i = \begin{pmatrix}
{\bf t}_1 & {\bf t}_2 & {\bf n} \,
\end{pmatrix}.
\end{align}

With the help of the rotation matrix, new degrees of freedom $ \tilde{{\bf u}}_i $ are defined in the tangential coordinate system:

\begin{align}
{\bf u}_i = {\bf O}_i \tilde{{\bf u}}_i \,.
\end{align}

Equation~\eqref{eqn:matrixsystem} may then be formulated in terms of the new degrees of freedom:

\begin{align}
\begin{pmatrix}
{\bf A}_{11} & {\bf A}_{12}  & \dots \\
{\bf A}_{21} & {\bf A}_{22}  & \dots \\
\vdots & \vdots & \ddots &
\end{pmatrix} 
\begin{pmatrix}
{\bf O}_1 & 0 &   \dots \\
0 & {\bf O}_2 &  \dots \\
 \vdots & \vdots  & \ddots 
\end{pmatrix}
\begin{pmatrix}
\tilde{{\bf u}}_1\\ \tilde{{\bf u}}_2 \\ \vdots
\end{pmatrix} = \begin{pmatrix}
{\bf b}_1 \\ {\bf b}_2  \\ \vdots
\end{pmatrix} \,.
\end{align}

\noindent With 

\begin{align}
\begin{pmatrix}
{\bf O}_1 & 0 &  \dots \\
0 & {\bf O}_2 & \dots \\
 \vdots & \vdots & \ddots 
\end{pmatrix}^{-1} = 
\begin{pmatrix}
{\bf O}_1^{-1} & 0 &  \dots \\
0 & {\bf O}_2^{-1} &  \dots \\
 \vdots & \vdots &  \ddots 
\end{pmatrix}  
= \begin{pmatrix}
{\bf O}_1^{T} & 0 & \dots \\
0 & {\bf O}_2^{T} &  \dots \\
 \vdots & \vdots &  \ddots 
\end{pmatrix},
\end{align}

\noindent we can solve the system using:

\begin{align}
\begin{pmatrix}
{\bf O}_1^{T} & 0 &   \dots \\
0 & {\bf O}_2^{T} &  \dots \\
 \vdots & \vdots &  \ddots 
\end{pmatrix}
{\bf A} 
\begin{pmatrix}
{\bf O}_1 & 0 & \dots \\
0 & {\bf O}_2 &  \dots \\
 \vdots & \vdots  & \ddots 
\end{pmatrix}
\begin{pmatrix}
\tilde{{\bf u}}_1\\ \tilde{{\bf u}}_2 \\ \vdots
\end{pmatrix} =
\begin{pmatrix}
{\bf O}_1^{T} & 0  & \dots \\
0 & {\bf O}_2^{T}  & \dots \\
 \vdots & \vdots  & \ddots 
\end{pmatrix}
\begin{pmatrix}
{\bf b}_1 \\ {\bf b}_2  \\ \vdots
\end{pmatrix} \,.
\end{align}

For analytical domain shapes (straight or tilted lines, circles), \cite{Behr2003a} already shows how the rotation matrix ${\bf O}$ can be formed. In order to cope with arbitrary domain shapes, we propose a boundary description that is based on NURBS (cf. Section \ref{s-NURBS} for the introduction to NURBS and the notation). The use of NURBS as a boundary description has the following advantages: (1) NURBS can represent arbitrary analytical and free-form shapes, (2) they can be extracted from CAD-models, (3) tangent and normal vectors can be analytically computed at any point on the NURBS. Due to point (3), the alignment of the system of equations of the flow solution can be performed in the standard way using a rotation matrix $ {\bf O} $ obtained as in Equation~\eqref{eqn:rotationmatrix}.

The first tangent vector is obtained as the first derivative of the NURBS surface (or curve) with respect to the first local parameter $\theta$; a second tangent vector can be obtained by taking the derivative with respect to the second local parameter $\iota$:

\begin{align}
{\bf t}_{1,NURBS} = \frac{\partial {\bf S}(\theta, \iota)}{\partial \theta} \,; \quad {\bf t}_{2,NURBS} = \frac{\partial {\bf S}(\theta, \iota)}{\partial \iota} \,.
\end{align}

Note, however, that in physical space, these two tangent vectors are not orthogonal on each other. This only holds for the parameter space spanned by $\theta$ and $\iota$. Since the rotation matrix $ {\bf O} $ requires orthogonality in physical space, Gram-Schmidt orthogonalization is performed to obtain the two tangent vectors ${\bf t}_{1,n}$ and ${\bf t}_{2,n}$:

\begin{align}
{\bf t}_{1,n} &= \frac{{\bf t}_{1,NURBS}}{\parallel{\bf t}_{1,NURBS}\parallel} \,; \\ {\bf t}_{2,n} &= \frac{{\bf t}_{2,NURBS}-({\bf t}_{2,NURBS} \cdot {\bf t}_{1,n}) {\bf t}_{1,n}}{\parallel {\bf t}_{2,NURBS}-({\bf t}_{2,NURBS} \cdot {\bf t}_{1,n}) {\bf t}_{1,n} \parallel} \,.
\label{eqn:orthonormaltangents}
\end{align}

From the two tangent vectors, the normal vector can be computed via cross product. 

Applying a similar procedure to the mesh deformation scheme is slightly more complex. The aim is that the boundary nodes of the mesh slide along the arbitrarily shaped boundary in order to accommodate the deformation imposed through the free surface. The individual displacement for each node $\bm{\upsilon}$ is to be determined via the mesh deformation equation~\eqref{eq:sopt-emum}. However, as for the flow solution, we require a slip boundary condition. In addition to this condition, which restrict the nodal displacement to the local tangent direction, the requirement that the nodes remain on the NURBS at all times needs to be fulfilled. In total, the following steps are performed:

\begin{enumerate}
\item Rotate the mesh equation into the local tangential coordinate systems.
\item Drop the equation responsible for the normal component.
\item Compute the tangential displacement.
\item From the displacement in ${\bf x}$-coordinates, deduce a displacement on the NURBS in the local NURBS coordinates $\theta$ and $\iota$.
\end{enumerate}

The steps 1. through 3. are performed within the solution of the mesh deformation equation. Step 4. involves a transformation of the displacement given in tangential coordinates into the local NURBS coordinates. If the displacement is directly available in NURBS coordinates, it is ensured that all nodes remain on the NURBS. The following relation holds for the displacement from one iteration ($it$) to the next ($it+1$), where an iteration can be both a  linearization step (such as in the Newton-Raphson algorithm) or a time step:

\begin{align}
{\bf x}_{it+1} - {\bf x}_{it} = {\bf S} ({\bm{\theta}}_{it+1}) - {\bf S} ({\bm{\theta}}_{it}) \underbrace{\approx}_{\mathrm{linearization}} {\bf J} \Delta {\bm{\theta}} \,.
\label{eqn:NURBSdisplacement}
\end{align}

$ \Delta {\bm{\theta}} $ refers to the displacement in NURBS-coordinates and $ {\bf J} $ is the Jacobian of $ {\bf S} $ with respect to $ {\bm{\theta}} $. From Equation~\eqref{eqn:NURBSdisplacement},  $ \Delta \bm{\theta} $ can be computed.

\begin{figure}[h]
\centering
\psfrag{a}{${\bf n}$}
\psfrag{b}{${\bf t}_{1,n}$}
\psfrag{c}{$\bm{\upsilon}$}
\psfrag{d}{$\Delta \bm{\theta} $}
\includegraphics[scale=1.0]{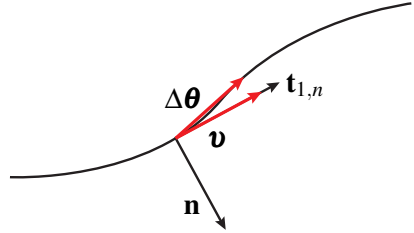}
\caption{The nodal displacement $ \bm{\upsilon} $ is computed in terms of the tangential coordinate system spanned by ${\bf t}_{1,n}, {\bf t}_{2,n}, {\bf n}$. This displacement is transformed into a displacement expressed in the NURBS parametric coordinates $\bm{\theta}$ based on a local linearization of the NURBS curve.} 
\label{NURBSdisplacement}
\end{figure}

\noindent Figure~\ref{NURBSdisplacement} illustrates the resulting displacement. 

In 2D, $ {\bf J} = \frac{\partial C}{\partial \theta}$ and $\Delta {\theta}$ can be readily obtained. In 3D, $ {\bf J}^{-1} $ cannot be determined easily and the (overdetermined) equation system needs to be solved. Starting point is step 3., where the tangential displacement $\Delta t_1, \Delta t_2$ has been obtained. The displacement reads

\begin{align}
\begin{pmatrix} {\bf t}_{1,n} &  {\bf t}_{2,n}  \end{pmatrix} \begin{pmatrix}
\Delta t_1 \\ \Delta t_2
\end{pmatrix} \,.
\end{align}

This displacement can be exactly represented as a displacement in the $ {\bf t}_{1} , {\bf t}_{2}$-coordinate system \--- i.e., the coordinate system aligned with the NURBS local coordinates. What remains is to determine the appropriate magnitude of the displacement, $\Delta \theta, \Delta \iota$, which fulfils the relation

\begin{align}
\begin{pmatrix} {\bf t}_{1,n} &  {\bf t}_{2,n}  \end{pmatrix} \begin{pmatrix}
\Delta t_1 \\ \Delta t_2
\end{pmatrix} = \begin{pmatrix} {\bf t}_{1} &  {\bf t}_{2}  \end{pmatrix} \begin{pmatrix}
\Delta \theta \\ \Delta \iota
\end{pmatrix} \,.
\end{align}

Note that $ \begin{pmatrix} {\bf t}_{1} &  {\bf t}_{2}  \end{pmatrix} $ is also the Jacobian matrix ${\bf J}$. The above relation can be rewritten as

\begin{align}
\Delta t_1 {\bf t}_{1,n} + \Delta t_2 {\bf t}_{2,n} = \Delta \theta {\bf t}_{1} + \Delta \iota {\bf t}_{2}  \,.
\label{eqn:equalityofdisplacement}
\end{align}

With the aim of comparing coefficients, the expressions \eqref{eqn:orthonormaltangents} are inserted in Equation~\eqref{eqn:equalityofdisplacement} 

\begin{align}
\Delta t_1 \frac{{\bf t}_1}{\parallel{\bf t}_1\parallel} + \Delta t_2 \frac{{\bf t}_{2}-({\bf t}_{2} \cdot {\bf t}_{1,n}) {\bf t}_{1,n}}{\parallel {\bf t}_{2}-({\bf t}_{2} \cdot {\bf t}_{1,n}) {\bf t}_{1,n} \parallel} = \Delta \theta {\bf t}_{1} + \Delta \iota {\bf t}_{2}  \,.
\end{align}

\noindent leading to the final result of

\begin{align}
\Delta \theta &=  \frac{1}{\parallel{\bf t}_1\parallel} \left( \Delta t_1 - \frac{\Delta t_2 ({\bf t}_{2} \cdot {\bf t}_{1,n}) }{\parallel {\bf t}_{2}-({\bf t}_{2} \cdot {\bf t}_{1,n}) {\bf t}_{1,n} \parallel} \right) \\
\Delta \iota &= \frac{\Delta t_2 }{\parallel {\bf t}_{2}-({\bf t}_{2} \cdot {\bf t}_{1,n}) {\bf t}_{1,n} \parallel}
\end{align}

In such cases, one final remark needs to be made about the intersection between the free surface and the wall ($\Gamma_s$ in Figure~\ref{fig:contactline}). Here, it is important that the displacement of the free surface is chosen compatible with the displacement along the wall.

The effect of the new implementation of the slip condition is detailed by means of a numerical example. 

\begin{figure}[Htbp]
\center
\subfigure[$t=4.0 s$]{
\includegraphics[width=10.0cm]{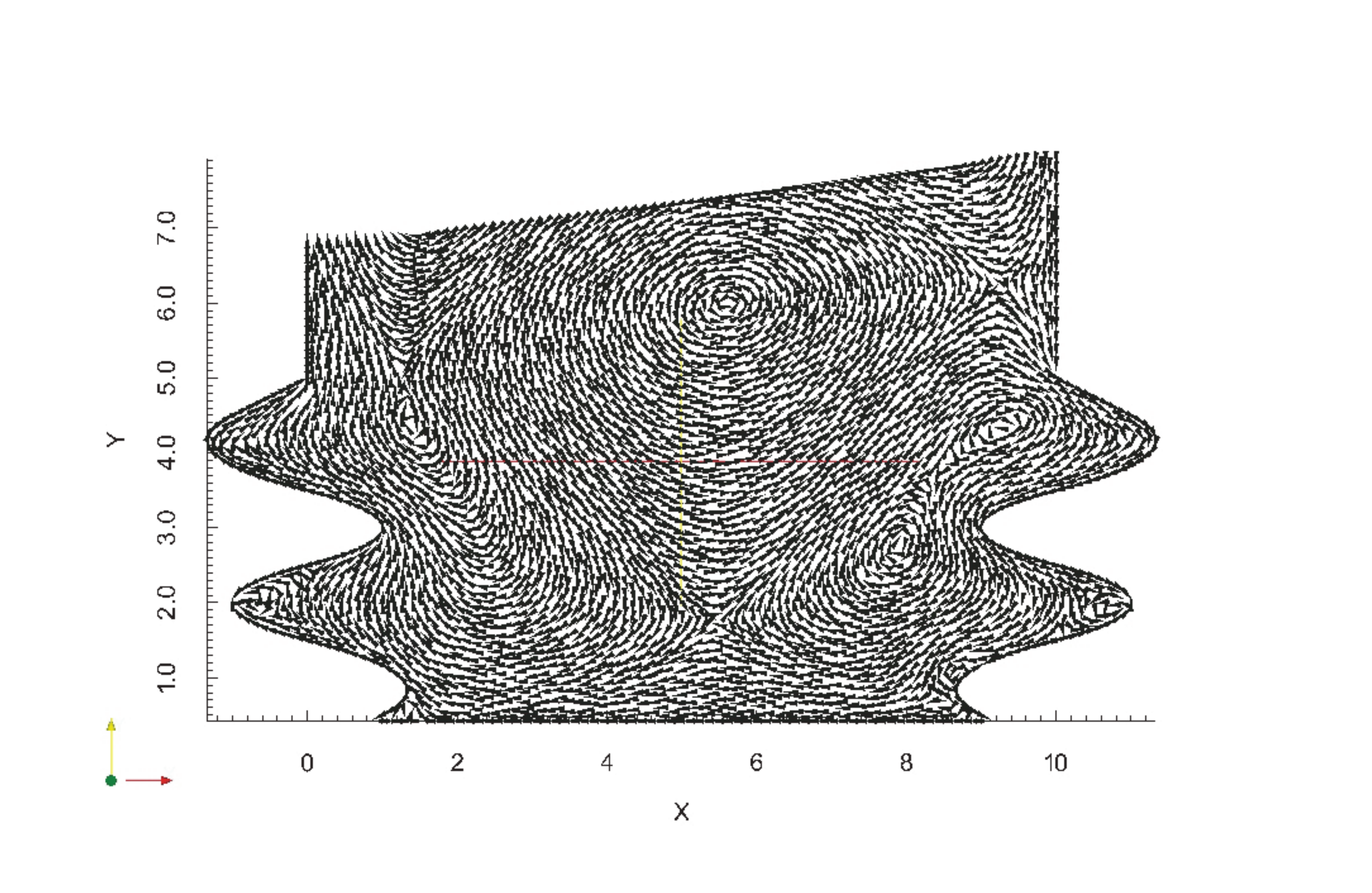}
}
\hfil
\subfigure[$t=13.5 s$]{
\label{resistance}
\includegraphics[width=10cm]{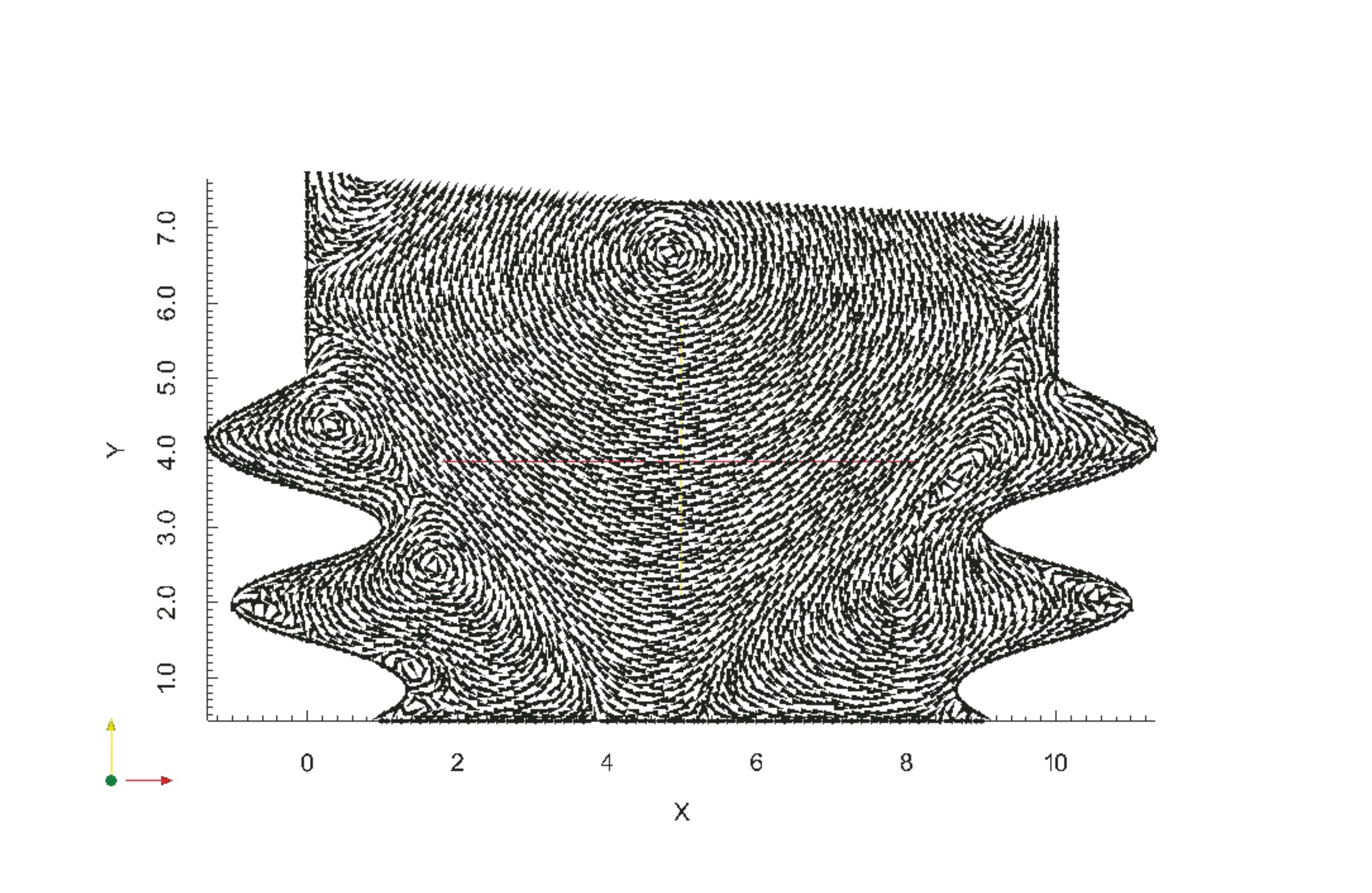}
}
\caption{Sloshing in a tank with curved side walls: the tank wall is described by a quadratic NURBS curve with knot vector $[0\ 0\ 0\ 1/16\  2/16\  3/16\ 4/16\ 5/16\  6/16\ $ $7/16\ 8/16\ 9/16\ 10/16\ 11/16\ 12/16\ 13/16\ 14/16\ 15/16\ 1\ 1\ 1]$ and control points $(10,10); (10,5); (10,5); (12,4); (8,3); (12,2); (8,1); (10,0); (10,0); (0,0); (0,0); (2,1);$ $ (-2,2); (2,3); (-2,4); (0,5); (0,5); (0,10) $. The velocity vectors of the flow solution are depicted after $4 s$ and $13.5 s$. Note how the velocity aligns with the boundary. \cite{Sliwiak2014}}
\label{fig:tank2D}
\end{figure}

In the test case, we consider an undeformable tank with curved side walls as indicated in Figure~\ref{fig:tank2D}. The tank wall is defined through a quadratic NURBS curve with 18 control points. Along all walls, a slip boundary condition is imposed. The top boundary constitutes a free-surface. The tank is subjected to a downward gravitational acceleration $g = -0,98 m/s^2$ and a horizontal sinusoidal acceleration $s = -0.2 \sin(t)$. The considered fluid has a density of $\rho = 1000 kg/m^3$ and a viscosity of $\mu = 0.1 kg/m/s$. 
Figure~\ref{fig:tank2D} depicts the velocity vectors of the flow solution after $4 s$ and $13.5 s$. Note how the velocity aligns with the boundary. The mesh deformation is conform with the boundary at all times. 
 \section{Conclusion} \label{s-conclusion} 
 
This paper gave an overview of the challenges in the simulation of fluid flow in deformable domains. Five main methods have been introduced with their individual advantages and drawbacks. As was seen, some may be used more often than others, but each has its niche. A recent area for the Marker and Cell method has been found in flow simulations for movies, where a qualitatively correct solution is more important than quantitative accuracy. Level-set has its main advantage when topological changes are to be expected within the domain. The volume-of-fluid method is utilized frequently in commercial codes. In cases where boundary conditions are difficult to impose, the use of the phase-field method can be worthwhile.  The interface tracking approach may have a limited scope of applications, but comes with very high computational efficiency. 

Out of a very large number of possibilities, the applications of drops and sloshing tanks were selected to demonstrate some of the approaches when computing free-boundary problems.

\newpage

\subsection*{Acknowledgements}
The authors gratefully acknowledge support from the German Research Foundation (DFG) through the SFB 1120 ``Thermal Precision'', the Emmy-Noether-research group ``Numerical methods for discontinuities in continuum mechanics'', and the DFG program GSC 111 (AICES Graduate School). The computations were conducted on computing clusters provided by the J\"ulich Aachen Research Alliance (JARA). Furthermore, we would like to thank Philipp Knechtges for the infinite patience with which he shared his mathematical insight.

\bibliography{shorttitles,referencesARCME}

\end{document}